%% file: tropical-weierstrass-points.tex
\documentclass[11pt,leqno]{amsart}

\usepackage{amsmath}
\usepackage{amssymb}
\usepackage{amscd}
\usepackage{latexsym}
\usepackage{mathtools}
\usepackage{mathrsfs}
\usepackage{mathdots}
\usepackage{tabularx}
\usepackage{blkarray}
\usepackage{a4wide}
\usepackage[usenames]{color}
\usepackage{epsfig}
\usepackage{enumerate}
\usepackage{subfigure}
\usepackage{url}
\usepackage{tikz,tikz-cd}
\usetikzlibrary{decorations.markings}
\usetikzlibrary{calc}
\usepackage[matrix,arrow,tips,curve]{xy}
\usepackage{pb-diagram,pb-xy}
\usepackage{verbatim}
\usepackage{scalerel}
\usepackage{mleftright}

\definecolor{cadmiumgreen}{rgb}{0.0, 0.42, 0.24}
\usepackage[
colorlinks, citecolor=cadmiumgreen,
pdfauthor={Omid Amini, Lucas Gierczak, Harry Richman}, 
pdftitle={Tropical Weierstrass points},
pdfstartview ={FitV},
bookmarksdepth=3
]{hyperref}

\usepackage[
msc-links,
nobysame,
lite,
alphabetic,
]{amsrefs} 

\usepackage{todonotes} 

\newtheorem{thm}{Theorem}[section]

\newtheorem{lemma}[thm]{Lemma}

\newtheorem{prop}[thm]{Proposition}

\newtheorem{cor}[thm]{Corollary}

\newtheorem{question}[thm]{Question}

\theoremstyle{definition}
\newtheorem{defii}[thm]{Definition}
\newenvironment{defi}
	{\pushQED{\qed}\defii}
	{\popQED\enddefii}
\newtheorem{remm}[thm]{Remark}
\newenvironment{remark}
	{\pushQED{\qed}\remm}
	{\popQED\endremm}
\newtheorem{exx}[thm]{Example}
\newenvironment{example}
	{\pushQED{\qed}\exx}
	{\popQED\endexx}

\newcommand{\cE}{\mathcal E}

\newcommand{\cO}{\ssub{\mathcal O}!}
\newcommand{\cD}{\mathcal D}
\newcommand{\cW}{\mathcal W}
\newcommand{\cL}{\ssub{\mathcal L}!}
\newcommand{\cA}{\mathcal A}
\newcommand{\cC}{\mathcal C}
\newcommand{\cK}{\mathcal K}

\newcommand{\cF}{\mathcal F}
\newcommand{\Wron}{\ssub[-1pt]{\mathrm{Wr}}!}

\newcommand{\ssunif}{\ssub{{\mathfrak t}}!}
\newcommand{\ssg}{\ssub{g}!}
\newcommand{\ssf}{\ssub{f}!}

\newcommand{\Q}{\mathbb{Q}}
\newcommand{\QQ}{\mathbb{Q}}
\newcommand{\Z}{\mathbb{Z}}

\newcommand{\R}{\mathbb{R}}
\newcommand{\RR}{\mathbb{R}}

\newcommand{\K}{\mathbb K}
\newcommand{\FF}{\mathbb{F}}

\newcommand{\fC}{\mathfrak{C}}
\newcommand{\fD}{\mathfrak D}
\newcommand{\fK}{\mathfrak K}
\newcommand{\fX}{\ssub{\mathfrak X}}
\newcommand{\fS}{\mathfrak S}
\newcommand{\fA}{\mathfrak A}
\newcommand{\g}{\mathfrak g}
\newcommand{\f}{\mathfrak f}

\newcommand{\T}{\mathrm{T}}
\renewcommand{\k}{\kappa}

\newcommand{\ord}{\mathrm{ord}}

\newcommand{\st}{\bigm|}

\DeclareMathOperator{\Div}{Div}

\DeclareMathOperator{\Pic}{Pic}
\DeclareMathOperator{\Eff}{Eff}
\DeclareMathOperator{\Spec}{Spec}

\DeclareMathOperator{\an}{an}

\renewcommand{\div}{\mathrm{div}}

\newcommand{\W}{{\scaleto{W}{4pt}}}
\newcommand{\Wloc}{L_\W}

\newcommand{\wfinite}{\text{W-finite}}
\newcommand{\slope}{\mathrm{sl}}
\newcommand{\weight}{\mu_\W}
\newcommand{\hatweight}{\hat\mu_\W}
\DeclareMathOperator\val{val}
\DeclareMathOperator{\outval}{outval}
\DeclareMathOperator{\Rat}{Rat}
\DeclareMathOperator{\KRat}{KRat}
\newcommand{\PLZ}{\Rat}

\newcommand{\tropic}[1]{\mathrm{trop}\mleft(#1\mright)}

\newcommand{\degrest}[2]{\deg\mleft({#1} \rest{#2} \mright)}
\newcommand{\rest}[1]{\raisebox{-1pt}{$\vert$}_{#1}}
\newcommand{\abs}[1]{\mleft|#1\mright|}

\newcommand{\partialout}{\partial^{^{\scaleto{\mathrm{out}}{3pt}}}}
\newcommand{\partialin}{\partial^{^{\scaleto{\mathrm{in}}{3pt}}}}

\newcommand{\generic}{{\hspace{-.02cm}\scaleto{\mathrm{gen}}{3.3pt}}} 

\newcommand{\alg}{{^{\hspace{-.02cm}\mathrm{clls}}}}

\newcommand{\Ratgen}{\Rat^{^\generic}}
\newcommand{\Wlocgen}{L_\W^{^\generic}}
\newcommand{\weightgen}{\weight^\generic}

\makeatletter
\newsavebox\myboxA
\newsavebox\myboxB
\newlength\mylenA

\newcommand*\overbar[2][0.75]{%
	\sbox{\myboxA}{$\m@th#2$}%
	\setbox\myboxB\null
	\ht\myboxB=\ht\myboxA%
	\dp\myboxB=\dp\myboxA%
	\wd\myboxB=#1\wd\myboxA
	\sbox\myboxB{$\m@th\overline{\copy\myboxB}$}
	\setlength\mylenA{\the\wd\myboxA}
	\addtolength\mylenA{-\the\wd\myboxB}%
	\ifdim\wd\myboxB<\wd\myboxA%
		\rlap{\hskip 1\mylenA\usebox\myboxB}{\usebox\myboxA}%
	\else
		\hskip -0.5\mylenA\rlap{\usebox\myboxA}{\hskip 0.5\mylenA\usebox\myboxB}%
	\fi}

\renewcommand\part{%
	\vspace*{2.5 ex} 
	\if@noskipsec \leavevmode \fi
	\@afterindentfalse
 \secdef \@part\@spart}

\def\@part[#1]#2{%
	\ifnum \c@secnumdepth >\m@ne
		\refstepcounter{part}%
		\addcontentsline{toc}{part}{\partname \nobreakspace \thepart:\hspace{1em}#1} 
	\else
		\addcontentsline{toc}{part}{#1}%
	\fi 
	{\parindent \z@ \centering 
	 \interlinepenalty \@M
	 \normalfont
	 \ifnum \c@secnumdepth >\m@ne
	\LARGE \partname \nobreakspace\thepart:
		\fi 
		\LARGE {\Fontskrivan{#2}}%
	 \par}%
	\nobreak
	\vskip 5ex
	\@afterheading}
\def\@spart#1{%
	{\parindent \z@ \centering 
	 \interlinepenalty \@M
	 \normalfont
	 \huge \bfseries #1\par}%
	 \nobreak
	 \vskip 5ex
	 \@afterheading}

\makeatother

\renewcommand{\thepart}{\Roman{part}} 

\newcommand{\comp}[1]{\overbar[.5]{#1}} 

\newcommand{\Trop}{\mathrm{trop}} 
\newcommand{\trop}{{^{\hspace{-.02cm}\scaleto{\Trop}{4.4pt}}}} 

\newcommand{\mgbar}{\comp{\mathscr M}} 
\newcommand{\mgbarg}[1]{\mgbar_{{\hspace{-.09cm}#1}}^{ }} 
\newcommand{\unicurve}[1]{\mathscr{C}_{{{\hspace{-.05cm}#1}}}^{}}

\NewDocumentCommand{\ssub}{O{0pt} O{.8} m t! e{_^}}{
	#3%
	\IfValueT{#5}{
		\IfBooleanTF{#4}{\sb{\hspace{#1}\scaleobj{#2}{#5}}}{\sb{#5}}
	}
	\IfValueT{#6}{
	\IfBooleanTF{#4}{\sp{\hspace{#1}\scaleobj{#2}{#6}}}{\sp{#6}}
}
}

\ExplSyntaxOn
\NewDocumentCommand{\tossub}{o o m}{
	\expandafter\let\csname old\cs_to_str:N #3\endcsname#3
	\renewcommand#3%
	{\ssub[#1][#2]{\csname old\cs_to_str:N #3\endcsname}}
}
\ExplSyntaxOff

\tossub\omega
\tossub\Omega

\newcommand{\valuation}{\mathfrak v}
\newcommand{\vR}{\mathrm{R}}
\newcommand{\fm}{\mathfrak{m}}
\newcommand{\rmC}{\mathrm{C}}
\newcommand{\sstildef}{\ssub{\tilde f}!}

\newcommand{\sstildealpha}{\ssub{\tilde \alpha}!}
\newcommand{\sstildea}{\ssub{\tilde a}!}

\newcommand{\sstildeH}{\ssub{\widetilde H}!}
\newcommand{\ssa}{\ssub{a}!}
\newcommand{\ssp}{\ssub{p}!}

\newcommand{\B}{\mathbb B}
\newcommand{\ssA}{\ssub{\mathbb A}!}

\newcommand{\varT}{\ssub{\scaleto{{T}}{6.5pt}}!}
\newcommand{\vart}{t}
\renewcommand{\setminus}{\smallsetminus}

\newcommand{\Der}[1]{\mathrm{D}^{^{(#1)}}}

\title{Tropical Weierstrass points and Weierstrass weights}

\author{Omid Amini}
\address{CNRS - Laboratoire de Mathématiques d'Orsay, Université Paris-Saclay}
\email{\href{omid.amini@universite-paris-saclay.fr}{omid.amini@universite-paris-saclay.fr}}

\author{Lucas Gierczak}
\address{Institut de mathématiques de Marseille, Universit\'e d'Aix-Marseille, CNRS, Centrale Marseille}
\email{\href{lucas.gierczak-galle@univ-amu.fr}{lucas.gierczak-galle@univ-amu.fr}}

\author{Harry Richman}
\address{Department of Mathematics, University of Washington --- Mathematics Division, National Center for Theoretical Sciences}
\email{\href{hrichman@ncts.ntu.edu.tw}{hrichman@ncts.ntu.edu.tw}}

\subjclass[2020]{
\href{https://mathscinet.ams.org/msc/msc2020.html?t=14H55}{14H55},
\href{https://mathscinet.ams.org/msc/msc2020.html?t=14T20}{14T20},
\href{https://mathscinet.ams.org/msc/msc2020.html?t=05C22}{05C22},
\href{https://mathscinet.ams.org/msc/msc2020.html?t=14H10}{14H10},
\href{https://mathscinet.ams.org/msc/msc2020.html?t=14H20}{14H20},
\href{https://mathscinet.ams.org/msc/msc2020.html?t=14T15}{14T15}}

\begin{document}

\begin{abstract}
	In this paper, we study tropical Weierstrass points. These are the analogues for tropical curves of ramification points of line bundles on algebraic curves.

	For a divisor on a tropical curve, we associate intrinsic weights to the connected components of the locus of tropical Weierstrass points. These are obtained by analyzing the slopes of rational functions in the complete linear series of the divisor. We prove that for a divisor $D$ of degree $d$ and rank $r$ on a genus $g$ tropical curve, the sum of weights is equal to $d - r + rg$. We establish analogous statements for tropical linear series.
	
	In the case $D$ comes from the tropicalization of a divisor, these weights control the number of Weierstrass points that are tropicalized to each component.
	Our results provide answers to open questions originating from the work of Baker on specialization of divisors from curves to graphs.
	
	We conclude with multiple examples that illustrate interesting features appearing in the study of tropical Weierstrass points, and raise several open questions.
\end{abstract}

\maketitle
\setcounter{tocdepth}{1}

\tableofcontents

\section{Overview} \label{sec:overview}
	
	Weierstrass points have a rich history in the development of algebraic geometry as they provide an important tool for the study of smooth algebraic curves and their moduli spaces. It is natural to ask how their theory can be extended to stable curves, which correspond to boundary points in the Deligne--Mumford compactification $\mgbarg g$ of the moduli space of genus $g$ smooth curves. One strategy is to take the {\em limit Weierstrass points} induced by a one-parameter family $(X_t)_{t \neq 0}$ of smooth curves degenerating to a stable curve $X_0$; there will be $g^3 - g$ limit Weierstrass points on $X_0$ when counted with appropriate weights. However, the limit points generally depend on the chosen family, and a stable curve $X_0$ has many possible smoothings corresponding to paths in $\mgbarg g$ that end at the point representing $X_0$.
	
	Tropical geometry provides a new perspective on degeneration methods in algebraic geometry by enriching it with polyhedral geometry. Given the successes of tropical methods in the past two decades in the study of algebraic curves and their moduli spaces, it is natural to ask whether tropical geometry can be used to gain insight about the limiting behavior of Weierstrass points on degenerating families of curves. 
	In the tropical perspective, the data of a stable curve $X_0$ is replaced by the data of its dual graph, endowed with edge lengths. The collection of all stable curves having the same dual graph forms a stratum of $\mgbarg{g}$. This gives a correspondence between the strata of $\mgbarg{g}$ and the set of stable graphs of genus $g$~\cite{caporaso}.
	
	The prototype of what we can expect to address using tropical techniques is the following natural question. 
	
	\begin{question} \label{question:clls_vs_tropical}
		Given a stratum of $\mgbarg{g}$, what can be said about the limit Weierstrass points of a smooth family $(X_t)_{t \neq 0}$ degenerating to a stable curve in that stratum?
	\end{question}
	
	The arithmetic geometric version of the above question is formulated as follows: Given a smooth proper curve over the field $\Q_p$ of $p$-adic numbers with stable reduction lying in a given stratum of $\mgbarg{g}$ (over the algebraic closure of the residue field $\FF_p$), what can be said about the specialization of the Weierstrass points?
	
	
	Previously, there has been much work making incremental progress on the above question, e.g.~\cite{EH87-WP, EM02, ES07, Dia85, BL12, Ami14, Gen21} for the geometric setting, and~\cite{Ogg78, LN64, Atk67, AP03, baker2008specialization} for the arithmetic setting. We refer to Section~\ref{subsec:previous_work} for a more thorough discussion.
	
	Our aim in this paper is to provide an answer to the above questions from the point of view of tropical geometry. This is done by introducing new tools that allow us to solve problems related to the tropical geometry of curves, whose origin goes back to the beginning of the use of tropical methods in the study of algebraic curves.
	
	Our answer to Question~\ref{question:clls_vs_tropical} can be summarized as follows: we can specify how many Weierstrass points degenerate to each component and to each node of a stable curve $X_0$ lying in the given stratum. This is done without specifying their precise position within each irreducible component, giving instead a more precise location of those degenerating to a node by specifying their position on the dual metric graph of the family $(X_t)$. Our result also applies to limits of ramification points of arbitrary line bundles, in addition to the case of the canonical bundle.
	
	Similarly, we answer the arithmetic geometric version of Question~\ref{question:clls_vs_tropical} by specifying where Weierstrass points specialize when reducing modulo $p$. Moreover, these results lead to an effective way of locating the limit Weierstrass points.
	
	\smallskip
	
	We next give an overview of our results.
	
	\subsection{Tropical perspective}
		
		The central concept studied in this paper is that of tropical Weierstrass points. The definition is based on divisor theory on metric graphs, and we refer to the survey paper~\cite{BJ16} and the references there for more details.
		
		Let $\Gamma$ be a metric graph, and let $D$ be a divisor of degree $d$ on $\Gamma$. We denote by $r$ the \emph{rank} of $D$, defined as the maximum integer $r \geq -1$ such that for any effective divisor $E$ of degree $r$, there exists an element $f \in M$ with $\div(f) + D - E \geq 0$. In the present paper, we consider only divisors of nonnegative rank since, otherwise, the Weierstrass locus will be empty.
		
		\begin{defi}[Weierstrass points] \label{def:weierstrass_points_intro}
			A point $x$ in $\Gamma$ is called a {\em Weierstrass point}, or {\em ramification point}, for $D$ if there exists an effective divisor $E$ in the linear system of $D$ whose coefficient at $x$ is at least $r + 1$. The {\em (tropical) Weierstrass locus} of $D$, denoted by $\Wloc(D)$, is the set of all such points in $\Gamma$. Clearly, it is preserved by linear equivalence of divisors.
		\end{defi}
		
		The set $\Wloc(D)$ is a closed subset of $\Gamma$ that can be infinite, in contrast with the classical setting of algebraic curves. In this regard, Baker comments in~\cite[Rem.~4.14]{baker2008specialization}, regarding the case when $D$ is canonical divisor, that ``it is not clear if there is an analogue for metric graphs of the classical fact that the total weight of all Weierstrass points on a smooth curve of genus $g$ is $g^3 - g$.'' More generally, we can ask the following question.
		
		\begin{question} \label{question:intrinsinc_weight}
			Is it possible to associate intrinsic tropical weights to the connected components of $\Wloc(D)$? What is the total sum of weights associated to these components?
		\end{question}
		
		Note that this question is interesting even when the locus of Weierstrass points is finite.
		
		Our aim in this paper is to provide an answer to the above question. In order to streamline the presentation that follows, we first discuss our results in the case of non-augmented metric graphs. From the geometric perspective, this corresponds to the situation of a {\em totally degenerate} stable curve, that is, a stable curve whose irreducible components are all projective lines.
		This is the same as requiring that the arithmetic genus of the stable curve is equal to the genus of the dual graph.
		We have an analogue of these statements for augmented metric graphs (respectively, arbitrary stable curves), see the discussion that follows below.
		
		In order to solve Question~\ref{question:intrinsinc_weight}, we make the following definition.
		
		\begin{defi}[Intrinsic Weierstrass weight of a connected component] \label{def:w_weight_component_intro}
			Let $D$ be a divisor of rank $r$, and let $C$ be a connected component of the Weierstrass locus $\Wloc(D)$.
			For a rational function $f$ on $\Gamma$, a point $x \in \Gamma$ and a unit tangent direction $\nu \in \T_x(\Gamma)$, we denote by $\slope_\nu f(x)$ the slope of $f$ along the tangent direction $\nu$ based at $x$.
			We define the {\em tropical Weierstrass weight} of $C$ as
			\begin{equation} \label{eq:w_weight_component}
				\weight(C; D) \coloneqq \degrest{D}{C} + \mleft(g(C) - 1\mright) \, r - \sum_{\nu \in \partialout C} s_0^\nu(D)
			\end{equation}
			where
			\begin{itemize}
				\item $\degrest{D}{C}$ is the total degree of $D$ in $C$, defined by $\degrest{D}{C} = \sum_{x \in C} D(x)$;
				
				\item $g(C)$ is the genus of $C$, i.e., its first Betti number $\dim H_1(C, \RR)$; 
				
				\item $\partialout C$ is the set of outgoing unit tangent directions from $C$; and
				
				\item $s_0^\nu(D)$ is the minimum slope $\slope_\nu f(x)$ along tangent direction $\nu$ of any rational function $f$ on $\Gamma$ with $\div(f) + D \geq 0$.
	
			\end{itemize}
			We abbreviate $\weight(C; D)$ simply as $\weight(C)$ if $D$ is understood from the context.
		\end{defi}
		
		Although it is not obvious from the definition, we will show in Theorem~\ref{thm:positivity} that the tropical Weierstrass weight of any component is positive. Note as well that a connected component of $\Wloc(D)$ is always a metric subgraph of $\Gamma$, see Proposition~\ref{prop:well_defined}. Moreover, the tropical Weierstrass weights are independent of the choice of the divisor $D$ in its linear equivalence class, see Proposition~\ref{prop:tropical_weight_independent_linear_equivalence}.
		
		We say that $D$ is {\em Weierstrass finite} or simply {\em \wfinite{}} if the tropical Weierstrass locus $\Wloc(D)$ has finite cardinality. In this case, connected components of $\Wloc(D)$ are isolated points in $\Gamma$, and we define the {\em tropical Weierstrass divisor} $W(D)$ as the effective divisor
		\[
			W(D) \coloneqq \sum_{x \in \Wloc(D)} \weight(x) \, (x)
		\]
		where $\weight(x) \coloneqq \weight(\{x\})$. The support $|W(D)|$ of the tropical Weierstrass divisor is exactly the tropical Weierstrass locus $\Wloc(D)$.
		The tropical Weierstrass weight of $x$ can be identified as $\weight(x) = D_x(x) - r$, with $D_x$ denoting the unique $x$-reduced divisor in the linear system of $D$, see Remark~\ref{rem:multiplicity_singleton}.
		
		This gives the following geometric meaning to the Weierstrass weights, in the spirit of the classical definition on algebraic curves. The coefficient of the reduced divisor at a point $x \in \Gamma$ corresponds precisely to the maximum order of vanishing at $x$ of any global section of the line bundle $\mathcal O(D)$ defined by the divisor. The Weierstrass weight of the point $x$ is thus obtained by comparing this quantity to $r$, which would be the expected minimum value, over points $y \in \Gamma$, of the largest order of vanishing of global sections at $y$. (Note, however, that $r$ is not always equal to the actual minimum largest order of vanishing, as examples in Section~\ref{subsec:graph_all_weierstrass} show.) That being said, the definition differs from the algebraic setting, where we need to take into account {\em all} the orders of vanishing of global sections of the line bundle at a given point (and then compare them with the standard sequence, the one obtained for a point in general position on the curve).
		
		\smallskip
		
		The following theorem answers Question~\ref{question:intrinsinc_weight}. It is proved in Section~\ref{subsec:consequences_theorem}.
		
		\begin{thm}[Total weight of the Weierstrass locus] \label{thm:number_w_points_intro}
			Let $\Gamma$ be a metric graph of genus $g$, and let $D$ be an effective divisor of degree $d$ and rank $r$ on $\Gamma$. Then, the total sum of weights of the connected components of $\Wloc(D)$ is equal to $d - r + rg$. In particular, if $D$ is \wfinite{}, then we have $\deg(W(D)) = d - r + rg$.
		\end{thm}
		
		The proof of this theorem will imply in particular the following local topological constraint on the locus of tropical Weierstrass points. The proof is given in Section~\ref{subsec:consequences_theorem}.
		
		\begin{thm} \label{thm:cycle_weierstrass_point}
			If the rank $r$ of $D$ is at least one, then every cycle in $\Gamma$ intersects the tropical Weierstrass locus $\Wloc(D)$. In particular, if $\Gamma$ has genus at least two, then every cycle intersects the Weierstrass locus of the canonical divisor $K$.
		\end{thm}
		
		In~\cite{baker2008specialization}, Baker proves that the tropical Weierstrass locus of the canonical divisor is nonempty if $\Gamma$ has genus at least two. This earlier tropical result is obtained as a consequence of the analogous algebraic statement, using the specialization lemma. In contrast, our theorem above states that tropical Weierstrass points obey a stronger ``local'' existence condition, which has seemingly no algebraic analogue. In the case that the canonical divisor of $\Gamma$ is \wfinite{}, our result implies that for an arbitrary family $(X_t)_{t \neq 0}$ of smooth curves tropicalizing to $\Gamma$, every cycle in $\Gamma$ contains a limit Weierstrass point of the family.
		
		To prove Theorem~\ref{thm:number_w_points_intro}, we will show that in fact~\eqref{eq:w_weight_component} defines a consistent notion of Weierstrass weight when applied to any connected, closed subset of $\Gamma$ whose boundary points are not in the interior of $\Wloc(D)$; see Section~\ref{subsec:W_measure} and Theorem~\ref{thm:weight_measure} for the precise statement.
	
	\subsection{Comparison results and extensions}
		
		We further justify our definition of weights by making a precise link to tropicalizations of Weierstrass points on algebraic curves.
		
		Suppose that $\Gamma$ and $D$ come from geometry; that is, let $X$ be a smooth proper curve of genus $g$ over an algebraically closed non-Archimedean field $\K$ of characteristic zero with a non-trivial valuation and a residue field of arbitrary characteristic. Let $\cL = \cO(\cD)$ be a line bundle of degree $d$ on $X$. Assume that $\Gamma$ is a skeleton of the Berkovich analytification $X^{\an}$ of $X$. Denote by $\tau$ the tropicalization map from $X$ to $\Gamma$, and suppose that $D = \tau_*\mleft(\cD\mright)$ is the tropicalization of $\cD$ on $\Gamma$ where $\tau_* \colon \Div(X) \to \Div(\Gamma)$ the induced map on divisors.
		
		Denote by $\cW(\cD)$ the Weierstrass divisor of $\cD$ on $X$, and by $\tau_*\mleft(\cW(\cD)\mright)$ its tropicalization on $\Gamma$. The following result is proved in Section~\ref{subsec:proofs_tropicalization_theorems}.
		
		\begin{thm}[Algebraic versus tropical Weierstrass weights] \label{thm:weierstrass_tropicalization_intro}
			Assume that $D$ and $\cD$ have the same rank $r$, and let $A \subseteq \Gamma$ be a closed, connected subset whose intersection with the Weierstrass locus of $D$ is a union of components of the Weierstrass locus. Then, the total weight of Weierstrass points of $\cW(\cD)$ tropicalizing to points in $A$ is precisely
			\[
				\degrest{\cW(\cD)}{\tau^{-1}(A)} = (r + 1) \mleft(\degrest{D}{A} + r\mleft(g(A) - 1\mright) - \sum_{\nu \in \partialout A} s_0^\nu(D) \mright).
			\]
			
			In particular, if, in addition, $D$ is \wfinite{}, then we have the equality
			\[
				\tau_*\mleft(\cW(\cD)\mright)	= (r + 1) \, W(\tau_*\mleft(\cD\mright)).
			\]
		\end{thm}
		
		The expression on the right-hand side of the equation is $r + 1$ times the tropical Weierstrass weight in Definition~\ref{def:w_weight_component_intro}, extended to more general subsets of $\Gamma$, see Section~\ref{subsec:W_measure} for details.
		
		This statement, which involves the metric of $\Gamma$ in a crucial way, gives an essentially complete description of the behavior of Weierstrass points in the tropical limit. In particular, if the limit divisor is \wfinite{}, then for every family $(X_t)_{t \neq 0}$ of smooth proper curves approaching a stable curve with dual metric graph $\Gamma$, the limit Weierstrass points are precisely described by the tropical Weierstrass divisor. This rigidity type theorem on the limiting behavior of Weierstrass points allows us to give a precise count of the number of Weierstrass points going to the nodes or to the smooth parts of a limit stable curve $X_0$ on the given stratum of $\mgbarg{g}$ along the given direction from which the family $(X_t)_{t \neq 0}$ approaches $X_0$. Moreover, as a special case, the theorem also applies in the context of arithmetic geometry in which the curve $X$ is defined over a finite extension of $\QQ_p$. 
		As we will show in Section~\ref{subsec:clls_tropicalization}, this theorem holds as well over a field $\K$ of positive characteristic provided that the gap sequence of $\cL$, defined as the sequence of orders of vanishing of the global sections of $\cL$ at a general point of $X$, is the standard sequence $0, 1, \dots, r$. (In this case, $\cL$ is called classical~\cite{Lak81, Nee84}.)
		
		\smallskip
		
		We provide natural extensions and refinements of the above results to non-necessarily complete linear series and to the setting of augmented metric graphs, which, from the degeneration perspective, corresponds to the situation where the limit stable curve has irreducible components of possibly positive genus. The data of genera of the components is encoded in a genus function $\g$ defined on vertices and extended by zero to the whole metric graph. Since a given vertex of positive genus hides information about the geometry of the corresponding component, it turns out that there will be an ambiguity when talking about the Weierstrass locus of a divisor $D$. In fact, the right setup in this context is a divisor $D$ endowed with the data of a closed sub-semimodule $M$ of $\Rat(D)$, which plays the role of a (not necessarily complete) linear series on the augmented metric graph.
		
		In this regard, first, we use the weights defined in Definition~\ref{def:w_weight_component_intro} with a relevant notion of divisorial rank associated to the sub-semimodule that we further modify by including the data of the genus function. We get Theorem~\ref{thm:number_w_points_M}, which provides a global count of weights in this setting.
		
		To the question of {\em whether it is still possible to associate a natural Weierstrass locus to a divisor in the augmented setting}, we provide an answer by introducing two special classes of semimodules, the {\em generic semimodule} associated to any divisor (see Section~\ref{subsec:weierstrass_points_augmented_graphs_gen}), and the {\em canonical semimodule} associated to the canonical divisor on an augmented metric graph (see Section~\ref{subsec:weierstrass_points_augmented_graphs_can}). Both of them require some level of genericity, which we properly justify in Section~\ref{subsec:justification_definition_semimodule_genus_function} using the framework of metrized complexes.
		
		The case of the canonical divisor on an augmented metric graph is particularly interesting as it reveals new facets of divisor theory in the augmented setting. We associate a canonical linear series to any augmented metric graph, show that it has the appropriate rank, and study its Weierstrass locus. To justify the definition and prove these results, we use the setting of metrized complexes and their divisor theory from~\cite{AB15}. Using that framework, we show that the canonical linear series on an augmented metric graph is the tropical part of the canonical linear series on any metrized complex with that underlying augmented metric graph, provided that the markings associated to edges on the curves of the metrized complex are in general position. It is interesting to note that this is the assumption made in the works by Esteves and coauthors~\cite{EM02, ES07}, and our results here complement these works by developing the tropical part of the story in greater generality.
		
		As we show in Theorem~\ref{thm:weierstrass_tropicalization_extended}, the statement of Theorem~\ref{thm:weierstrass_tropicalization_intro} remains valid in these settings (when including the genera of points of $A$ on the right-hand side of the stated equality). The following theorem is a direct application of our results on Weierstrass weights for an augmented metric graph. We use the setting of tropicalization preceding Theorem~\ref{thm:weierstrass_tropicalization_intro}.
		
		\begin{thm} \label{thm:bound_weierstrass_tropicalization_intro}
			Suppose $\cD$ is a divisor on an algebraic curve $X$ over an algebraically closed non-Archimedean field $\K$ of characteristic zero with a non-trivial valuation and a residue field of arbitrary characteristic. Let $(\Gamma, \g)$ be an (augmented) skeleton of $X^{\an}$. Let $H$ be a vector space of global sections of $\cO(\cD)$ of rank $r$ and denote by $\cW(H)$ the Weierstrass divisor of $H$. Let $M$ be the tropicalization of $H$. Then, for any connected, closed subset $A \subseteq \Gamma$ whose intersection with the Weierstrass locus of $(M, \g)$ is a union of components of the Weierstrass locus, we have the bound
			\[
				\degrest{\cW(H)}{\tau^{-1} (A)} \geq \mleft(r^2 + r\mright) \mleft(g(A) + \sum_{x \in A} \g(x) \mright).
			\]
		\end{thm}
		
		The proof of this theorem will be given in Section~\ref{subsec:proofs_tropicalization_theorems}. As in the case of Theorem~\ref{thm:weierstrass_tropicalization_intro}, the statement holds as well over a field $\K$ of positive characteristic provided the gap sequence of $H$ is the standard sequence.
		
		In the case the tropical Weierstrass of $(M, \g)$ is finite, this inequality holds for any closed subset $A \subseteq \Gamma$. In particular, we have the following application to stable curves: suppose $X_0$ is a stable curve with dual augmented graph $(G, \g)$, and suppose $(X_t)$ is a family degenerating to $X_0$ with tropicalization $(\Gamma, \g)$.
		If the locus of canonical Weierstrass points of $(\Gamma, \g)$ is finite, then for every connected subgraph $A$ of $G$, the number of limit Weierstrass points lying on components and nodes of $X_0$ that correspond to vertices and edges of $A$, respectively, is at least $(g^2 - g) \mleft(g(A) + \sum_{v \in A} \g(v) \mright)$.
		
		\medskip
		
		Semimodules inside $\Rat(D)$ that come from the tropicalization of linear series verify an extra set of properties. These are thoroughly studied in recent works~\cite{AG22} and~\cite{JP22} that develop a combinatorial theory of (limit) linear series. In particular, such a semimodule $M$ of rank $r$ satisfies the following:
		\begin{enumerate}
			\item[($\star$)] \label{item:star}
			For each point $x$ in $\Gamma$ and any unit tangent direction $\nu \in \T_x(\Gamma)$, the set of slopes taken by functions in $M$ has size $r + 1$.
		\end{enumerate}
 		(We refer to Section~\ref{app:slopes} for more details.)
		
		In Section~\ref{sec:weierstrass_locus_combinatorial_linear_series}, we associate a refined notion of Weierstrass divisor to any divisor $D$ and any closed sub-semimodule $M \subseteq \Rat(D)$ that verifies the above property. In particular, in this refined setting, there are only finitely many Weierstrass points, and the sum of coefficients of the Weierstrass points lying on any connected component $C$ of the Weierstrass locus of $M$, according to the earlier definition, is the Weierstrass weight of $C$ multiplied by $r + 1$. The definition takes into account the higher orders of vanishing of the combinatorial limit linear series, and is closer to the spirit of the algebraic definition of Weierstrass weights on curves. In contrast, the refined Weierstrass divisor in this level of generality might have negative coefficients in case the sub-semimodule does not come from geometry.
		
		Using this together with the results proved in Section~\ref{sec:tropicalization_Weierstrass_points}, discussed below, we provide a proof of Theorem~\ref{thm:weierstrass_tropicalization_intro} and its extensions to the augmented and incomplete settings.
	
	\subsection{Tropicalization of Weierstrass divisors}
		
		The proof of our comparison results, Theorem~\ref{thm:weierstrass_tropicalization_intro} and its extension Theorem~\ref{thm:weierstrass_tropicalization_extended}, makes use of the results proved in Section~\ref{sec:tropicalization_Weierstrass_points}. An earlier version of these results was written by the first author around 2014.
		
		Let $X$ be a smooth proper curve defined over $\K$. Let $\cD$ be a divisor of degree $d$ on $X$ and let $\cL = \cO(\cD)$ be the corresponding line bundle. Let $H \subseteq H^0(X, \cL)$ be a space of sections of dimension $r + 1$ and denote by $\cW = \cW(\cD, H)$ the corresponding Weierstrass divisor. We assume that the gap sequence of $H$ is the sequence $0, 1, \dots, r$, that is, for a general point $x \in X(\K)$, the orders of vanishing of sections of $\cL$ in $H$ are $0, 1, \dots, r$. Let $\tau$ be the tropicalization map from $X$ to $\Gamma$. We describe the tropicalization $W = \tau_*(\cW)$. The divisor $\cW$ is equal to $(r + 1)\cD +\div(\Wron_{\cF})$ for a section $\Wron_{\cF}$ of the sheaf $\omega!_X^{\otimes r(r + 1)/2}$ called the Wronskian. The sheaf $\omega!_X^{\otimes r(r + 1)/2}$ admits a natural norm; using this norm, we can tropicalize the section $\Wron_{\cF}$, and define a rational function $F = \Trop(\Wron_{\cF}) \colon \Gamma \to \R$. Denote by $K$ the canonical divisor of $(\Gamma, \g)$. Using the slope formula for sections of powers of the canonical sheaf, Lemma~\ref{lem:slope_diff}, it is shown in Theorem~\ref{thm:redW_general} that for any $x \in \Gamma$, we have
 		\[ W(x) = (r + 1) D(x) + \frac{r(r + 1)}2 K(x) - \sum_{\nu \in \T_x(\Gamma)} \slope_{\nu} F. \]
		It follows from the results proved in Section~\ref{sec:annuli}, that for a point $x \in \Gamma$ and $\nu \in \T_x(\Gamma)$, if
		\begin{itemize}
			\item either, the residue field $\k$ is of characteristic zero,
				
			\item or, the sequence $s^\nu_0, \dots, s^{\nu}_r$ forms an interval, that is, $s^\nu_j = s^\nu_0 + j$,
		\end{itemize}
		then the slope $\slope_{\nu} F$ is given by the sum $s^\nu_0 + \dots + s^\nu_r$, see Proposition~\ref{prop:reduction_inteval} and Theorem~\ref{thm:redW_general}. This result is needed to prove in Section~\ref{sec:weierstrass_locus_combinatorial_linear_series} our comparison results between tropical and algebraic Weierstrass loci.
		
		Note that over a field $\K$ of equicharacteristic zero, the first item in the above condition is verified, and we get all the coefficients $W(x)$, for $x \in \Gamma$, as
		\[ W(x) = (r + 1) D(x) + \frac{r(r + 1)}{2} K(x) - \sum_{\nu \in \T_x(\Gamma)} \sum_{i = 0}^r s_i^\nu, \] 
		see Theorem~\ref{thm:redW_zero}.
	
	\subsection{Previous work} \label{subsec:previous_work}
		
		The study of Weierstrass points from a tropical perspective was initiated by Baker~\cite[\S~4]{baker2008specialization}.
		Baker defines Weierstrass points for graphs and metric graphs, and uses his Specialization Lemma~\cite[Lem.~2.8]{baker2008specialization} to prove an essential compatibility with Weierstrass points on stable curves---namely, that the tropicalization of the algebraic Weierstrass locus is a subset of the tropical Weierstrass locus.
		To be more precise, for a divisor $\cD$ on a non-Archimedean curve with Weierstrass divisor $\cW(\cD)$, if $\tau_*(\cD)$ has the same rank as $\cD$ and $\cD$ has classical gap sequence, then we have an inclusion
		\[
			\abs{\tau_*\mleft(\cW(\cD)\mright)} \subseteq \Wloc(\tau_*\mleft(\cD\mright)),
		\]
		which may be strict in general. (This is stated for the canonical divisor in {\em loc. cit.}, but the proof works in greater generality.)
		This statement has strong implications for the behavior of Weierstrass points on a family of degenerating Riemann surfaces, and for $p$-adic reduction of curves over $\QQ_p$, discussed earlier in the introduction.
		Indeed, Baker motivates his study of Weierstrass points on graphs with several results from the arithmetic geometry of modular curves, in particular, as a way to decide whether certain cusps are Weierstrass points, c.f.~\cite{Ogg78, LN64, Atk67, AP03}.
		
		The question of how to determine the tropicalization of Weierstrass points on a non-Archimedean curve was settled in \cite{Ami14}; these results appear in Section~\ref{sec:tropicalization_Weierstrass_points} and are used to prove our comparison results. Regarding Question~\ref{question:intrinsinc_weight}, a previous work by Brugallé and López de Medrano~\cite{BL12} treated the case of tropical plane curves, by associating weights to inflection components (using Newton polygons), and proving a version of Theorem~\ref{thm:weierstrass_tropicalization_intro} for plane curves (it is not clear to the authors how the results are related in this case). The question of determining tropical Weierstrass loci and their weights, and the way to properly count them in the tropical setting remained however open in general. The work~\cite{Richman18} by the third-named author studies Weierstrass points on tropical curves.
		Although the tropical Weierstrass locus may be infinite in general, \cite{Richman18} shows that for a generic divisor class (i.e., lying in a nonempty open subset of $\Pic^d$), this locus is finite, and moreover computes its cardinality.
		It is worth mentioning that important divisor classes such as the canonical divisor are non-generic, so they are not covered by the methods of~\cite{Richman18}.
		The way tropical Weierstrass points distribute when the degree of divisor classes tend to infinity is studied in~\cite{Ami14, Richman18}.
		For a more thorough discussion of how divisor theory on graphs is connected to the degeneration of smooth curves to nodal curves, with various applications, see the survey by Baker--Jensen \cite{BJ16}, in particular Section 12.
		
		
		For an extensive and informative survey describing the history and applications of Weierstrass points, starting with Weierstrass and H\"{u}rwitz~\cite{Weierstrass67, Hur92} in the 1800s, see Del Centina~\cite{Del08}.
		The study of Weierstrass points on stable curves was initiated by Eisenbud and Harris~\cite{EH87-WP}, who proved results on nodal curves of {\em compact type}, i.e., curves whose dual graph is a tree.
		This work served as an application of their newly-developed theory of limit linear series~\cite{EH86}. They moreover raised the question of constructing a moduli space parametrizing all possible limit Weierstrass divisors of a given stable curve, a problem that has been widely open since then.
		
		Moving beyond stable curves of compact type, Lax~\cite{Lax87-rational} studied Weierstrass points on stable curves consisting of one rational component with nodes; in this case, the dual graph is a single vertex with self-loops.
		(The term {\em tree-like} is used in the literature to describe curves whose dual graph consists of a tree after removing self-loops.)
		A further breakthrough came with Esteves--Medeiros~\cite{EM02} who worked with stable curves with two components, i.e., curves whose dual graph is a dipole graph.
		(We refer to Section~\ref{sec:dipole_augmented} for a discussion of our results applied to dipole graphs and the connection to~\cite{EM02}.)
		Esteves--Salyehan~\cite{ES07} studied further cases of nodal curves, including when the dual graph is a complete graph.
		Cumino--Esteves--Gatto~\cite{CEG08} studied limits of {\em special} Weierstrass points on certain stable curves, i.e., Weierstrass points with weight at least two. The problem of describing limits of Weierstrass points away from the nodes in a given one-parameter family in characteristic zero is addressed in~\cite{Esteves98}.
		
		Although not directly related to the results of this paper, we mention that other works treat the case of irreducible Gorenstein curves, and associate Weierstrass weights to their singular points, see e.g.~\cite{LW90, CS94, GL95, BG95} and the references there. It might be possible to use tropical geometry to describe these weights.
		
		
		Weierstrass points have appeared in other interesting work on moduli spaces of curves.
		Arbarello~\cite{Arb74} studied subvarieties of the moduli space of curves cut out by Weierstrass points;
		further results were found in Lax~\cite{Lax75} and Diaz~\cite{Dia85}.
		Eisenbud--Harris~\cite{EH87-KD} showed that the moduli space of curves has positive Kodaira dimension, using loci of Weierstrass points as part of their argument.
		Cukierman(--Fong)~\cite{Cuk89, CF91} found the coefficients for the Weierstrass locus in the universal curve $\unicurve g$ of genus $g$, in a standard basis for the Picard group of $\unicurve g$.
	 
	\subsection{Organization of the text}
		 We begin with non-augmented metric graphs, then treat refinements, a choice intended to simplify the presentation and to increase readability.
		
		We define slope sets and prove Theorem~\ref{thm:consecutive_slopes_intro} in Section~\ref{sec:slope_sets}.
		
		Section~\ref{sec:proof_number_w_points} develops Weierstrass weights and the associated measure on metric graphs. We prove Theorem~\ref{thm:weight_measure}, which provides a description of the Weierstrass measure using the slopes, deduce Theorem~\ref{thm:number_w_points_intro}, establish positivity of Weierstrass weights, and discuss combinatorial graphs. The case of canonical Weierstrass locus is treated in Section~\ref{sec:canonical_locus}.
		
		
		Section~\ref{sec:weierstrass_locus_generalized_setup} extends the setting to incomplete linear series (by taking closed sub-semimodules of $\Rat(D)$) and to augmented metric graphs, justifying the definitions in the augmented setting and generalizing Theorems~\ref{thm:weight_measure} and~\ref{thm:number_w_points_intro}.
		
		
		Section~\ref{sec:weierstrass_locus_combinatorial_linear_series} associates Weierstrass divisors to combinatorial limit linear series, including the interesting case where the Weierstrass locus becomes infinite after forgetting the slopes, and shows compatibility with previous definitions.
		
		
		Section~\ref{sec:tropicalization_Weierstrass_points} establishes the needed results on tropicalization of Weierstrass divisors. Building on these, we prove Theorem~\ref{thm:weierstrass_tropicalization_intro} and its generalizations, linking tropical and classical Weierstrass divisors.
		
		
		Finally, Section~\ref{sec:further_examples_discussions} offers further examples to clarify the constructions and discusses related results and questions.
		
	
	\subsection*{Acknowledgments}
		
		We are grateful to the referees for their careful reading of the paper and their comments and suggestions, which helped improve the presentation.
		
		This material is based upon research supported by the Transatlantic Research Partnership of the Embassy of France in the United States and the FACE Foundation. We thank them, as well as the Cultural Services of the French Embassy in the U.S., and the administrations of the CNRS and the University of Washington at Seattle, for support.
		
		O.\thinspace A. thanks Math+, the Berlin Mathematics Research Center and TU Berlin, where this research was carried out. L.\thinspace G. thanks the Mathematical Foundation Jacques Hadamard (FMJH) and the French \emph{Agence Nationale de la Recherche} (project ANR-22-CE40-0014) for financial support.
	
	\subsection{Basic notation} \label{subsec:basic_notation}
		
		A {\em (combinatorial) graph} $G = (V, E)$ is defined by a set of vertices $V$ and a set of edges $E$ between certain vertices. In the current paper, graphs will always be taken to be finite and connected. Moreover, they will allow loops and multiple edges.
		
		A {\em metric graph} is a compact, connected metric space $\Gamma$ verifying the following properties:
		
		\begin{enumerate}[(i)]
			\item For every point $x \in \Gamma$, there exist a positive integer $n_x$ and a real number $r_x > 0$ such that the $r_x$-neighborhood of $x$ is isometric to the star of radius $r_x$ with $n_x$ branches.
					
			\item The metric on $\Gamma$ is given by the path metric, i.e., for points $x$ and $y$ in $\Gamma$, the distance between $x$ and $y$ is the infimum (in fact minimum) length of any path from $x$ to $y$.
		\end{enumerate}
		
		The integer $n_x$ above is called the {\em valence} of $x$ and is denoted by $\val(x)$.
		
		Given a graph $G = (V, E)$ and a length function $\ell \colon E \to (0, +\infty)$ assigning to every edge of $G$ a positive length, we can build from this data a metric graph $\Gamma$ by gluing a closed interval of length $\ell(e)$ between the two endpoints of the edge $e$, for every $e \in E$, and endowing $\Gamma$ with the path metric. The space $\Gamma$ is then called the {\em geometric realization} of the pair $(G, \ell)$.
		
		A {\em model} of a metric graph $\Gamma$ is a pair $(G, \ell)$ consisting of a graph $G = (V, E)$ and a length function $\ell \colon E \to (0, +\infty)$ such that $\Gamma$ is isometric to the geometric realization of $(G, \ell)$. By an abuse of notation, we also call $G$ a model of $\Gamma$.
		
		For a metric graph $\Gamma$ and a point $x \in \Gamma$, the tangent space $\T_x(\Gamma)$ is defined as the set of all unit outgoing tangent vectors to $\Gamma$ at $x$. This is a finite set of cardinality $\val(x)$. If $G = (V, E)$ is a loopless model for $\Gamma$ such that $x \in V$, then $\T_x(\Gamma)$ is in one-to-one correspondence with the edges of $G$ incident to $x$. Through this natural bijection, a tangent direction $\nu$ is said to be {\em supported by} the corresponding edge $e \in E$.
		
		Each edge $e$ supports two tangent directions, which belong to either endpoint of $e$, respectively. If $\nu$ is one of those tangent directions, the opposite direction is denoted by $\overline \nu$.
		For $\nu \in \T(\Gamma)$, we denote by $x_\nu$ the point $x$ with $\nu \in \T_x(\Gamma)$.
		
		In this paper, the rank $r$ of divisors will always be nonnegative, and all the semimodules will be assumed to be nonempty.
			
\section{Slope sets} \label{sec:slope_sets}
	
	The main result of this section is Theorem~\ref{thm:consecutive_slopes_intro}. We first recall some terminology for divisors and functions on metric graphs.
	
	Given a metric graph $\Gamma$, let $\Div(\Gamma)$ denote the {\em group of divisors} of $\Gamma$, which is the free abelian group generated by points $x \in \Gamma$. Let $\PLZ(\Gamma)$ denote the set of real-valued piecewise linear functions on $\Gamma$ whose slopes are all integers. Given a function $f \in \PLZ(\Gamma)$, let $\div(f)$ denote the {\em principal divisor} of $f$, defined as
	\[
		\div(f) \coloneqq \sum_{x \in \Gamma} a_x (x) \qquad\text{where} \qquad a_x = - \sum_{\nu \in \T_x(\Gamma)} \slope_\nu f(x).
	\]
	
	Let $D$ be a divisor of rank $r$ on $\Gamma$. Recall that the rank of $D$ is defined as the maximum integer $r \geq -1$ such that for any effective divisor $E$ of degree $r$, there exists an element $f \in M$ with $\div(f) + D - E \geq 0$. Let $\Rat(D)$ denote the set of {\em rational functions in the complete linear series of $D$} defined as
	\[
		\Rat(D) \coloneqq \mleft\{f \in \PLZ(\Gamma) \st D + \div(f) \geq 0\mright\}.
	\]
	Given a point $x \in D$, there is a unique representative $f_x$ of the linear series of $D$ defined by
	\begin{equation} \label{eq:defition_reduced_function}
		f_x \coloneqq \min_{\substack{f \in \Rat(D) \\ f(x) = 0}} f.
	\end{equation}
	The corresponding divisor $D + \div(f_x)$, denoted by $D_x$, is the (unique) {\em $x$-reduced divisor} linearly equivalent to $D$. This statement is a consequence of the maximum principle, see e.g.~\cite[Lem.~4.11]{BS13}.
	
	\begin{defi}[Slope sets and minimum slopes] \label{defi:slopes}
		Let $D$ be a divisor of nonnegative rank on $\Gamma$. Given a point $x \in \Gamma$ and a tangent direction $\nu \in \T_x(\Gamma)$, let $\fS^{\nu}(D)$ denote the {\em slope set}
		\[
			\fS^{\nu}(D) \coloneqq \mleft\{\slope_\nu f(x) \st f \in \Rat(D)\mright\}.
		\]
		Let $s_0^\nu(D)$ denote the {\em minimum slope} in the slope set $\fS^\nu(D)$, i.e.,
		\[
			s_0^{\nu}(D) \coloneqq \min \mleft\{\slope_\nu f(x) \st f \in \Rat(D)\mright\}.
		\]
		When the divisor $D$ is clear from context, we will simply use $s_0^\nu$ to denote $s_0^\nu(D)$.
	\end{defi}

	The following result allows to retrieve information about the slopes of rational functions in the linear series $\Rat(D)$ along tangent directions in $\Gamma$.
	
	\begin{thm} \label{thm:consecutive_slopes_intro}
		Let $D$ be a divisor of rank $r$ on $\Gamma$. We take a model for $\Gamma$ whose vertex set contains the support of $D$. Let $x \in \Gamma$ be a point and $\nu \in \T_x(\Gamma)$ be a tangent direction. Recall that we denote by $\slope_\nu f(x)$ the slope of $f$ along the tangent direction $\nu$ based at $x$.
		\begin{enumerate}[(a)]
			\item If the open interval $(x, x + \varepsilon \nu)$ is disjoint from $\Wloc(D)$ for some $\varepsilon > 0$, then the set of slopes $\mleft\{\slope_{\nu} f(x) \st f \in \Rat(D)\mright\}$ consists of $r + 1$ consecutive integers $\{s^\nu_0, s^\nu_0 + 1, \dots, s^\nu_0 + r\}$.
			
			\item If the open interval $(x, x + \varepsilon \nu)$ is contained in $\Wloc(D)$ for some $\varepsilon > 0$, then the set of slopes $\mleft\{\slope_\nu f(x) \st f \in \Rat(D)\mright\}$ is a set of consecutive integers of size at least $r + 2$.
		\end{enumerate}
	\end{thm}
	
	In preparation for the proof, we start by the following useful fact.
	
	\begin{remark} \label{rem:fact_set_slopes_constant_closed_subsemimodule}
		Since $\Rat(D)$ is finitely generated (see~\cite[Thm.~6]{haase2012linear}), the slope set is directionally locally constant, in the following sense. For any point $x \in \Gamma$ and any unit tangent direction $\nu \in \T_x(\Gamma)$, there is a positive length interval $[x, y]$ supporting $\nu$, with $\nu$ pointing toward $y$, on which the slope set $\fS^{\nu_z}(D)$ is constant, for $z \in [x, y)$, with $\nu_z \in \T_z(\Gamma)$ parallel to $\nu$ and pointing toward $y$.
		
		The same holds for any closed sub-semimodule $M \subseteq \Rat(D)$. We omit the proof.
	\end{remark}
	
	\begin{lemma} \label{lem:minimum_number_slopes}
		Suppose $D$ is a divisor of rank $r \geq 0$. Then, for every $x \in \Gamma$ and every $\nu \in \T_x(\Gamma)$, there are at least $r + 1$ integers in the set of slopes $\mleft\{\slope_{\nu}f(x) \st f \in \Rat(D)\mright\}$.
	\end{lemma}
	
	\begin{proof}
	 	Let $x_1, \dots, x_r$ be a set of distinct points in the branch incident to $x$ in the direction of $\nu$ sufficiently close to $x$. There exists a function $f \in \Rat(D)$ such that
	 	\[
	 		D + \div(f) \geq (x_1) + \dots + (x_r). 
	 	\]
	 	The function $f$ changes slope at the points $x_1, \dots, x_r$. Thanks to Remark~\ref{rem:fact_set_slopes_constant_closed_subsemimodule}, each of the slopes taken at $x_j$ in the direction of $\nu$ can be obtained as the slope of a function in $\Rat(D)$ at $x$ along $\nu$.
	\end{proof}
	
	The minimum slope $s^\nu_0(D)$ is related to the reduced divisor $D_x$ at $x$.
	
	\begin{lemma} \label{lem:formula_reduced_divisor}
		Let $D$ be an effective divisor on $\Gamma$, and $x$ a point of $\Gamma$. Let $D_x$ be the $x$-reduced divisor linearly equivalent to $D$. \hfill
		\begin{enumerate}[(a)]
			\item 
			Let $f_x$ be the rational function defined in~\eqref{eq:defition_reduced_function} and satisfying $\div(f_x) + D = D_x$, then, for any outgoing tangent vector $\nu \in \T_x(\Gamma)$,
			\[
				s_0^\nu(D) = \slope_\nu f_x(x).
			\]
			
			\item
			The coefficient of $D_x$ at $x$ satisfies
			\begin{equation*}
		 		D_x(x) = D(x) - \sum_{\nu \in \T_x(\Gamma)} s_0^\nu(D) .
			\end{equation*}
		\end{enumerate}
	\end{lemma}
	
	\begin{proof}
	 	The first result is obtained using Remark~\ref{rem:fact_set_slopes_constant_closed_subsemimodule} again, by observing that $f_x = \min h$ for $h \in \Rat(D)$ verifying $h(x) = 0$. The second result is a direct consequence of (a) and the definition of the principal divisor $\div(f_x)$.
	\end{proof}
	
	We now turn to the proof of Theorem~\ref{thm:consecutive_slopes_intro}. 
	Let $D$ be a divisor of rank $r$ on $\Gamma$. 
	Recall (Definition~\ref{def:weierstrass_points_intro}) that the Weierstrass locus of $D$, denoted by $\Wloc(D)$, is the subset of $\Gamma$ formed by the points $x$ such that there exists an effective divisor $E \sim D$ with $E(x) \geq r + 1$. 
	Equivalently, $\Wloc(D)$ is defined in terms of reduced divisors as
	\begin{equation*}
		\Wloc(D) = \mleft\{x \in \Gamma \st D_x(x) > r\mright\},
	\end{equation*}
	where $D_x$ denotes the $x$-reduced divisor linearly equivalent to $D$.
	
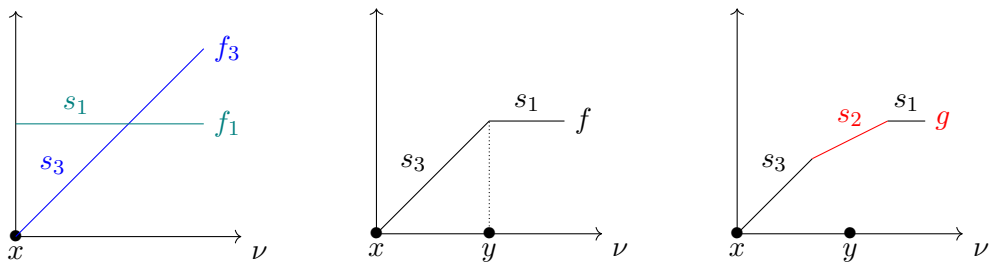
\begin{figure}[h!]
		\begin{minipage}{0.3\textwidth}
			\centering
			\begin{tikzpicture}[scale=.75]
				\draw[->] (0,0) -- (3,0);
				\draw[->] (0,0) -- (0,3);
				\draw (0,0) node[below]{$x$} node{$\bullet$};
				\draw (3,0) node[below right]{$\nu$};
				\draw[color = blue] (0,0) -- (2.5,2.5);
				\draw[color = teal] (0,1.5) -- (2.5,1.5);
				\draw (2.5,2.5) node[right, color = blue]{$f_3$};
				\draw (2.5,1.5) node[right, color = teal]{$f_1$};
				\draw(0.8,0.7) node[above left, color = blue]{$s_3$};
				\draw(0.8,1.5) node[above, color = teal]{$s_1$};
			\end{tikzpicture}
		\end{minipage}
		\begin{minipage}{0.3\textwidth}
			\centering
			\begin{tikzpicture}[scale=.75]
				\draw[->] (0,0) -- (3,0);
				\draw[->] (0,0) -- (0,3);
				\draw (0,0) node[below]{$x$} node{$\bullet$};
				\draw (3,0) node[below right]{$\nu$};
				\draw (0,0) -- (1.5,1.5) -- (2.5,1.5);
				\draw (2.5,1.5) node[right]{$f$};
				\draw (1.5,0) node[below]{$y$} node{$\bullet$};
				\draw[densely dotted] (1.5,0) -- (1.5,1.5);
				\draw (0.8,0.7) node[above left]{$s_3$};
				\draw (2, 1.5) node[above]{$s_1$};
			\end{tikzpicture}
		\end{minipage}
		\begin{minipage}{0.3\textwidth}
			\centering
			\begin{tikzpicture}[scale=.75]
				\draw[->] (0,0) -- (3,0);
				\draw[->] (0,0) -- (0,3);
				\draw (0,0) node[below]{$x$} node{$\bullet$};
				\draw (3,0) node[below right]{$\nu$};
				\draw (0,0) -- (1,1);
				\draw[color = red] (1,1) -- (2,1.5);
				\draw (2,1.5) -- (2.5,1.5);
				\draw (2.5,1.5) node[right, color = red]{$g$};
				\draw (1.5,0) node[below]{$y$} node{$\bullet$};
				\draw (0.8,0.7) node[above left]{$s_3$};
				\draw (1.5,1.3) node[above, color = red]{$s_2$};
				\draw (2.25,1.5) node[above]{$s_1$};
			\end{tikzpicture}
		\end{minipage}
		\caption{Construction of the functions $f$ and $g$ using functions $f_1$ and $f_3$ taking slopes $s_1 < s_3$.}
		\label{fig:intermediate_slope}
	\end{figure}

	\begin{proof}[Proof of Theorem~\ref{thm:consecutive_slopes_intro}]
		We first assume that the open interval $(x, x + \varepsilon \nu)$ is disjoint from $\Wloc(D)$. 
		Along the branch incident to $x$ in the direction of $\nu$, there is a small segment on which $s^\nu_0$ is the slope of a function in $\Rat(D)$, and it is the smallest slope taken by a function of $\Rat(D)$ on this segment.
		If a slope of $s^\nu_0 + r + 1$ or larger is achieved at $x$, then, again on a small segment, it will be achieved at any point of that segment. This means that on the interior of this segment, the two minimum outgoing slopes at every point are $s^\nu_0$ and $- s^\nu_0 - s'$ with $s' \geq r + 1$. Therefore, by Lemma~\ref{lem:formula_reduced_divisor}, we infer that on the interior of this segment, the reduced divisor at each point has coefficient at least $r + 1$. This contradicts the assumption that a neighborhood $(x, x + \varepsilon \nu)$ is disjoint from $\Wloc(D)$, and shows that the highest possible slope is $s^\nu_0 + r$. Combining this with Lemma~\ref{lem:minimum_number_slopes}, the slopes achieved at $x$ along $\nu$ must be precisely $s^\nu_0, s^\nu_0 + 1, \ldots, s^\nu_0 + r$. This proves (a).
		
		We now assume that $(x, x + \varepsilon \nu) \subset \Wloc(D)$ and that $\varepsilon$ is small enough so that the set of slopes of functions of $\Rat(D)$ along $\nu$ is constant on this interval. By Lemma~\ref{lem:formula_reduced_divisor}, this means that on the interior of a small segment starting at $x$, the two minimum outgoing slopes at every point are $s^\nu_0$ and $-s^\nu_0 - s'$ with $s' \geq r + 1$. Therefore, close to $x$, a slope of at least $s^\nu_0 + r + 1$ is achieved by a function in $\Rat(D)$. To prove (b), it is thus sufficient to show that the set of slopes $\slope_\nu f(x)$ of functions $f \in \Rat(D)$ is always made up of consecutive integers. Take $s_1 < s_2 < s_3$ to be three integers, and suppose that for $i \in \{1, 3\}$ there exists a function $f_i \in \Rat(D)$ such that $\slope_\nu f_i(x) = s_i$. Using $f_1$, $f_3$ and tropical operations, it is easy to construct a function $f$ taking slopes $s_3$ and then $s_1$ away from $x$, changing slope at a point we denote by $y$ (see Figure~\ref{fig:intermediate_slope}).
		We can then ``chop up'' the graph of $f$ to construct a function $h$ equal to $f$ everywhere except on a small interval around $y$ where it takes slope $s_2$.
		Since $(x, x + \varepsilon \nu)$ is disjoint from the support of $D$, $h$ still belongs to $\Rat(D)$.
		The assumption made on $\varepsilon$ at the beginning ensures that in fact there exists a function $f_2 \in \Rat(D)$ taking slope $s_2$ at $x$ along $\nu$, which concludes the argument.
	\end{proof}
	
	\begin{remark}
		In particular, note that along a given unit tangent vector $\nu$ attached to a point $x$, the slopes $\slope_\nu f(x)$ for $x \in \Rat(D)$ always form a set of consecutive integers. Moreover, if $t$ is a positive integer such that for every $x \in e$, for $e$ an edge of some model of $\Gamma$, the $x$-reduced divisor $D_x$ satisfies $D_x \geq t \, (x)$, then for any $x \in \mathring e$, the set of slopes $\mleft\{\slope_\nu f(x) \st f \in \Rat(D)\mright\}$ contains at least $t + 1$ consecutive integers. This claim is analogous to~\cite[Thm.~14]{amini2013reduced} and is proved using Theorem~3 of the same paper, which gives a concrete description of the variations of the reduced divisor $D_x$ with respect to $x$. See also~\cite[\S~6.6]{AG22}.
	\end{remark}
	
	For future use, we note the following generalization of part (b) of Lemma~\ref{lem:formula_reduced_divisor}.
	
	\begin{prop} \label{prop:bound_red_degree}
		Suppose $D$ is a divisor of rank $r \geq 0$. Then, for any closed, connected subset $A \subseteq \Gamma$, we have
		\[
			\degrest{D}{A} - \sum_{\nu \in \partialout A} s_0^\nu(D) \geq r.
		\]
	\end{prop}
	
	\begin{proof}
		Let $E$ be an effective divisor of degree $r$, with support contained in $A$. Since $D$ has rank $r$, there exists a function $f \in \Rat(D)$ such that $D + \div(f) \geq E$. Evaluating the respective degrees restricted to $A$ yields
		\[
			\degrest{D}{A} - \sum_{\nu \in \partialout A} \slope_\nu f(x_\nu) \quad \geq \quad \degrest{E}{A} = r,
		\]
		where, we recall, $x_\nu$ is the point $x$ of $\Gamma$ with $\nu \in \T_x(\Gamma)$. By definition of the minimum slope $s_0^\nu(D)$, we have $s_0^\nu(D) \leq \slope_\nu f(x_\nu)$ for each $\nu \in \partialout A$, so the result follows.
	\end{proof}

\section{Weierstrass weights} \label{sec:proof_number_w_points}
	
	Using the structure of slope sets in $\Rat(D)$, we prove Theorem~\ref{thm:number_w_points_intro}, which will follow from the more general Theorem~\ref{thm:weight_measure}.
 
	\subsection{Definition of weights and basic properties of the Weierstrass locus}
			
		We start by establishing basic properties of Weierstrass loci.
		(Definition~\ref{def:weierstrass_points_intro})
		The Weierstrass locus $\Wloc(D)$ is defined as the set of points $x$ in $\Gamma$ such that there exists an effective divisor $E$ in the linear system of $D$ whose coefficient at $x$ is at least $r + 1$. This is equivalent to requiring that $D_x(x) \geq r + 1$. Let us now recall Definition~\ref{def:w_weight_component_intro} from the introduction. Given a connected component $C$ of the Weierstrass locus $\Wloc(D)$, the tropical Weierstrass weight of $C$ is defined as
		\begin{equation*}
			\weight(C) = \weight(C; D) \coloneqq \degrest{D}{C} + \mleft(g(C) - 1\mright) \, r - \sum_{\nu \in \partialout C} s_0^\nu(D)
		\end{equation*}
		where $\degrest{D}{C} = \sum_{x \in C} D(x)$ is the degree of $D$ in $C$, $g(C) = \dim H_1(C, \RR)$ is the genus of $C$, $\partialout C$ is the set of outgoing unit tangent directions from $C$, and $s_0^\nu(D)$ is the minimum slope at $x$ along a tangent direction $\nu$, as defined in Definition~\ref{defi:slopes}.
		
		We already noted that the tropical Weierstrass locus $\Wloc(D)$ is independent of the choice of the divisor $D$ in its linear equivalence class. In fact, this is also the case, though not as obviously, for the tropical Weierstrass weights.
		
		\begin{prop} \label{prop:tropical_weight_independent_linear_equivalence}
			Let $D$ be a divisor on $\Gamma$ and $h \in \Rat(D)$. Then the following equality holds:
			\[ \weight(C; D + \div(h)) = \weight(C; D). \]
		\end{prop}
		
		\begin{proof}
			Given a point $x \in \Gamma$ and a tangent direction $\nu \in \T_x(\Gamma)$, the slope set introduced in Definition~\ref{defi:slopes} satisfies
			\[ \fS^\nu(D + \div(h)) = \fS^\nu(D) - \slope_\nu h(x). \]
			Therefore, the minimum slope along $\nu$ also transforms nicely:
			\[ s_0^\nu(D + \div(h)) = s_0^\nu(D) - \slope_\nu h(x). \]
			Summing this equality over all points $x$ belonging to the boundary of a given connected component $C$ of $\Wloc(D)$ and over all outgoing tangent directions $\nu$, we get
			\[ \sum_{\nu \in \partialout C} s_0^\nu(D + \div(h)) = \sum_{\nu \in \partialout C} s_0^\nu(D) - \sum_{\nu \in \partialout C} \slope_\nu h(x_\nu), \]
			where $x_\nu$ is the point at which $\nu$ is based.
			
			On the other hand, when computing the degree of the restriction $\div(h) \rest{C}$, all slopes of $h$ along tangent vectors pointing toward $C$ cancel out, and we get
			\[ \degrest{\div(h)}{C} = -\sum_{\nu \in \partialout C} \slope_\nu h(x_\nu), \]
			which proves that $\weight(C; D + \div(h)) = \weight(C; D)$, as needed.
		\end{proof}
		
		The following proposition shows that $\Wloc(D)$ is topologically nice.
		
		\begin{prop} \label{prop:well_defined}
			The Weierstrass locus $\Wloc(D)$ is closed and has finitely many components. Each connected component is a metric graph.
		\end{prop}
		
		\begin{proof}
			By the continuity of variation of reduced divisors proved in~\cite[Thm.~3]{amini2013reduced}, the function $x \mapsto D_x(x)$ is upper semicontinuous, which implies that the subset $\Wloc(D)$ is closed. Then, since $\Gamma$ is compact, $\Wloc(D)$ is itself compact. Therefore, it has a finite number of connected components. The last statement follows as any connected component of a closed subset in a metric graph is itself a metric graph.
		\end{proof}
		
		\begin{remark} \label{rem:geometric_description_Weierstrass_locus}
			We have the following geometric construction of $\Wloc(D)$, which gives another proof for Proposition~\ref{prop:well_defined}. Let $\Pic^d(\Gamma)$ denote the space of divisor classes of degree $d$ on $\Gamma$, and let $\Eff^d(\Gamma)$ denote the space of effective divisor classes of degree $d$. Let $\varphi \colon \Gamma \to \Pic^{d - r - 1}(\Gamma)$ be the map defined by $\varphi(x) = [D - (r + 1) \, (x)]$.
	
			The condition that $D_x(x) > r$ is equivalent to $D_x \geq (r + 1) \, (x)$. This is in turn equivalent to the condition that the divisor class $[D - (r + 1) \, (x)]$ has an effective representative.
			Using this observation and the above terminology,
			$\Wloc(D) = \varphi^{-1} \mleft( \varphi(\Gamma) \cap \Eff^{d - r - 1}(\Gamma) \mright)$.
			In other words, $\Wloc(D)$ is described by the following pullback diagram.
			\[
				\begin{tikzcd}
					\Wloc(D) \arrow[d, hook] \arrow[r] \arrow[dr, phantom, "\scalebox{1.5}{$\lrcorner$}", very near start] & \Eff^{d - r - 1}(\Gamma) \ar[d, hook] \\
					\Gamma \ar[r, "\varphi"] & \Pic^{d - r - 1}(\Gamma)
				\end{tikzcd}
			\]
			Both $\Eff^{d - r - 1}(\Gamma)$ and $\varphi(\Gamma)$ are polyhedral subsets of $\Pic^{d - r - 1}(\Gamma)$ with finitely many facets. Thus their intersection has finitely many components, and each component is a union of finitely many closed intervals.
		\end{remark}
		
		\begin{remark} \label{rem:multiplicity_singleton}
			As before, let $D_x$ denote the $x$-reduced divisor linearly equivalent to $D$.
			Since the Weierstrass locus $\Wloc(D)$ is defined as $\mleft\{x \in \Gamma \st D_x(x) - r > 0\mright\}$, the expression $D_x(x) - r$ is a natural ``naive'' candidate for defining a tropical Weierstrass weight. In fact, this ends up being the correct definition when $x$ is an isolated component of $\Wloc(D)$. When $x$ is not an isolated component, our more technical definition of weight is required.
			
			If the singleton $\{x\}$ is a connected component of $\Wloc(D)$, then we verify that the weight of $x$ is simply given by $D_x(x) - r$. Since the genus of the component $\{x\}$ is zero, Definition~\ref{def:w_weight_component_intro} states that
			\begin{align*}
				\weight(x) = D(x) - r - \sum_{\nu \in \T_x(\Gamma)} s^\nu_0(D),
			\end{align*}
			and Lemma~\ref{lem:formula_reduced_divisor} states that
			$D_x(x)= D(x) - \sum_{\nu \in \T_x(\Gamma)} s^\nu_0(D)$. This verifies the claim.
			
			Note that this applies for every connected component of $\Wloc(D)$ if $D$ is \wfinite{}.
		\end{remark}
		
		We now give two examples of metric graphs and their Weierstrass loci. The first Weierstrass locus is finite whereas the second one is infinite.
		
		\begin{example} \label{ex:complete_graph_4}
			Suppose $\Gamma$ is the complete graph on four vertices with unit edge lengths; see Figure~\ref{fig:complete_4}. This graph has genus three, and the rank of the canonical divisor $K$ is $r = g - 1 = 2$.
			The Weierstrass locus $\Wloc(K)$ is finite and consists of the four vertex points. At a vertex $v$, the reduced divisor at $v$ is $K_v = 4 \, (v)$. Thus, $\weight(v) = K_v(v) - r = 4 - 2 = 2$.

			\begin{figure}[h!]
				\begin{minipage}{0.45\textwidth}
				\centering
				\begin{tikzpicture}[scale=.65]
					\coordinate (A) at (0,0);
					\coordinate (B) at (3.2,0);
					\coordinate (C) at (1.6,2.4);
					\coordinate (D) at (1.6,0.8);
					
					\draw (A) -- (B) -- (C) -- cycle;
					\draw (A) -- (D) -- (B);
					\draw (C) -- (D);
					
					\foreach \c in {A,B,C,D} {
						\fill (\c) circle (3pt);
					}
				\end{tikzpicture}
				\end{minipage}
				\begin{minipage}{0.45\textwidth}
				\centering
				\begin{tikzpicture}[scale=.65]
					\coordinate (A) at (0,0);
					\coordinate (B) at (3.2,0);
					\coordinate (C) at (1.6,2.4);
					\coordinate (D) at (1.6,0.8);
					
					\draw (A) -- (B) -- (C) -- cycle;
					\draw (A) -- (D) -- (B);
					\draw (C) -- (D);
					
					\foreach \c in {A,B,C,D} {
						\fill[color=red] (\c) circle (3pt);
					}
					\node[above left] at (A) {$2$};
					\node[above right] at (B) {$2$};
					\node[above left] at (C) {$2$};
					\node[above left] at (D) {$2$};
				\end{tikzpicture}
				\end{minipage}
				\caption{Complete graph on four vertices, and its Weierstrass locus $\Wloc(K)$.}
				\label{fig:complete_4}
			\end{figure}
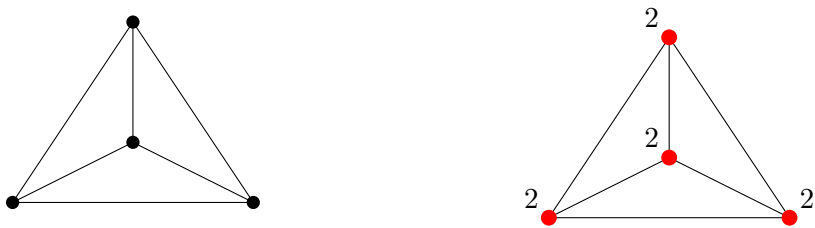
			
			We will treat the example of the complete graph on five or more vertices in Section~\ref{subsec:graph_all_weierstrass}.
		\end{example}
		
		\begin{example} \label{ex:barbell_graph}
			Suppose $\Gamma$ is the ``barbell graph'' consisting of two cycles joined by a bridge edge; see Figure~\ref{fig:barbell}. (The edge lengths may be arbitrary.) This graph has genus two, and the canonical divisor $K$ has rank $r = g - 1 = 1$.
			
			The Weierstrass locus $\Wloc(K)$ consists of the middle edge and the outer midpoint on each cycle. The latter have weight one. Besides, if $v$ is either endpoint of the bridge edge, and if $\nu \in \T_v(\Gamma)$ is one of the two tangent vectors based at $v$ that are not supported on the bridge, then $\fS^\nu(K) = \{0, 1\}$. This implies that the weight of the middle edge is also one.
			
			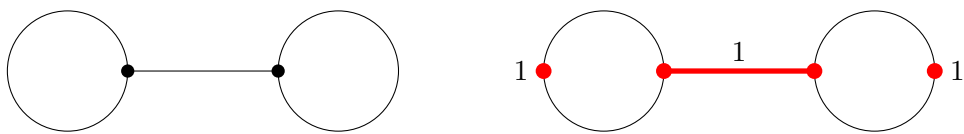
\begin{figure}[h!]
				\begin{minipage}{0.45\textwidth}
				\centering
				\begin{tikzpicture}[scale=.65]
					\coordinate (A) at (0,0);
					\coordinate (B) at (.8,0);
					\coordinate (C) at (2.8,0);
					\coordinate (D) at (3.6,0);
					
					\draw (B) -- (C);
					\draw[radius=0.8] (A) circle;
					\draw[radius=0.8] (D) circle;
					
					\foreach \c in {B,C} {
						\fill (\c) circle (3pt);
					}
				\end{tikzpicture}
				\end{minipage}
				\begin{minipage}{0.45\textwidth}
				\centering
				\begin{tikzpicture}[scale=.65]
					\coordinate (A) at (0,0);
					\coordinate (B) at (.8,0);
					\coordinate (C) at (2.8,0);
					\coordinate (D) at (3.6,0);
					\coordinate (E) at (-0.8,0);
					\coordinate (F) at (4.4,0);
					
					\draw[color=red,line width=1.2pt] (B) -- (C);
					\draw[radius=0.8] (A) circle;
					\draw[radius=0.8] (D) circle;
					
					\foreach \c in {B,C,E,F} {
						\fill[color=red] (\c) circle (3pt);
					}
					\node[above] at (1.8,0) {$1$};
					\node[left=2pt] at (E) {$1$};
					\node[right=2pt] at (F) {$1$};
				\end{tikzpicture}
				\end{minipage}
				\caption{The barbell graph and its Weierstrass locus $\Wloc(K)$.}
				\label{fig:barbell}
			\end{figure}
			
			We will show in Section~\ref{subsec:bridge_edges} that, more generally, if $e$ is a bridge edge of $\Gamma$ such that each component of $\Gamma \setminus \mathring e$ has positive genus, then $e$ is contained in the canonical Weierstrass locus.
		\end{example}
		
	\subsection{Positivity of Weierstrass weights} \label{sec:positivity}
		
		We now prove the following theorem.
		\begin{thm} \label{thm:positivity}
			Let $D$ be a divisor on $\Gamma$ with non-negative rank $r$, and let $C$ be a connected component of the Weierstrass locus $\Wloc(D)$. Then, the weight $\weight(C)$ given in Definition~\ref{def:w_weight_component_intro} is positive.
		\end{thm}
		
		\begin{proof}
			We use the notation introduced previously. Let $x$ be a point in the connected component $C$, and let $D_x$ be the $x$-reduced divisor equivalent to $D$. By definition of the Weierstrass locus $\Wloc(D)$, we have $ D_x(x) - r > 0.$ Let $f_x$ be a rational function such that $D_x = D + \div(f_x)$. We have
			\[
				D_x(x) = D(x) - \sum_{\nu \in \T_x(\Gamma)} \slope_\nu f_x(x).
			\]
			Let $A$ be any connected subgraph of $\Gamma$, and recall that $\degrest{D_x}{A}$ denotes the sum $\sum_{y \in A} D_x(y)$. For a tangent vector $\nu \in \partialout A$, let, as before, $x_\nu$ denote the associated boundary point. We have
			\[
				\degrest{D_x}{A} = \degrest{D}{A} - \sum_{\nu \in \partialout A} \slope_\nu f_x(x_\nu)
			\]
			by applying Stokes theorem to the derivative of $f_x$ on the region $A$.
			
			Because $x \in C$ and $D_x$ is effective, we have $\degrest{D_x}{C} \geq D_x(x) > r$.
			For each tangent direction $\nu \in \partialout C$, the minimum slope $s_0^\nu(D)$ satisfies $s_0^\nu(D) \leq \slope_\nu f_x (x_\nu)$ by definition (Definition~\ref{defi:slopes}). Therefore,
			\begin{align*}
				\weight(C) &= \degrest{D}{C} + (g(C) - 1) \, r - \sum_{\nu \in \partialout C} s_0^\nu(D)
				\geq \degrest{D}{C} - r - \sum_{\nu \in \partialout C} s_0^\nu(D) \\
				&\geq \degrest{D}{C} - r - \sum_{\nu \in \partialout C} \slope_\nu f_x(x_\nu)\quad = \quad \degrest{D_x}{C} - r > 0
			\end{align*}
			as claimed.
		\end{proof}
		
		The proof of Theorem~\ref{thm:positivity} shows the stronger bound $\weight(C) > g(C) \, r$.
		Besides, Theorem~\ref{thm:positivity} is addressed later, in greater generality, in Corollary~\ref{cor:weight_genus_bound}.
		
	\subsection{Weierstrass measure} \label{subsec:W_measure}
		
		We prove Theorem~\ref{thm:weight_measure} below, which will imply Theorem~\ref{thm:number_w_points_intro}.
		
		\begin{defi}
			Fix a divisor $D$ on a metric graph $\Gamma$, with Weierstrass locus $\Wloc(D)$. A subset $A \subseteq \Gamma$ is {\em $\Wloc(D)$-measurable} if $A$ is a Borel set and, for every component $C$ of the Weierstrass locus $\Wloc(D)$, we have either
			\[
				C \subseteq A \qquad \text{or} \qquad C \subseteq \Gamma \setminus A.
			\]
			Let $\cA = \cA(D)$ denote the $\sigma$-algebra of $\Wloc(D)$-measurable subsets of $\Gamma$.
		\end{defi}
		
		In other words, given a Weierstrass locus $\Wloc(D) \subseteq \Gamma$, we can construct the quotient map $\pi \colon \Gamma \to \Gamma_0$ in which each component $C_i \subseteq \Wloc(D)$ is contracted to a single point. Then, the $\Wloc(D)$-measurable sets of $\Gamma$ are the preimages of Borel sets of $\Gamma_0$. If the divisor $D$ is \wfinite{}, then all Borel sets in $\Gamma$ are $L_W(D)$-measurable.
		
		\begin{defi}[Weierstrass measure] \label{def:weight_measure}
			Keeping the same notation as above, let $D$ be an effective divisor of rank $r$ on $\Gamma$, and let $\cA$ denote the $\sigma$-algebra of $\Wloc(D)$-measurable subsets of $\Gamma$. We define the {\em Weierstrass measure} $\hatweight$ as the ``weighted counting measure'' on $\Gamma$ whose atoms are the connected components in the Weierstrass locus $\Wloc(D)$. More precisely, $\hatweight$ is the measure on $(\Gamma, \cA)$ defined by
			\[ \hatweight(A) \coloneqq \sum_{C \subseteq A} \weight(C), \]
			where the sum is taken over components of $\Wloc(D)$ contained in $A$, and $\weight(C)$ is given by~\eqref{eq:w_weight_component}.
		\end{defi}

		We have the following description of the Weierstrass measure.
		\begin{thm} \label{thm:weight_measure}
	 		Keeping the same notation as above, for any closed connected $A \in \cA$, we have
			\begin{equation} \label{eq:weight_general}
				\hatweight(A) = \degrest{D}{A} + \mleft(g(A) - 1\mright) \, r - \sum_{\nu \in \partialout A} s_0^\nu(D).
			\end{equation}
		\end{thm}
		
\begin{figure}[h!]
			\centering
			\scalebox{.35}{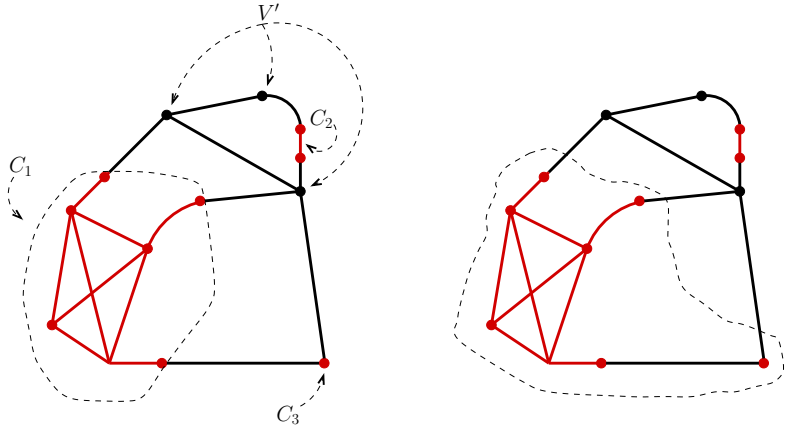}
			\caption{The part in red in the left figure is the (hypothetical) locus of Weierstrass points, and consists of three connected components. Red thickened points are on the boundary of the Weierstrass locus. Black vertices are those belonging to $V'$, that is, outside the Weierstrass locus. They are three in number. The right figure is an example of a set $A$ appearing in $\cA$. There is no vertex in $A$ outside the Weierstrass locus, so $m = 0$. There are two connected components of the Weierstrass locus in $A$, so $n = 2$. The subset $A \setminus |\widetilde \fA|$ consists of four intervals. This means $k = 4$.}
			\label{fig:example_weierstrass_measure}
		\end{figure}

		\begin{proof}
			Let $\fA = \{C_1, \ldots, C_n\}$ denote the set of components of $\Wloc(D)$ contained in $A$. Let $G = (V, E)$ be a model for $\Gamma$ whose vertex set $V$ contains the support of $D$, and let $V \cap (A \setminus L_W(D)) = \{v_1, \ldots, v_m\}$ denote the set of non-Weierstrass vertices in $A$.
			For each such vertex $v_i$, let $C_{n + i} = \{v_i\}$ denote the corresponding singleton, and let $\widetilde\fA$ denote the union 
			\[
				\widetilde\fA = \fA \cup \mleft\{ \{ v_1 \}, \ldots, \{ v_m \} \mright\} = \{ C_1, \ldots, C_n, C_{n + 1}, \ldots, C_{\tilde n} \} 
				\qquad \text{where} \qquad \tilde n = n + m.
			\]
			Finally, let $|\widetilde \fA| = \bigcup_{i = 1}^{\tilde n} C_i $ be the underlying subset of $ \Gamma$. Note that $|\widetilde \fA| \subseteq A$, and $A \setminus |\widetilde \fA|$ consists of a union of finitely many open intervals; let $k$ denote their number.
			
			Let $V' \coloneqq V \setminus \Wloc(D)$, as in Figure~\ref{fig:example_weierstrass_measure}.
			For each $v \in V'$, we have $D_v(v) = r$, so
			\[
				\weight(\{v\}) = D(v) - r - \sum_{\nu \in \T_v(\Gamma)} s_0^\nu(D) = D_v(v) - r = 0.
			\]
			Thus, the ``components'' $C_{n + i} = \{ v_i\}$ inside $\widetilde \fA \setminus \fA$ do not contribute to the total weight, so it suffices to show that
			$
				\displaystyle \sum_{C_i \in \widetilde \fA} \weight(C_i)
			$
			satisfies~\eqref{eq:weight_general}.
			
			From Definition~\ref{def:w_weight_component_intro}, we have
			\begin{equation} \label{eq:3}
				\sum_{i = 1}^{\tilde n} \weight(C_i) = \sum_{i = 1}^{\tilde n} \degrest{D}{C_i} + r \sum_{i = 1}^{\tilde n}(g(C_i) - 1) - \sum_{i = 1}^{\tilde n} \mleft(\sum_{\nu \in \partialout C_i} s_0^{\nu}(D) \mright).
			\end{equation}
			
			We treat separately the three terms appearing on the right-hand side of~\eqref{eq:3}. The first term $\sum_i \degrest{D}{C_i}$ is equal to $\degrest{D}{A}$, since the vertex set $V$ was chosen to contain the support of $D$.
			
			For the second term, we apply the identity
			\[
				\degrest{K}{B} = 2 g(B) - 2 + \outval(B)
				\qquad \text{for} \quad B \subseteq \Gamma \quad \text{closed and connected}
			\]
			(see Lemma~\ref{lem:canonical_sum}) twice to obtain
			\begin{align*}
				r \sum_{i = 1}^{\tilde n}(g(C_i) - 1) &= \frac r 2 \sum_{i = 1}^{\tilde n} \mleft(\degrest{K}{C_i} - \outval(C_i)\mright) \\
				&= \frac r 2 \mleft(\degrest{K}{A} - \outval(A) - 2 k\mright) = r \mleft(g(A) - 1\mright) - rk,
			\end{align*}
			where $k$, we recall, denotes the numbers of edges of $\Gamma \setminus |\widetilde \fA|$ whose endpoints are both in $|\widetilde \fA|$. Here, where $\outval(B) \coloneqq \abs{\partialout B}$ is the number of outgoing branches from $B$.
			
			For the third term, the collection of all tangent directions $\bigcup_{C_i \in \widetilde \fA} \mleft\{\nu \in \partialout C_i\mright\}$ can be partitioned into ``paired'' directions, if following $\nu$ leads to another component in $\widetilde \fA$, and ``unpaired'' directions, if following $\nu$ leads out of $A$.
			For any paired tangent direction $\nu \in \partialout C_i$, there is a matching opposite direction $\overline \nu \in \partialout C_j$ (see Section~\ref{subsec:basic_notation}) and their minimum slopes satisfy $s_0^{\nu}(D) + s_0^{\overline \nu}(D) = - r$.
			For any unpaired tangent direction $\nu \in \partialout C_i$, the minimum slope $s_0^\nu(D)$ is equal to $s_0^{\nu'}(D)$ for some parallel tangent direction $\nu' \in \partialout A$. Moreover, this gives a bijection between $\partialout A$ and the unpaired tangent directions.
			Using this, we have
			\begin{align*}
				\sum_{i = 1}^{\tilde n} \mleft( \sum_{\nu \in \partialout C_i} s_0^{\nu}(D) \mright) &= \sum_{\text{unpaired }\nu} s_0^\nu(D) + \sum_{\text{paired }\nu} s_0^\nu(D) \\
				&= \sum_{\nu \in \partialout A} s_0^\nu(D) + \sum_{\ell = 1}^k \mleft(s_0^{\nu_\ell}(D) + s_0^{\overline \nu_\ell}(D)\mright) = \sum_{\nu \in \partialout A} s_0^\nu(D) - rk.
			\end{align*}
			Combining the above identities shows that $\hatweight(A)$ satisfies~\eqref{eq:weight_general}. 
		\end{proof}
			
		\begin{remark}
			For a closed subset $A \in \cA$ with a finite number of connected components, the weight $\hatweight(A)$ can be expressed equivalently as
			\[
				\hatweight(A) = \degrest{D}{A} + \mleft(g(A) - c(A)\mright) \, r - \sum_{\nu \in \partialout A} s_0^\nu(D)
			\]
			where $c(A) = h_0(A)$ denotes the number of connected components of $A$. Note that $g(A) = h_1(A)$, so that in terms of the Euler characteristic $\chi$, the middle term is $- r \cdot \chi(A)$.
		\end{remark}
	
		The following result can be obtained by the same method.
	
		\begin{thm} \label{thm:number_w_points_more_general}
			Let $U$ be a connected open subset of $\Gamma$ that is $\Wloc(D)$-measurable. The Weierstrass weight $\hatweight(U)$ can be recovered from the slopes around the incoming branches as the sum
			\[
				\hatweight(U) = \degrest{D}{U} + \mleft(g(U) - 1\mright) \, r + \sum_{\nu \in \partialin U} s_{\mathrm{max}}^\nu(D)
			\]
			where $\partialin U$ denotes the set of incoming unit tangent vectors from the boundary of $U$, and $s_{\mathrm{max}}^\nu(D)$ the maximum slope along the incoming tangent vector $\nu$ of any rational function in $\Rat(D)$.
		\end{thm}
		
		Note that since $U$ is open and $\Wloc(D)$ is closed, every $\nu \in \partialin U$ is tangent to an open interval on $\Gamma$ that is outside $\Wloc(D)$ and thus $s^\nu_{\mathrm{max}}(D) = s^\nu_0(D) + r$ (see Theorem~\ref{thm:consecutive_slopes_intro}).
		Thus we have
		\[ \hatweight(U) = \degrest{D}{U} + \mleft(g(U) - 1 + \mathrm{inval}(U)\mright) \, r + \sum_{\nu \in \partialin U} s_{0}^{\nu}(D). \]
		
		\begin{proof}[Proof of Theorem~\ref{thm:number_w_points_intro}]
			We apply Theorem~\ref{thm:weight_measure} with $A = \Gamma$. The statement about \wfinite{} divisors follows from the first statement and Remark~\ref{rem:multiplicity_singleton}.
		\end{proof}
		
		We end this section by raising a rather vague question about the geometric description of the tropical Weierstrass locus provided in Remark~\ref{rem:geometric_description_Weierstrass_locus}.
		
		\begin{question}
			It is possible to express the Weierstrass weights $\Wloc(C)$ using the geometric description of Remark~\ref{rem:geometric_description_Weierstrass_locus} in a meaningful way?
		\end{question}
	
	\subsection{Consequences} \label{subsec:consequences_theorem}
		
		We now provide some direct consequences of the above results, starting with the following remark.
		
		\begin{remark}
			Theorem~\ref{thm:number_w_points_intro} and~\cite[Thm.~A]{Richman18} together imply that a generic divisor $D$ of degree $d \geq g$ has a finite Weierstrass locus made up of $g \, (d - g + 1)$ points, all of weight one. 
			Indeed, the cardinality of $\Wloc(D)$ given by~\cite[Thm.~A]{Richman18} is $g \, (d - g + 1)$, whereas the total weight given by Theorem~\ref{thm:number_w_points_intro} is $d - r + rg$. But $r = d - g$ generically and in this case we have $g \, (d - g + 1) = g \, (r + 1) = d - r + rg$.
		\end{remark}
		
		
		\begin{cor} \label{cor:weight_genus_bound}
			Suppose $D$ is a divisor of rank $r$. For any closed, connected, $\Wloc(D)$-measurable subset $A \subseteq \Gamma$, we have
			\[
				\hatweight(A) \geq g(A) \, r,
			\]
			with equality if and only if $A$ does not intersect the tropical Weierstrass locus.
		\end{cor}
		
		\begin{proof}
			The displayed inequality is a direct consequence of Theorem~\ref{thm:weight_measure} and Proposition~\ref{prop:bound_red_degree}. If now $A$ does not intersect the tropical Weierstrass locus, then both sides are zero. Otherwise, it contains a connected component of the Weierstrass locus and the strict inequality follows from the remark right after Theorem~\ref{thm:positivity}.
		\end{proof}
		
		\begin{cor}[Theorem~\ref{thm:cycle_weierstrass_point}]
			Suppose that the rank $r$ of $D$ is at least one. Then, the complement of the Weierstrass locus $\Wloc(D)$ is a disjoint union of (open) metric trees. In other words, every cycle in $\Gamma$ intersects the tropical Weierstrass locus.
		\end{cor}
		
		\begin{proof}
			For the sake of a contradiction, suppose that $A$ is a cycle in $\Gamma$ disjoint from the Weierstrass locus $\Wloc(D)$. Then, $A$ is $\Wloc(D)$-measurable, and by definition (Definition~\ref{def:weight_measure}), $\hatweight(A) = 0$. However, Corollary~\ref{cor:weight_genus_bound} states that $\hatweight(A) \geq g(A) \, r > 0$, which gives a contradiction.
		\end{proof}
	
	\subsection{Special cases of weights}
		
		Here, we point out some special cases of the weight formula.
		
		\begin{enumerate}[(i)]
			\item If a divisor $D$ has rank $r = 0$, then $\hatweight(\Gamma) = d$.
			Suppose $D$ is effective in its linear equivalence class.
			For any tangent direction $\nu$ outside the Weierstrass locus, the slope set $\fS^\nu(D)$ contains a single slope, and this slope must be zero since $D$ is effective.
			Thus, a component $C$ of the Weierstrass locus has weight $\displaystyle \weight(C) = \degrest{D}{C}$.
			
			\item If the genus $g = 0$, then for any effective divisor $\hatweight(\Gamma) = d - r = 0$. (In general $0 \leq d - r \leq g$.) In particular, this implies that the Weierstrass locus $\Wloc(D)$ is empty.
			
			\item If the genus $g = 1$, then for a divisor $D$ of degree $d$, the total Weierstrass weight is $\hatweight(\Gamma) = d$.
			The Weierstrass locus consists of $d$ separate components, and it follows that every component has weight $1$.
			Regarding the number of components, this is because the Weierstrass points are those satisfying $[d \cdot x] = [D]$, and therefore the Weierstrass locus restricted to the central circle of $\Gamma$ form a coset of the subgroup of $d$-torsion points.
			
			\item If the rank satisfies $r = d - g$, then $\hatweight(\Gamma) = d - r + rg = g \, (r + 1)$.
			In particular, this holds for a generic divisor class with degree $d \geq g$, and for every divisor with degree $d \geq 2g - 1$.
			
			\item If $D = K$ is the canonical divisor, then $d = 2g - 2$ and $r = g - 1$, so $\hatweight(\Gamma) = g^2 - 1$.
			See Section~\ref{sec:canonical_locus} below for more discussion of this case.
		\end{enumerate}
	
	\subsection{Combinatorial graphs}
	
		In this section we assume $\Gamma$ is a combinatorial graph. By this we mean $\Gamma$ admits a model $(G = (V, E), \ell)$ that has unit edge lengths. We assume the divisor $D$ is supported on the vertex set $V$.
		
		\begin{thm}
			Suppose $e = uv$ is an edge in $G$ whose interior $\mathring e$ is $\Wloc(D)$-measurable. Let $f_{uv}$ be a rational function that satisfies $\div(f_{uv}) = D_u - D_v$. Let $\nu$ be the unit tangent vector at $v$ along $e$, toward $u$.
			Then, the Weierstrass weight of the interior of $e$ is
			\[
				\hatweight(\mathring e) = r - \slope_\nu(f_{uv}).
			\]
		\end{thm}
		
		\begin{proof}
			Let $U \coloneqq \mathring e$. Since $\Wloc(D)$ is closed, we can take the open interval $U$ a little bit smaller so that its extremities are distinct from $u$ and $v$ and $U$ still contains the same components of $\Wloc(D)$. Theorem~\ref{thm:number_w_points_more_general} states that the sum of Weierstrass weights on $U = \mathring e$ is equal to 
			\[
				\hatweight(U) = \degrest{D}{U} + (g(U) - 1) \, r + \sum_{\nu \in \partialin U} s_\mathrm{max}^\nu(D).
			\]
			Since $D$ is supported on the vertex set, we have $\degrest{D}{U} = 0$, and we also have $g(U) = 0$. Thus, the expression simplifies to
			\[
				\hatweight(U) = - r + \mleft(s_\mathrm{max}^{v, \nu}(D) + s_\mathrm{max}^{u, \overline{\nu}}(D)\mright)
			\]
			where $\nu$ and $\overline{\nu}$ are tangent directions toward $u$ and $v$, respectively.
			If $f_u$ and $f_v$ satisfy
			\[
				\div(f_u) = D_u - D
				\qquad \text{and} \qquad
				\div(f_v) = D_v - D,
			\]
			then we have
			\[
				\slope_{\overline{\nu}} f_u(u) = s_0^{u, \overline \nu}(D) = s_\mathrm{max}^{u, \overline \nu}(D) - r
				\qquad\text{and}\qquad
				\slope_{\nu} f_v(v) = s_0^{v, \nu}(D) = s_\mathrm{max}^{v, \nu}(D) - r,
			\]
			and the relation $f_{uv} = f_u - f_v$ implies
			\begin{align*}
				\slope_{\nu} f_{uv}(u) &= \slope_{\nu}(f_u - f_v) (u)
				= - \mleft(s_\mathrm{max}^{u, \overline\nu}(D) - r\mright) - \mleft(s_\mathrm{max}^{v, \nu}(D) - r\mright) \\
				&= 2 r - \mleft(s_\mathrm{max}^{v, \nu}(D) + s_\mathrm{max}^{u, \overline \nu}(D)\mright).
			\end{align*}
			Note that the slope of $f_{uv}$ is constant along the interior of $e$, since the reduced divisors $D_u$ and $D_v$ are supported on vertices.
		\end{proof}

	\subsection{Canonical Weierstrass locus} \label{sec:canonical_locus}
			
		In this section we discuss the case of the canonical divisor on a metric graph.
		
		\subsubsection{Weierstrass weight}
			
			The weight formula~\eqref{eq:w_weight_component} for $\weight(C; D)$ may be specialized to the case of the canonical divisor $D = K$. We need the following lemma.
			
			\begin{lemma} \label{lem:canonical_sum}
				Let $K = \sum_{x \in \Gamma} (\val(x) - 2) \, (x)$ denote the canonical divisor of $\Gamma$, and let $A \subseteq \Gamma$ be a closed connected subset. Then
				\[
					\degrest{K}{A} = 2 g(A) - 2 + \outval(A),
				\]
				where $\outval(A)$ is the number of outgoing branches from $A$.
			\end{lemma}
			
			\begin{proof}
				The proof can be obtained by direct calculation using an adapted graph model. The details are omitted.
			\end{proof}
			
			By direct summation, this result generalizes to closed subsets with finitely many connected components.
			
			\begin{thm} \label{thm:canonical_weight}
				Suppose $\Gamma$ is a metric graph of genus $g$, and let $K$ be its canonical divisor. Then for any closed, connected subset $A \subseteq \Gamma$ that is $\Wloc(K)$-measurable, we have
				\[
					\hatweight(A; K) = (g + 1)(g(A) - 1) - \sum_{\nu \in \partialout A} (s_0^\nu(K) - 1).
				\]
			\end{thm}
			
			\begin{proof}
				As usual, let $\partialout A$ denote the set of outgoing tangent directions from $A$ in $\Gamma$, and let $\outval(A)$ denote its cardinality. From Theorem~\ref{thm:weight_measure} we have
				\[
					\hatweight(A; K) = \degrest{K}{A} + r \, (g(A) - 1) - \sum_{\nu \in \partialout A} s_0^\nu(K).
				\]
				The canonical divisor $K$ has rank $r = g - 1$.
				By Lemma~\ref{lem:canonical_sum}, on a closed connected set $B \subseteq \Gamma$, the degree $\degrest{K}{B}$ satisfies $\degrest{K}{B} = 2g(B) - 2 + \outval(B)$.
				Therefore,
				\begin{align*}
					\hatweight(A; K) &= \degrest{K}{A} + (g - 1)(g(A) - 1) - \sum_{\nu \in \partialout A} s_0^\nu(K) \\
					&= 2(g(A) - 1) + \outval(A) + (g - 1)(g(A) - 1) - \sum_{\nu \in \partialout A} s_0^\nu(K) \\
					&= (g + 1)(g(A) - 1) - \sum_{\nu \in \partialout A} (s_0^\nu(K) - 1),
				\end{align*}
				which concludes.
			\end{proof}
		
			If we repeat the same computation for the pluricanonical divisor $nK$, where $n \geq 2$, we find that
			\[
				\hatweight(A; nK) = (2n - 1) g (g(A) - 1) - \sum_{\nu \in \partialout A} (s_0^\nu(nK) - n).
			\]
			This next corollary to Theorem~\ref{thm:canonical_weight} is also a direct consequence of Theorem~\ref{thm:number_w_points_intro}.
			
			\begin{cor} \label{cor:canonical_total_weight}
				Suppose $\Gamma$ is a genus $g$ metric graph.
				
				\begin{enumerate}[(a)]
					\item 
					The sum of Weierstrass weights over all components of $\Wloc(K)$ is equal to $g^2 - 1$.
					
					\item 
					For any integer $n \geq 2$, the sum of Weierstrass weights over all components of $\Wloc(nK)$ is equal to $(2n - 1)g(g - 1)$.
				\end{enumerate}
			\end{cor}
			
			The next result is a special case of Corollary~\ref{cor:weight_genus_bound}.
			
			\begin{cor} \label{cor:weight_bound_canonical}
				Suppose $\Gamma$ is a metric graph of genus $g$.
				For any closed, connected, $\Wloc(K)$-measurable subset $A \subseteq \Gamma$, we have
				\[
					\hatweight(A) \geq g(A) \, (g - 1).
				\]
			\end{cor}

		\subsubsection{Edge symmetry}
			
			We now discuss properties of some specific Weierstrass points under some symmetry condition; see as well Section~\ref{subsec:graph_all_weierstrass}.
			
			\begin{defi} \label{def:edge_reflexive}
				An edge $e$ of a metric graph $\Gamma$ is {\em reflexive} if there is an automorphism $\sigma \colon \Gamma \to \Gamma$ such that $\sigma(e) = \overline{e}$, i.e., $\sigma$ reverses the direction of $e$.
			\end{defi}
			
			\begin{thm}
				Suppose $\Gamma$ is a metric graph of genus $g \geq 2$, and let $K$ denote the canonical divisor of $\Gamma$. Suppose $e$ is a reflexive edge in $\Gamma$.
				\begin{enumerate}[(a)]
					\item If $g$ is even, then the midpoint of $e$ is in the Weierstrass locus $\Wloc(K)$.
					
					\item If $g$ is odd, then the midpoint of $e$ is in the Weierstrass locus $\Wloc(nK)$ for any integer $n \geq 2$.
				\end{enumerate}
			\end{thm}
			
			\begin{proof}
				Let $x$ denote the midpoint of the reflexive edge $e$. The tangent space $\T_x(\Gamma)$ contains two directions $\{\nu_1, \nu_2\}$, and the reflexive assumption implies that the minimum slopes are equal in both directions, i.e., $s_0^{\nu_1}(K) = s_0^{\nu_2}(K)$.
				If $x$ is outside the Weierstrass locus, then the singleton $\{x\}$ is $\Wloc(K)$-measurable and we may apply the weight formula from Theorem~\ref{thm:canonical_weight},
				\[
					\hatweight(x; K) = (g + 1)(-1) - 2(s_0^{\nu_1}(K) - 1) \equiv g + 1 \mod 2.
				\]
				Hence if $g$ is even, then $\hatweight(x)$ is nonzero, which contradicts our assumption that $x$ is outside the Weierstrass locus. This proves part (a).
							
				Now consider $D = nK$ for $n\geq 2$. By a similar argument, if $x$ is outside the Weierstrass locus $\Wloc(nK)$, then its Weierstrass weight is
				\[
					\hatweight(x; nK) = (2n - 1)g(-1) - 2 (s_0^{\nu_1}(nK) - n) \equiv g \mod 2.
				\]
				If $g$ is odd, then the weight $\hatweight(x)$ is nonzero, which again gives a contradiction. This proves part (b).
			\end{proof}
		
	\subsection{Effective computation}
		
		We discuss a concrete way of determining the Weierstrass locus and weights in a given metric graph.
		
		Let $D$ be an effective divisor on $\Gamma$. There is an algorithmic way for determining all the minimum slopes of functions in $\Rat(D)$ along unit tangent vectors in $\Gamma$. This is based on chip-firing on metric graphs. More precisely, \cite{luo2011rank} gives a generalization of Dhar's burning algorithm for metric graphs, which allows us to test whether a divisor is $x$-reduced for any point $x \in \Gamma$ and eventually to compute reduced divisors. See Definition~2.10, Algorithm~2.13 and Theorem~2.15 in~\cite{luo2011rank}.
		
		We can extract the minimum slopes from this procedure. Let $x$ be a point of $\Gamma$ and $\nu \in \T_x(\Gamma)$ be a tangent direction at $x$. At step $i$ of the algorithm, following the notation of~\cite[Def.~2.10]{luo2011rank}, we count the number $n_i$ of indices $1 \leq j \leq J$ such that $Q^{(1)}_j$ contains a segment of $\Gamma$ starting at $x$ and supporting the direction $\nu$. The number $n_i$ is either zero or one and represents the number of chips that go through this segment toward the point $x$ at step $i$. We denote by $n$ the sum of the $n_i$'s. It is the total number of chips that are brought to $x$ by Dhar's algorithm via the branch supporting $\nu$. This means that $s^\nu_0 = -n$, which shows that the minimum slope on $\nu$ can be computed using Dhar's algorithm.
			
\section{Generalizations} \label{sec:weierstrass_locus_generalized_setup}
	
	In this section, we generalize the setting of the previous sections to the case of augmented metric graphs, that is, in the presence of genera associated to the vertices.
	
	Since the genus of a given vertex hides information about the geometry of the component, it turns out that there will be an ambiguity when talking about the Weierstrass locus of a divisor $D$. In fact, the right setup in this context is a divisor $D$ endowed with the data of a closed sub-semimodule $M$ of $\Rat(D)$, which plays the role of a (not necessarily complete) linear series on the augmented metric graph. In what follows, we will explain how the preceding definitions and results extend from divisors to semimodules in the more general setting of augmented metric graphs. We then introduce two special classes of semimodules, the {\em generic} semimodule associated to any divisor, and the {\em canonical} semimodule associated to the canonical divisor. We properly justify both of them using the framework of metrized complexes.
	
	In the following, we assume all semimodules are nonempty unless specified otherwise.
	
	\subsection{Weierstrass loci of semimodules and augmented metric graphs} \label{subsec:augmented_metric_graphs_semimodules}
		
		\subsubsection{Semimodules}
			
			Let $\Gamma$ be a metric graph, and $D$ a divisor of degree $d$ on $\Gamma$. 
			The set of functions $\Rat(D)$ naturally has the structure of a semimodule on the tropical semifield; we refer to~\cite{haase2012linear, AG22} for a discussion on this semimodule structure. Let $M$ be a sub-semimodule of $\Rat(D)$.
			We endow $\Rat(D)$ with the topology induced by $\| \cdot \|_\infty$, and say $M \subseteq \Rat(D)$ is {\em closed} if it is closed with respect to this topology.
			The following is a direct extension to semimodules of the rank of divisors on graphs introduced by Baker and Norine~\cite{baker2007riemann}.
			
			\begin{defi}[Divisorial rank]
				The {\em divisorial rank} or simply {\em rank} of $M \subseteq \Rat(D)$ (also called the {\em rank of $D$ with respect to $M$}) is the greatest integer $r$ such that for any effective divisor $E$ on $\Gamma$ of degree $r$, there exists a function $f \in M$ verifying $D + \div(f) \geq E$. It is denoted by $r(M, D)$.
			\end{defi}
			
			Note that for every divisor $D'$ such that $D \leq D'$, the inclusion $M \subseteq \Rat(D')$ holds, and therefore the divisorial rank $r(M, D')$ can be defined. In fact, as the following statement shows, the divisorial rank will only depend on the sub-semimodule $M \subseteq \Rat(\Gamma)$, and not on the divisor $D$, if we additionally assume that $M$ is closed. Therefore, we will work only with closed semimodules in the following, and will denote their rank simply by $r(M)$.
			Note that any (nonempty) semimodule has rank $r(M) \geq 0$. Also note that by definition, we have the immediate inequality $r(M) \leq r(D)$.
			
			\begin{prop} \label{prop:divisorial_rank_independent_divisor}
				The divisorial rank $r(M, D)$ of a closed semimodule $M \subseteq \Rat(D)$ depends only on $M$, as a sub-semimodule $M \subseteq \Rat(\Gamma)$.
			\end{prop}
			
			\begin{proof}
				First note that there is a unique minimal divisor $D_{\mathrm{min}}$ such that $M \subseteq \Rat(D_{\mathrm{min}})$, which is obtained by taking the (point-wise) minimum of all such divisors.
				
				Then, we denote $r(M, D)$ by $r$ and $r(M, D_{\mathrm{min}})$ by $r_0$. It is clear from the inequality $D_{\mathrm{min}} \leq D$ that the inequality $r_0 \leq r$ holds. We thus prove that $r_0 \geq r$. We choose a model $G = (V, E)$ such that the vertex set contains the support of $D$.
				
				First, we suppose that $E$ is an effective divisor of degree $r$ on $\Gamma$ whose support is disjoint from the support of $D$. By definition of $r$, there exists $f \in M$ such that $D + \div(f) \geq E$. Since $M \subseteq \Rat(D_{\mathrm{min}})$ and $D$ coincides with $D_{\mathrm{min}}$ outside $V$, it follows that $D_{\mathrm{min}} + \div(f) \geq E$.
				
				Now, let $E$ be an effective divisor of degree $r$ on $\Gamma$ whose support may intersect that of $D$. Let $(E_n)_n$ be a sequence of divisors of degree $r$ converging to $E$, such that for each $n$, the support of $E_n$ is disjoint from $V$. By what precedes, for each $n$, there exists a function $f_n \in M$ such that $D_{\mathrm{min}} + \div(f_n) \geq E_n$. Without loss of generality, assume that $f_n(x_0) = 0$ for some $x_0 \in \Gamma$. Thanks to the boundedness of the slopes of functions in $\Rat(D_{\mathrm{min}})$ (see~\cite[Lem.~1.8]{gathmann2008riemann}), we can assume that $(f_n)_n$ converges uniformly to a function $f$, which satisfies $D_{\mathrm{min}} + \div(f) \geq E$ at the limit. The limit function $f$ is in $M$ by assumption that $M$ is closed, which concludes the argument.
			\end{proof}
			
			\begin{remark}
				In essence, the above proof shows that the complement of the support of $D$ is a ``rank-determining set'' for the semimodule $M$ in the sense of~\cite{luo2011rank}.
			\end{remark}
			
			The notion of minimum slopes naturally extends to closed sub-semimodules.
			
			\begin{defi}[Slope sets and minimum slopes] \label{defi:minimum_slope_semimodule}
				Let $M \subseteq \Rat(D)$ be a closed sub-semimodule. Given a point $x \in \Gamma$ and a tangent direction $\nu \in \T_x(\Gamma)$, let $\fS^{\nu}(M)$ denote the {\em slope set}
				\[
					\fS^{\nu}(M) \coloneqq \mleft\{\slope_\nu f(x) \st f \in M\mright\}.
				\]
				Let $s_0^\nu(M)$ denote the {\em minimum slope} along $\nu$ of functions in $M$. More generally, let $s_j^\nu(M)$ denote the $(j + 1)$-smallest slope along $\nu$ of functions in $M$, i.e.,
				\[
					s_0^{\nu}(M) = \min \{\fS^\nu(M)\},
					\qquad
					s_j^\nu(M) = \min \{s \in \fS^\nu(M),\, s > s_{j - 1}^\nu\}.
				\]
				When the semimodule $M$ is clear from context, we will simply use $s_j^\nu$ to denote $s_j^\nu(M)$.
			\end{defi}
			
			The following result is obtained similarly to Proposition~\ref{prop:bound_red_degree}, see Remark~\ref{rem:fact_set_slopes_constant_closed_subsemimodule}; we omit the details.
			
			\begin{prop} \label{prop:bound_red_degree_semimodule}
				Suppose $M \subseteq \Rat(D)$ is a closed semimodule of divisorial rank $r$. Then for any closed, connected subset $A \subseteq \Gamma$, we have
				\[
					\degrest{D}{A} - \sum_{\nu \in \partialout A} s_0^\nu(M) \geq r.
				\]
			\end{prop}
		
		\subsubsection{Reduced divisors} \label{sec:reduced_semimodule}
			
			For closed $M \subseteq \Rat(D)$, there is a well-defined and well-behaved notion of $x$-reduced divisor denoted $D_x^M$ linearly equivalent to $D$ with respect to $M$ for every $x \in \Gamma$. Simply, we define $f_x \colon \Gamma \to \R$ by setting
			\[ f_x(p) \coloneqq \inf_{\substack{f \in M \\ f(x) = 0}} f(p) \quad \forall p \in \Gamma. \]
			Using the boundedness of slopes~\cite[Lem.~1.8]{gathmann2008riemann}, which makes all the functions of $M$ uniformly Lipschitz, the infimum in the definition above turns out to be a minimum, and $f_x$ is the uniform limit of a sequence of elements in $M$ (see~\cite[Prop.~7.3]{AG22}). Therefore, $f_x \in M$. We set $D_x^M \coloneqq D + \div(f_x)$. It follows from the definition that $\slope_\nu f_x(x) = s^\nu_0$ for all $\nu \in \T_x(\Gamma)$, and $D_x^M(x) = D(x) -\sum_{\nu \in \T_x(\Gamma)}s^\nu_0$. Therefore, the analogue of Lemma~\ref{lem:formula_reduced_divisor} holds.
		
		\subsubsection{Augmented metric graphs}
			
			An {\em augmented metric graph} is a metric graph $\Gamma$ endowed with a model $(G = (V, E), \ell)$ and a genus function $\g \colon V \to \Z_{\geq 0}$. The {\em genus} of $(\Gamma, \g)$, denoted by $g(\Gamma, \g)$ or simply $g$, is defined by
			\[ g(\Gamma, \g) \coloneqq g(\Gamma) + \sum_{v \in V} \g(v). \qedhere \]
			This terminology follows~\cite{amini2015lifting1}; ``vertex-weighted graph'' is used in other places. Augmented metric graphs arise from the semistable reduction of smooth proper curves over a valued field, when remembering the genera $\g(v) = g(X_v)$ of the components $X_v$, for $v \in V$.
			
			Note that any metric graph is naturally an augmented metric graph, by declaring the genus function to be the zero function. This means that what we will discuss below applies equally to the setting of non-augmented metric graphs.
		
		\subsubsection{Weierstrass locus}
			
			We now extend the notion of tropical Weierstrass locus to semimodules in the general setting of augmented metric graphs. Let $(\Gamma, \g)$ be an augmented metric graph. Let $D$ be a divisor on $\Gamma$ and $M$ be a closed sub-semimodule of $\Rat(D)$ of divisorial rank $r \leq r(D)$.
			
			\begin{defi}[Tropical Weierstrass locus of a closed semimodule] \label{def:weierstrass_points_semimodules}
				 The {\em tropical Weierstrass locus} of $M$, denoted by $\Wloc(M, D, \g)$ (or $\Wloc(M, \g)$ if $D$ is clear from the context), is the set of all points $x \in \Gamma$ that verify $D_x^M(x) + (\g(x) - 1) \, r > 0$.
				 
				 In the case the genus function $\g$ is zero, we lighten the notation and simply write $\Wloc(M, D)$, instead of $\Wloc(M, D, 0)$. We abbreviate $\Wloc(M, D)$ as $\Wloc(M)$ if $D$ is clear from context.
			\end{defi}
			
			The set $\Wloc(M, \g)$ is a closed subset of $\Gamma$ that can in general be infinite. Note that for every $x \in \Gamma$, we have $D_x^M(x) \geq r$ and therefore $D_x^M(x) + (\g(x) - 1) \, r \geq \g(x) \, r \geq 0$. In particular, if $\g(x) > 0$ and $r > 0$, then $x$ belongs to the tropical Weierstrass locus.
			
			We now associate an intrinsic weight to each connected component of the Weierstrass locus. The definition is analogous to Definition~\ref{def:w_weight_component_intro}; here it is adapted to semimodules and depends on the genus function. 
			
			Let $D$ be a divisor of degree $d$ on $\Gamma$, and let $M \subseteq \Rat(D)$ be a closed sub-semimodule of divisorial rank $r$. We use the notation of Definition~\ref{def:w_weight_component_intro} for $\degrest{D}{C}$, $g(C)$, and $\partialout C$; $s_0^\nu(M)$ is introduced in Definition~\ref{defi:minimum_slope_semimodule}.
			
			\begin{defi}[Intrinsic Weierstrass weight of a connected component] \label{def:w_weight_component_M}
				Let $C$ be a connected component of the tropical Weierstrass locus $\Wloc(M, \g)$. The {\em Weierstrass weight} of $C$, denoted by $\weight(C; M, D, \g)$, is defined by
				\begin{equation} \label{eq:w_weight_component_M}
					\weight(C; M, D, \g) \coloneqq \degrest{D}{C} + \mleft(g(C) + \sum_{x \in C} \g(x) - 1\mright) \, r - \sum_{\nu \in \partialout C} s_0^\nu(M).
				\end{equation}
				It is also denoted simply by $\weight(C; M, \g)$ or $\weight(C; \g)$ if $M$ and $D$ are understood from the context.
				
				In the case the genus function is zero, we use $\weight(C; M, D)$, $\weight(C; M)$ or $\weight(C)$ for $\weight(C; M, D, 0)$.
			\end{defi}
			
			Proposition~\ref{prop:tropical_weight_independent_linear_equivalence} on the independence of weights with respect to linear equivalence generalizes easily to this context. Furthermore, this quantity is well-defined because any connected component of $\Wloc(M, \g)$ is a metric graph, a result that adapts directly from Proposition~\ref{prop:well_defined}. As in the case of divisors (Proposition~\ref{prop:well_defined}), $\Wloc(M, \g)$ has a finite number of connected components. And since Theorem~\ref{thm:positivity} extends directly, we get $\weight(C; M, \g) > 0$. We denote
			by $g(C, \g)$ the sum $g(C) + \sum_{x \in C} \g(x)$, that is, the genus of $C$ in the augmented metric graph $(\Gamma, \g)$.
				
			\begin{defi}[Tropical Weierstrass divisor] \label{def:weierstrass_divisor_M}
				We say that $(M, D, \g)$ is {\em Weierstrass finite} or simply {\em \wfinite{}} if the tropical Weierstrass locus $\Wloc(M, D, \g)$ is finite. In this case, we define the {\em tropical Weierstrass divisor} $W(M, D, \g)$ as the effective divisor
				\[
					W(M, D, \g) \coloneqq \sum_{x \in \Wloc(M, \g)} \weight(x; M, D, \g) \, (x). 
				\]
				The tropical weight of $x$ verifies $\weight(x; M, D, \g) = D_x^M(x) + (\g(x) - 1) \, r$. We abbreviate $W(M, D, \g)$ as $W(M, \g)$ if $D$ is clear from the context. Note that the support $\abs{(W(M, \g))}$ of the tropical Weierstrass divisor is exactly the tropical Weierstrass locus $\Wloc(M, \g)$.
				
				In the case the genus function is zero, we simply use $W(M, D)$ or $W(M)$ for $W(M, D, 0)$.
			\end{defi}
			
			\begin{remark}
				If we set $M = \Rat(D)$, and if the genus function is $\g = 0$, then we recover the definitions given in Section~\ref{sec:proof_number_w_points} for a complete linear series on a non-augmented metric graph. Namely,
				\begin{enumerate}[(i)]
					\item For every $x \in \Gamma$, we have $D_x^{\Rat(D)} = D_x$.
					
					\item We have $\Wloc(\Rat(D), 0) = \Wloc(D)$.
					
					\item For every connected component $C$ of $\Wloc(\Rat(D), 0)$, we have
					\[ \weight(C; \Rat(D), 0) = \weight(C; D). \]
					
					\item $D$ is \wfinite{} if, and only if, $\Rat(D)$ is so. In this case, $W(\Rat(D), 0) = W(D)$. \qedhere
				\end{enumerate}
			\end{remark}
			
			The following proposition, a direct consequence of the definitions, states how the Weierstrass locus and Weierstrass weights on an augmented graph are related to the non-augmented definition.
			
			\begin{prop} \label{prop:wloci_augmented}
				If $M \subseteq \Rat(D)$ is a closed semimodule of rank $r$, then the following equalities hold.
				\begin{enumerate}[(a)]
					\item 
					$\Wloc(M, \g) = \Wloc(M) \cup |\g|$.
					
					\item For every connected component $C$ of $\Wloc(M, \g)$, we have
					\[ \weight(C; M, D, \g) = \weight(C; M, D) + r \sum_{x \in C} \g(x). \]
				\end{enumerate}
			\end{prop}
		
		\subsubsection{Total sum of Weierstrass weights}
			
			The following theorem is an analogue of Theorem~\ref{thm:consecutive_slopes_intro} for closed sub-semimodules of $\Rat(D)$, and is proved using a natural analogue of Lemma~\ref{lem:formula_reduced_divisor}, given in Section~\ref{sec:reduced_semimodule}. The only difference is that in the case of semimodules, sets of slopes are no longer necessarily made up of consecutive integers.
		
			\begin{thm} \label{thm:description_slopes_semimodules}
				Let $D$ be a divisor on $\Gamma$ and $M$ be a closed sub-semimodule of $\Rat(D)$ of divisorial rank $r$. We take a model for $(\Gamma, \g)$ such that the support of $D$ is made up of vertices. Let $x \in \Gamma$ be a point and $\nu \in \T_x(\Gamma)$ be a tangent direction.
				\begin{enumerate}[(a)]
					\item If the open interval $(x, x + \varepsilon \nu)$ is disjoint from the Weierstrass locus $\Wloc(M, \g)$, for $\varepsilon > 0$, then the set of slopes $\mleft\{\slope_{\nu} f(x) \st f \in M\mright\}$ consists of $r + 1$ consecutive integers $\{s^\nu_0, s^\nu_0 + 1, \dots, s^\nu_0 + r\}$.
					
					\item If the open interval $(x, x + \varepsilon \nu)$ is contained in the Weierstrass locus $\Wloc(M, \g)$, then the set of slopes $\mleft\{\slope_\nu f(x) \st f \in M\mright\}$ consists of integers $\{s^\nu_0 < s^\nu_1 < \cdots < s^\nu_t\}$ with $t \geq r$ and $s^\nu_t - s^\nu_0 \geq r + 1$.
				\end{enumerate}
			\end{thm}
			
			Part (a) implies in particular that for any edge $e$ outside the Weierstrass locus of $M$, the number of slopes of functions on $e$ is $r + 1$ and these slopes are consecutive.
			
			As a corollary, following the same computation as in the case of a divisor, we get an analogue of Theorem~\ref{thm:weight_measure}.
			
			\begin{thm}[Sum of Weierstrass weights for an incomplete series on an augmented metric graph] \label{thm:number_w_points_M}
				Suppose $(\Gamma, \g)$ is a genus $g = g(\Gamma, \g)$ augmented metric graph, $D$ is a degree $d$ divisor, and $M \subseteq \Rat(D)$ is a closed semimodule of divisorial rank $r \geq 0$.
				
				Then, the total sum of weights associated to connected components of $\Wloc(M, \g)$ is equal to $d - r + rg$. In particular, if $M$ is \wfinite{}, then we have $\deg(W(M, \g)) = d - r + rg$.
				
				More generally, let $\cA$ denote the $\sigma$-algebra of $\Wloc(M, \g)$-measurable subsets of $\Gamma$ and $\hatweight$ the counting measure on $(\Gamma, \cA)$ associated to the weights $\weight(C; M, \g)$ given as above. Then, for any closed, connected $A \in \cA$, we have
				\begin{equation} \label{eq:weight_general_semimodule}
					\hatweight(A; M, \g) = \degrest{D}{A} + \mleft(g(A, \g) - 1\mright) \, r - \sum_{\nu \in \partialout A} s_0^\nu(M),
				\end{equation}
				where $g(A, \g)$ denotes $g(A) + \sum_{x \in A} \g(x)$.
			\end{thm}
			
			Theorem~\ref{thm:number_w_points_M} implies the following analogue of Theorem~\ref{thm:cycle_weierstrass_point}.
			
			\begin{thm} \label{thm:cycle_weierstrass_point_M}
				If the divisorial rank $r$ of $M$ is at least one, then every closed connected subset $A$ of $\Gamma$ with $g(A, \g) \geq 1$ contains a point of $\Wloc(M, \g)$.
			\end{thm}
			
			\begin{proof}
				Theorem~\ref{thm:number_w_points_M} and Proposition~\ref{prop:bound_red_degree_semimodule} imply that for any closed, connected, $\Wloc(M, \g)$-measurable subset $A \subseteq \Gamma$, we have
				\[
					\hatweight(A; M, \g) \geq g(A, \g) \, r,
				\]
				an analogue of Corollary~\ref{cor:weight_genus_bound} for closed semimodules. Then, the argument used in the proof of Theorem~\ref{thm:cycle_weierstrass_point} yields the result.
			\end{proof}
		
		\subsubsection{Coherence under inclusion of semimodules}
			
			We have the following coherence property for the Weierstrass loci and weights associated to semimodules.
			
			\begin{prop} \label{prop:coherence}
			Let $M \subseteq M' \subseteq \Rat(D)$ be two closed semimodules of rank $r$. Then, $\Wloc(M, \g) \subseteq \Wloc(M', \g)$ and any $\Wloc(M', \g)$-measurable subset $A$ of $\Gamma$ is $\Wloc(M, \g)$-measurable. Moreover, the equality $\hatweight(A; M, \g) = \hatweight(A; M', \g)$ holds.
			\end{prop}
			
			\begin{proof}
				Note that the inclusion $M \subseteq M'$ implies that we have $D^M_y(y) \leq D^{M'}_y(y)$ for every $y \in \Gamma$. This, in turn, implies that $\Wloc(M, \g) \subseteq \Wloc(M', \g)$. The claim that $A$ is $\Wloc(M, \g)$-measurable follows then, since $A$ is assumed to be $\Wloc(M', \g)$-measurable.
				
				To see that $\hatweight(A; M, \g) = \hatweight(A; M', \g)$, it suffices to show that $s_0^\nu(M) = s_0^\nu(M')$ for each $\nu \in \partialout A$.
				Suppose $\nu$ is such a tangent direction pointing out of $A$. By part (a) of Theorem~\ref{thm:description_slopes_semimodules}, there are exactly $r + 1$ consecutive slopes of functions $F \in M'$ along $\nu$. The same statement holds for $M$. Since $M\subseteq M'$, we infer that these slopes are the same. In particular, $s_0^\nu(M) = s_0^\nu(M')$, as desired.
			\end{proof}
			
			In the following two sections, we specialize the above constructions to two special families of closed semimodules $M$: the generic semimodule associated to any divisor $D$, and the canonical semimodule.
	
	\subsection{The generic semimodule associated to a divisor} \label{subsec:weierstrass_points_augmented_graphs_gen}
		
		Let $(\Gamma, \g)$ be an augmented metric graph. Denote by $\abs{\g}$ the support of $\g$.
		For any divisor $D$ on $\Gamma$, we define a closed semimodule $\Ratgen(D, \g) \subseteq \Rat(D)$.
		
		\begin{defi} \label{def:sub_semimodule_genus_function_gen}
			The {\em generic linear series} or {\em generic semimodule} $\Ratgen(D, \g)$ consists of all rational functions $f$ on $\Gamma$ such that for every $x \in \Gamma$, we have the inequality
			\[ D(x) + \div(f)(x) \geq \g(x). \qedhere \]
		\end{defi}
			
		Equivalently, we have the equality $\Ratgen(D, \g) = \Rat(D_0)$ for the divisor $D_0$ defined by $D_0(x) \coloneqq D(x) - \g(x)$, for every $x \in \Gamma$. (The claimed containment $\Ratgen(D, \g) \subseteq \Rat(D)$ is clear.)
		
		It follows that $\Ratgen(D, \g)$ is closed in the $\| \cdot \|_\infty$ topology of $\Rat(D)$.
		
		\begin{remark}
			The superscript ``gen'' stands for ``generic'' because, from the viewpoint of the degeneration of smooth projective curves, augmented metric graphs can be obtained from intermediate geometric objects called metrized complexes of curves. If this is the case, the above definition gives, precisely, the tropical part of the linear series of a divisor on the metrized complex in the case where the restriction of the divisor on every curve component of the metrized complex is generic. See Section~\ref{subsec:justification_definition_semimodule_genus_function} for more details.
		\end{remark}
		
		The following statement computes the divisorial rank of the generic semimodule associated to a divisor.
		
		\begin{prop} \label{rem:inequalities_ranks}
			Denote by $r$ the divisorial rank of the generic semimodule $\Ratgen(D, \g)$, let $D_0$ be the divisor defined by $D_0(x) \coloneqq D(x) - \g(x)$ for every $x \in \Gamma$, and let $r(D)$ and $r(D_0)$ denote the respective ranks of the two divisors $D$ and $D_0$ in $\Gamma$ without the genus function. We have the following (in)equalities.
			\begin{enumerate}[(a)]
				\item 
				$ r \leq r(D) $;
				
				\item 
				$ r = r(D_0) $.
			\end{enumerate}
		\end{prop}

		\begin{proof}
			(a) The inequality follows from the containment $\Ratgen(D, \g) \subseteq \Rat(D)$.
			
			(b) This follows from Proposition~\ref{prop:divisorial_rank_independent_divisor} applied to $M \coloneqq \Ratgen(D, \g) = \Rat(D_0)$.
		\end{proof}
		
		Now that we have a closed sub-semimodule $\Ratgen(D, \g)$ of $\Rat(D)$ with a well-known divisorial rank, we can apply the machinery developed above.
		
		\begin{defi}[Generic tropical Weierstrass weights and locus of a divisor] \label{def:weierstrass_locus_augmented_graph_gen}
			Keeping the above notation, let $D$ be a divisor on an augmented metric graph $(\Gamma, \g)$. The tropical Weierstrass locus, the Weierstrass weights, and the Weierstrass divisor (if it exists) are defined by plugging the semimodule $M \coloneqq \Ratgen(D, \g)$ into Definitions~\ref{def:weierstrass_points_semimodules}, \ref{def:w_weight_component_M} and~\ref{def:weierstrass_divisor_M}.
			
			To lighten the notation while stressing the choice of the generic semimodule and the dependence on $D$ and $\g$, we write:
			\begin{enumerate}[(i)]
				\item $\Wlocgen(D, \g)$ for $\Wloc \mleft(\Ratgen(D, \g), \g\mright)$;
				
				\item $\weightgen(C; D, \g)$ for $\weight \mleft(C; \Ratgen(D, \g), \g\mright)$; and
				
				\item $W^{^\generic}(D, \g)$ for $W \mleft(\Ratgen(D, \g), \g\mright)$.
			\end{enumerate}
			When $D$ is clear from context, we simply use $\weightgen(C; \g)$ for $\weightgen(C; D, \g)$.
		\end{defi}
		
		Note that when $\g$ is the zero function, we have the equality $\Ratgen(D, \g) = \Rat(D)$, and so the above definition recovers the one given in the previous sections for the Weierstrass divisor associated to a divisor.
		
		Proposition~\ref{prop:wloci_augmented} and a straightforward computation gives the following description of the generic Weierstrass locus.
		
		\begin{prop} \label{prop:weierstrass_locus_containement_augmented_graph}
			The following equalities hold:
			\begin{enumerate}[(a)]
				\item
				$\Wlocgen(D, \g) = \Wloc(D_0) \cup \abs \g$;
				
				\item
				$
					\weightgen(C; D, \g) = \weight \mleft(C; D_0\mright) + (r + 1) \, \sum_{x \in C} \g(x).
				$
			\end{enumerate}
		\end{prop}

		In the remainder of this section, we discuss the generic semimodule associated to the canonical divisor. We first recall the definition of the canonical divisor in the augmented setting.
		
		\begin{defi}[Canonical divisor on an augmented metric graph] \label{def:canonical_divisor_augmented_graph}
			Given an augmented metric graph $(\Gamma, \g)$, the {\em canonical divisor} $K$ on $(\Gamma, \g)$ is defined by
			\begin{equation} \label{eq:canonical_div_augmented}
				K(x) \coloneqq \val(x) - 2 + 2 \, \g(x)
			\end{equation}
			for each $x \in \Gamma$.
		\end{defi}
		
		\begin{remark} \label{rem:canonical_sum_augmented}
			In the context of augmented metric graphs, Lemma~\ref{lem:canonical_sum} becomes the following statement: for every closed connected subset $A \subseteq \Gamma$,
			\[
				\degrest{K}{A} = 2 g(A) - 2 + 2 \sum_{x \in A} \g(x) + \outval(A). \qedhere
			\]
		\end{remark}
		
		The following statement gives the rank of the semimodule $\Ratgen(K, \g)$, which is not $g - 1$, contrary to what one might expect. The reason for this discrepancy will become clear in Section~\ref{subsubsec:justification_augmented_graphs}, using metrized complexes.
		
		\begin{prop}[Rank of the generic semimodule $\Ratgen(K, \g)$] \label{prop:rank_canonical_divisor_generic_setting}
			If the genus function $\g$ is nontrivial, the semimodule $\Ratgen(K, \g)$ has rank $g - 2$.
		\end{prop}
		
		\begin{proof}
			The rank of $\Ratgen(K, \g)$ coincides with the rank of $K_0 \coloneqq K - \sum \g(x) \, (x)$ within the non-augmented metric graph $\Gamma$. Since the genus function is nontrivial, we have $\deg(K_0) = 2 \, g(\Gamma) - 2 + \sum_{x} \g(x) > 2 \, g(\Gamma) - 2$ with $g(\Gamma)$ the genus of the non-augmented metric graph, and so, by Riemann--Roch on $\Gamma$, we have $r(K_0) = \deg(K_0) - g(\Gamma) = g(\Gamma, \g) - 2$.
		\end{proof}
		
		In the next section, we define the canonical linear series for an augmented metric graph, and show it has the correct rank $g - 1$.
		
		\begin{example} \label{ex:augmented_cycle_one_point_gen_2}
			We compute the Weierstrass locus of the generic semimodule $\Ratgen(K, \g)$ on a cycle with one point of positive genus equal to two.
		
			Let $(\Gamma, \g)$ be the augmented metric graph where $\Gamma$ is the cycle of length one, parametrized by the interval $[0, \, 1]$, the single vertex $v$ coincides with the endpoints $v = 0 = 1$, and $\g(v) = 2$. The genus of this augmented metric graph is $g = 3$.
								
			We consider the canonical divisor $K$ and the associated generic semimodule $\Ratgen(K, \g)$, as defined in the present section (see Definition~\ref{def:sub_semimodule_genus_function_gen}). The rank is $r = g - 2 = 1$ according to Proposition~\ref{prop:rank_canonical_divisor_generic_setting}, and the total weight of the Weierstrass locus is $6$. The Weierstrass locus consists of the vertex $v$ and the point of coordinate $\frac 1 2$. It is easy to compute that the weights are $\weight^\generic(v; K, \g) = 5$ and $\weight^\generic \mleft(\frac 1 2; K, \g\mright) = 1$. Figure~\ref{fig:aug_cycle_gen_1} shows the augmented metric graph and its Weierstrass locus. A generalization for any value of $\g(v)$ is presented in Section~\ref{subsubsec:augmented_cycle_one_gen}. \qedhere
			
			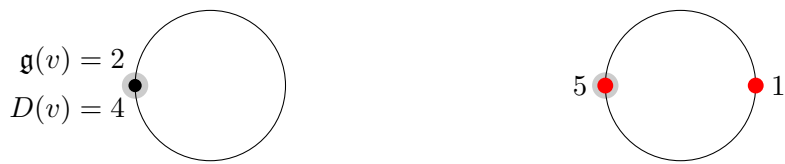
\begin{figure}[h!]
				\begin{minipage}{0.45\textwidth}
				\centering
				\begin{tikzpicture}[scale=.65]
					\coordinate (A) at (-1,0);
					\coordinate (B) at (0,0);
					
					\draw[radius=1.0] (B) circle;
					\fill (A) circle (2.5pt);
					\fill[opacity=0.2] (A) circle (5pt);
					
					\node[above left] at (A) {$\g(v) = 2$};
					\node[below left] at (A) {$D(v) = 4$};
				\end{tikzpicture}
				\end{minipage}
				\begin{minipage}{0.45\textwidth}
				\centering
				\begin{tikzpicture}[scale=.65]
					\coordinate (A) at (-1,0);
					\coordinate (B) at (0,0);
					\coordinate (C) at (1,0);
					
					\draw[radius=1.0] (B) circle;
					\fill[opacity=0.2] (A) circle (5pt);
					
					\foreach \c in {A,C} {
						\fill[color=red] (\c) circle (3pt);
					}
	
					\node[left=3pt] at (A) {$5$};
					\node[right=2pt] at (C) {$1$};
				\end{tikzpicture}
				\end{minipage}
				\caption{An augmented cycle graph with one point of genus two, the canonical divisor and its Weierstrass locus $\Wlocgen(D, \g)$.}
				\label{fig:aug_cycle_gen_1}
			\end{figure}
		\end{example}
		
	\subsection{The canonical linear series on an augmented metric graph} \label{subsec:weierstrass_points_augmented_graphs_can}
	
		Consider the augmented metric graph $(\Gamma, \g)$ and its canonical divisor $K$, defined by $K(x) = \val(x) - 2 + 2 \, \g(x)$ for each $x \in \Gamma$. In this section, we define the linear series $\KRat(\g)$ associated to $K$, that we call the {\em canonical linear series} or {\em canonical semimodule}.
		
		\begin{defi} \label{def:sub_semimodule_genus_function_can}
			We define the {\em canonical semimodule} $\KRat(\g)$ as the set of all functions $f \in \Rat(\Gamma)$ that verify the following conditions:
					
			\begin{enumerate}
				\item For every $x \in \Gamma$, we have $K(x) + \div(f)(x) \geq \g(x) - 1$.
				
				\item If $x$ has a tangent direction $\nu \in \T_x(\Gamma)$ such that $\slope_\nu f(x) \leq 0$, then $K(x) + \div(f)(x) \geq \g(x)$.
				\qedhere
			\end{enumerate}
		\end{defi}
		
		The following set of conditions is equivalent to that of Definition~\ref{def:sub_semimodule_genus_function_can}.
		
		\begin{enumerate}
			\item (local-minimum condition) If $x \in \Gamma$ is an {\em isolated local minimum} of $f$, i.e., $\slope_\nu f(x) \geq 1$ for every $\nu \in \T_x(\Gamma)$, then we impose
			$
				K(x) + \div(f)(x) \geq \g(x) - 1.
			$

			\item (generic condition) For all other points $x \in \Gamma$, we impose the stricter condition $K(x) + \div(f)(x) \geq \g(x)$.
		\end{enumerate}
		
		Note that according to the above definition, if a point $x$ has $\g(x) = 0$, then $x$ cannot be an isolated local minimum of $f \in \KRat(\g)$. Indeed, an isolated local minimum of $f$ satisfies $\div(f)(x) \leq - \val(x)$, and so $K(x) + \div(f)(x) \leq - 2$ assuming $\g(x) = 0$, which would violate both conditions. This means that, for any $x \in \Gamma$ and $f \in \KRat(\g)$, we have $K(x) + \div(f)(x) \geq 0$, which implies that $\KRat(\g)$ is a subset of $\Rat(K)$. (It is easy to see that it is in fact a semimodule.) This shows, moreover, that the above definition is equivalent to Definition~\ref{def:sub_semimodule_genus_function_gen} outside of the support of $\g$. Also note that we have the inclusion of semimodules $\Ratgen(K, \g) \subseteq \KRat(\g)$.
		
		\begin{remark}
			The definition of the canonical semimodule differs from the generic semimodule $\Ratgen(K, \g)$ given by Definition~\ref{def:sub_semimodule_genus_function_gen}. This is because the earlier definition, suitable for every divisor $D$ on $\Gamma$, assumed $D$ has ``generic support'' in the vertices with ``hidden genus.'' The canonical divisor, however, is not generic. Its specific properties suggest a distinct definition for the complete linear series of $K$.	
			The relevance of the above modification compared to Definition~\ref{def:sub_semimodule_genus_function_gen} will be further clarified in Section~\ref{subsec:justification_definition_semimodule_genus_function}.
		\end{remark}
		
		We have the following theorem, which justifies the name given to the linear series $\KRat(\g)$. Recall that $g = g(\Gamma, \g)$.
		
		\begin{thm} \label{thm:rank_semimodule_can}
			The divisorial rank of the semimodule $\KRat(\g)$ is $g - 1$.
		\end{thm}
		
		\begin{proof}
			The proof of this theorem will be given in Section~\ref{subsubsec:justification_canonical_divisor}.
		\end{proof}
		
		We have a closed sub-semimodule $\KRat(\g)$ of $\Rat(K)$ of divisorial rank $r = g - 1$, and we can apply the machinery developed for semimodules on augmented metric graphs.
				
		\begin{defi}[Canonical tropical Weierstrass weights and locus] \label{def:weierstrass_locus_augmented_graph_can}
			Keeping the above notation, the canonical tropical Weierstrass locus, the Weierstrass weights, and the Weierstrass divisor on an augmented metric graph are defined by plugging the semimodule $M \coloneqq \KRat(\g)$ into Definitions~\ref{def:weierstrass_points_semimodules}, \ref{def:w_weight_component_M} and~\ref{def:weierstrass_divisor_M}.
			
			To lighten the notation while stressing the choice of the canonical semimodule and the dependence on $\g$, we write:
			\begin{enumerate}[(i)]
				\item $\Wloc(K, \g)$ for $\Wloc(\KRat(\g), \g)$;
				
				\item $\weight(C; K, \g)$ for $\weight(C; \KRat(\g), \g)$; and
				
				\item $W(K, \g)$ for $W(\KRat(\g), \g)$. \qedhere
			\end{enumerate}
		\end{defi}
		
		\begin{example} \label{ex:augmented_cycle_one_point_can_2}
			In this example, we compute the canonical Weierstrass locus on an augmented cycle with a point of genus two. For the case of the generic Weierstrass locus associated to the same divisor $K$, see Example~\ref{ex:augmented_cycle_one_point_gen_2}.
			
			Let $(\Gamma, \g)$ be the augmented metric graph where $\Gamma$ is the cycle of length one, parametrized by the interval $[0, \, 1]$, the single vertex $v$ coincides with the endpoints $v = 0 = 1$, and $\g(v) = 2$. The genus of this augmented metric graph is $g = 3$.
			
			We consider the canonical divisor $K$ and the associated canonical semimodule $\KRat(\g)$, as defined in the present section (see Definition~\ref{def:sub_semimodule_genus_function_can}). The rank is $r = g - 1 = 2$ according to Theorem~\ref{thm:rank_semimodule_can}, and the total weight of the Weierstrass locus is $g^2 - 1 = 8$. The Weierstrass locus consists of the vertex $v$ and the points of coordinates $\frac 1 3$ and $\frac 2 3$. The Weierstrass weights are $\weight(v; K, \g) = 6$ and $\weight \mleft(\frac 1 3; K, \g\mright) = \weight \mleft(\frac 2 3; K, \g\mright) = 1$. Figure~\ref{fig:aug_cycle_1} shows the locus of Weierstrass points. A generalization for any value of $\g(v)$ is presented in Section~\ref{subsubsec:augmented_cycle_one_can}. \qedhere
			
			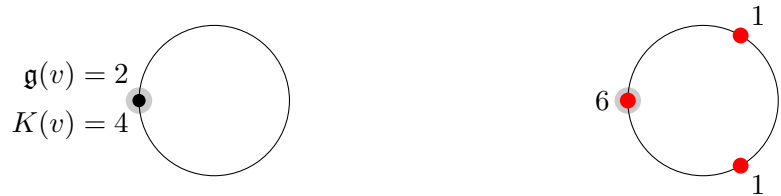
\begin{figure}[h!]
				\begin{minipage}{0.45\textwidth}
				\centering
				\begin{tikzpicture}[scale=.65]
					\coordinate (A) at (-1,0);
					\coordinate (B) at (0,0);
					
					\draw[radius=1.0] (B) circle;
					\fill (A) circle (3pt);
					\fill[opacity=0.2] (A) circle (5pt);
					
					\node[above left] at (A) {$\g(v) = 2$};
					\node[below left] at (A) {$K(v) = 4$};
				\end{tikzpicture}
				\end{minipage}
				\begin{minipage}{0.45\textwidth}
				\centering
				\begin{tikzpicture}[scale=.65]
					\coordinate (A) at (-1,0);
					\coordinate (B) at (0,0);
					\coordinate (C) at (60:1);
					\coordinate (D) at (300:1);
					
					\draw[radius=1.0] (B) circle;
					\fill[opacity=0.2] (A) circle (5pt);
					
					\foreach \c in {A,C,D} {
						\fill[color=red] (\c) circle (3pt);
					}
	
					\node[left=3pt] at (A) {$6$};
					\node[above right] at (C) {$1$};
					\node[below right] at (D) {$1$};
				\end{tikzpicture}
				\end{minipage}
				\caption{An augmented cycle graph, the canonical divisor and its Weierstrass locus $\Wloc(K, \g)$.}
				\label{fig:aug_cycle_1}
			\end{figure}
		\end{example}
		
		In the rest of the paper, when handling the canonical divisor $K$ on an augmented metric graph, the semimodule $\KRat(\g)$ will be preferred over $\Ratgen(K, \g)$, unless explicitly specified otherwise.
		
	\subsection{Justification of the definition of Weierstrass loci for augmented metric graphs, in the generic and canonical case} \label{subsec:justification_definition_semimodule_genus_function}
		In this section, we provide a justification for the definitions we gave in Sections~\ref{subsec:weierstrass_points_augmented_graphs_gen} and~\ref{subsec:weierstrass_points_augmented_graphs_can}. This will be through divisor theory on metrized complexes, that we recall first.
		
				
		\subsubsection{Divisor theory on a metrized complex of curves} \label{subsubsec:weierstrass_points_metrized_complex_curves}
			
			We fix $\kappa$ an algebraically closed field. A metrized complex of curves is, roughly speaking, the (metric realization of the) data of an augmented metric graph $(\Gamma, \g)$ endowed with a model $G = (V, E)$ and, for every $v \in V$, of a smooth, proper, connected, marked $\kappa$-curve $\cC_v$ of genus $\g(v)$ with marked points $A_v$ in bijection with the branches of $\Gamma$ incident to $v$. That is, a metrized complex of curves is a hybrid refinement of an augmented metric graph. For a full definition, see~\cite[Def.~2.17]{AB15}.
			
			Let $\fC$ be a metrized complex of curves. A {\em divisor} $\fD$ on $\fC$ is a formal sum with integer coefficients of a finite number of points in $\fC$. We denote its hybrid rank on $\fC$ by $r_{\fC}(\fD)$. By the forgetful projection map from $\fC$ to $\Gamma$, this gives rise to a divisor $D$ on $\Gamma$ of the same degree. Moreover, by restriction to each curve $\cC_v$, for $v \in V$, we get a divisor $\cD_v$ on $\cC_v$. A rational function $\mathfrak f$ on $\fC$ consists of a rational function $f$ on $\Gamma$ and, for every $v \in V$, a nonzero rational function $f_v$ on $\cC_v$. The space of such functions $\f = \mleft(f, f_v \st v \in V\mright)$ is denoted by $\Rat(\fC)$.
			
			Let now $\fC$ be a metrized complex of curves, with underlying metric graph $\Gamma$. Let $\fD$ be a divisor on $\fC$. We follow~\cite{AB15} and consider the linear series $\Rat(\fD, \fC)$ defined as the subset of $\Rat(\fC)$ consisting of all rational functions $\mathfrak f = (f \in \Rat(\Gamma); f_v \in \kappa(\cC_v), v \in V)$ on $\fC$ that verify $\fD + \div(\mathfrak f) \geq 0$. This means that the following two conditions are satisfied:
			\begin{enumerate}[(a)]
				\item \label{enumitem:tropical_divisor_effective} $D + \div(f)$ is effective on $\Gamma$; and
				
				\item \label{enumitem:geometric_divisor_effective} for each $v \in V$, the divisor $\cD_v - \sum_{\nu \in \T_v(\Gamma)} \slope_\nu f(v) \, (x^\nu_v) + \div(f_v)$ is effective on $\cC_v$.
			\end{enumerate}
			Here, $x^\nu_v$ is the marked point on $\cC_v$ that corresponds to $\nu$.
			
			\begin{defi}
				We define $\Rat^\trop(D, \fC)$ to be the subset of $\Rat(D)$ consisting of the tropical parts of all functions $\mathfrak f \in \Rat(\fD, \fC)$.
			\end{defi}
			
			We omit the proof of the following result.
			
			
			\begin{prop}
				$\Rat^\trop(D, \fC)$ is a closed sub-semimodule of $\Rat(D)$.
			\end{prop}
			
			We have the following comparison result, whose proof is direct from the definition of rank of divisors.
			
			\begin{prop} \label{prop:specialization_equality}
				Let $r(D, \fC)$ be the divisorial rank of the semimodule $\Rat^\trop(D, \fC)$. Then we have the inequality 
				\[ r_{\fC}(\fD) \leq r(D, \fC). \]
			\end{prop}
			
			
			The inequality in the above proposition can be strict in general. However, in some situations, e.g., for generic divisors on $\fC$ and for the canonical divisor, when the marked curves $(C_v, A_v)$, for all $v \in V$, are generic in their moduli, we have the equality, as we explain now.
		
		\subsubsection{The case of augmented metric graphs with generic divisors} \label{subsubsec:justification_augmented_graphs}
			
			Condition~\eqref{enumitem:geometric_divisor_effective} in the definition of $\Rat(\fD, \fC)$ in the previous section justifies Definition~\ref{def:sub_semimodule_genus_function_gen}. Indeed, take a rational function $f$ on $\Gamma$ such that for every $x \in \Gamma$, we have $D(x) + \div(f)(x) \geq \g(x)$. Assume that the augmented metric graph $(\Gamma, \g)$ comes from a metrized curve complex $\fC$. Let $v \in \Gamma$ be a point underlying a curve $\cC_v$. On the curve $\cC_v$, the divisor $\cD_v - \sum_{\nu \in \T_x(\Gamma)} \slope_\nu f(v) \, (x_v^\nu)$ has degree $\geq \g(v)$ by assumption. Therefore, by the Riemann--Roch theorem, its rank is non-negative, which is precisely Condition~\eqref{enumitem:geometric_divisor_effective} in the definition of $\Rat(\fD, \fC)$.
			Now, in the other direction, if $\cD_v$ is generic in the Picard group of $\cC_v$ of relevant degree, then the divisor $\cD_v - \sum_{\nu \in \T_v(\Gamma)} \slope_\nu f(v) \, (x^\nu_v)$ on $\cC_v$ appearing in the second condition has non-negative rank only if it has degree at least $\g(v)$. This means that Definition~\ref{def:sub_semimodule_genus_function_gen} is equivalent to the definition given for metrized complexes with a generic choice of divisors on components.
		
		\subsubsection{The case of canonical divisor in augmented metric graphs} \label{subsubsec:justification_canonical_divisor}
			
			We now justify Definition~\ref{def:sub_semimodule_genus_function_can} using the terminology of Section~\ref{subsubsec:weierstrass_points_metrized_complex_curves}, and also prove Theorem~\ref{thm:rank_semimodule_can}.
	
			Let $G = (V, E)$ be a model of $\Gamma$ whose vertex set contains all the points of valence different from two, and the support of $\g$. Let $\fC$ be a metrized complex with underlying augmented metric graph $(\Gamma, \g)$. Denote by $\fK$ a canonical divisor for $\fC$ given by the collection of divisors $\cK_{\cC_v} + A_v = \cK_{\cC_v} + \sum_{\nu \in \T_v(\Gamma)} (x^\nu_v)$ on $\cC_v$, where $\cK_{\cC_v}$ denotes a canonical divisor on $\cC_v$, i.e., $\cO(\cK_{\cC_v}) = \omega!_{\cC_v}$. The following claim justifies our definition of the canonical semimodule. We denote by $\Rat(\fK)^{\trop}$ the tropical part of $\Rat(\fK)$.
			
			\begin{prop} \label{prop:equality_seminodules_can}
				Keeping the same notation as above, we have $\KRat(\g) \subseteq \Rat(\fK)^{\trop}$. Moreover, if the markings $A_v$ on the curves $\cC_v$ are in general position (i.e., belong to a well-chosen Zariski-dense open subset of $\cC_v$), for all $v \in V$, then we have the equality $\KRat(\g) = \Rat(\fK)^{\trop}$.
			\end{prop}
			
			\begin{remark}
				This general position condition is the same as the one imposed in the work by Esteves and coauthors~\cite{EM02,ES07} in the special case of stable curves with two irreducible components, and stable curves in which any pair of components intersect. We will treat examples of augmented dipole graphs in Section~\ref{sec:further_examples_discussions}.
			\end{remark}
			
			\begin{proof}[Proof of Proposition~\ref{prop:equality_seminodules_can}]
				We first prove the inclusion $\KRat(\g) \subseteq \Rat(\fK)^{\trop}$.
				Consider an element $f \in \KRat(\g)$. We claim the existence of rational functions $f_v$ on $\cC_v$, for each $v \in V$, such that the collection $\{f, f_v, v \in V\}$ forms a rational function in $\Rat(\fK)$. This will prove the claim. Let $v \in V$, and consider the divisor $\cD$ on $\cC_v$ defined by $f$ as follows:
				\[ \cD_v \coloneqq \cK_{\cC_v} + \sum_{\nu \in \T_v(\Gamma)} (x^\nu_v) - \sum_{\nu \in \T_v(\Gamma)} \slope_\nu f(v) \, (x^\nu_v). \]
				Note that the degree of $\cD_v$ is precisely $K(v) + \div(f)(v)$. If the genus of $v$ is zero, then by the condition $K(v) + \div(f) \geq 0$, the degree of $\cD_v$ is non-negative and so there exists a rational function $f_v$ on $\cC_v$ such that $\cD_v+ \div(f_v) \geq 0$. If $\g(v) \geq 1$ and $v$ is not an isolated local minimum, then by the definition of $\KRat(\g)$, we have $\deg(\cD_v) \geq \g(v)$. By Riemann--Roch, this implies the existence of a function $f_v$ such that $\cD_v+ \div(f_v) \geq 0$.	Let $v \in \Gamma$ be a vertex of $\Gamma$ such that $\g(v) > 0$ and that is an isolated local minimum of $f$. In this case, by the definition of $\KRat(\g)$, we have $\deg(\cD_v)\geq \g(v) - 1$. The divisor $\cD_v$ can be rewritten as $\cK_{\cC_v} - E$, where
				\[
					E \coloneqq \sum_{\nu \in \T_v(\Gamma)} \mleft(\slope_\nu f(v) - 1\mright) \, (x^\nu_v) 
				\]
				is effective because $v$ is an isolated local minimum of $f$. The Riemann--Roch theorem on $\cC_v$, combined with the inequality $r(E) \geq 0$, thus yields
				\[ r(\cD_v) = r \mleft(\cK_{\cC_v} - E\mright) = r(E) + \deg(\cD_v) - \g(v) + 1 \geq 0. \]
				That is, there exists a function $f_v$ such that $\cD_v + \div(f_v) \geq 0$. The rational function $\f = \mleft(f, f_v \st v \in V\mright)$ on $\fC$ verifies $\fK + \div(\f) \geq 0$, as desired.
				
				We now prove the inclusion $\Rat(\fK)^{\trop} \subseteq \KRat(\g)$ provided that the markings $A_v$ on the curves $\cC_v$, for all $v \in V$, are generic. First, we observe that $\Rat(\fK)^{\trop} \subseteq \Rat(K)$. Combining this with the results we proved in Section~\ref{sec:slope_sets}, it follows that the slopes taken by functions in $\Rat(\fK)^\trop$ are bounded. Let $f$ be an element of $\Rat(\fK)^\trop$. We claim that under the general position assumption, we have $f \in \KRat(\g)$. Let $v$ be a vertex of $\Gamma$. Resuming the notation introduced above, we write $\cD_v$ for the divisor on $\cC_v$ induced by $f$, and write it in the form $\cD_v = \cK_{\cC_v} - E$.
				
				First consider the case where $v$ is an isolated local minimum of $f$. In this case, $E$ is an effective divisor. We need to show that $\deg(E) \leq \g(v) - 1$. By contradiction, assume that $\deg(E) \geq \g(v)$. If the points $x^\nu_v$, for $\nu \in \T_v(\Gamma)$, are in general position on $\cC_v$, we will get $r(\cD_v) \leq r(\cK_{\cC_v}) - \g(v) = -1$, which contradicts the assumption that $f \in \Rat(\fK)^\trop$.
				
				Consider the other case, where $v$ is not an isolated minimum. In this case, the divisor $E$ is not effective. We write $E = E_+ - E_-$ where $E_+$ and $E_-$ are the positive and negative parts of $E$, respectively. Note that $E_+$ and $E_-$ are effective and they have disjoint support. Since $E$ is not effective, $E_-$ is non-zero, and so by Riemann--Roch, we have
				\[ r(\cK_{\cC_v} + E_-) = 2 \g(v) - 2 + \deg(E_-) - \g(v) = \g(v) - 2 + \deg(E_-). \]
				Now, we write
				\[ \cD_v = \cK_{\cC_v} - E = \cK_{\cC_v} + E_- - E_+ \]
				and observe, by the general position assumption on the points of $A_v$, that
				\[ r(\cD_v) = \max \mleft\{-1, r(\cK_{\cC_v} + E_-) - \deg(E_+)\mright\}. \]
				Combining the two observations, we get
				\begin{align*}
					r(\cD_v) = \max \mleft\{-1, \g(v) - 2 + \deg(E_-) - \deg(E_+)\mright\} &= \max \mleft\{-1, \deg(\cK_{\cC_v} - E) - \g(v)\mright\} \\
					&= \max \mleft\{-1, \deg(\cD_v) - \g(v)\mright\}.
				\end{align*}
				If $\deg(\cD_v) < \g(v)$, we get $r(\cD_v) < 0$, which would be a contradiction to the assumption that $f \in \Rat(\fK)^\trop$. We conclude that $\deg(\cD_v) \geq \g(v)$, which leads to the inclusion $\Rat(\fK)^{\trop} \subseteq \KRat(\g)$.
			\end{proof}
			
 			We now show that $\KRat(\g)$ has the expected rank $g - 1$.
 			
			\begin{proof}[Proof of Theorem~\ref{thm:rank_semimodule_can}]
				We keep the above notation. We denote by $r$ the divisorial rank of $\KRat(\g)$.
				
				It will be enough to show that if the markings $A_v$ on the curves $\cC_v$, for all $v \in V$, are in general position, then we have $r_\fC(\fK) = r$. By Riemann--Roch for metrized complexes proved in~\cite{AB15}, we then obtain the equality $r = g - 1$, as desired.
				
				The inequality $r \geq r_\fC(\fK)$ follows from the case of equality $\Rat(\fK)^{\trop} = \KRat(\g)$ proved in the previous proposition, and by the definition of the rank in the metrized complex.
				
				It remains to show the inequality $r_\fC(\fK) \geq r$. Let $\cE$ be an effective divisor of degree $r$ on $\fC$, and let $E$ be the corresponding divisor on $\Gamma$. There exists a function $f \in \KRat(\g) = \Rat(\fK)^\trop$ such that $E + \div(f) \geq 0$. If $E$ has support outside the vertices of $\Gamma$, that is, $\cE$ is entirely supported at the interior of edges of $\Gamma$, then using the arguments we used in the first part of Proposition~\ref{prop:equality_seminodules_can}, we deduce the existence of rational functions $f_v$ on $\cC_v$, for all $v \in V$, such that the rational function $\f = \mleft(f, f_v \st v \in V\mright)$ on $\fC$ gives $\fK - \cE + \div(\f) \geq 0$, as desired.			
				
				Otherwise, if $\cE$ has support in some of the curves $\cC_v$, for $v \in V$, we write $E$ as a limit of effective divisors $E_n$, for $n \geq 0$, of the same degree with support outside the vertices of $\Gamma$, and find elements $f_n$ in $\Rat(\fK)^\trop = \KRat(\g)$ that verify $K - E_n + \div(f_n) \geq 0$. Going to a subsequence, and using the boundedness of the slopes in $\KRat(\g)$, we can suppose that all the $f_n$ have the same slopes along tangent directions at $v$, for each vertex $v \in V$. Moreover, changing the function $f \in \KRat(\g) = \Rat(\fK)^\trop$ under the constraint that $E + \div(f) \geq 0$ if necessary, we can assume furthermore that $f_n$ converges to $f$ as $n$ tends to infinity.
				
				Denote by $s^\nu_v$ the slope of the $f_n$ along the tangential direction $\nu \in \T_v(\Gamma)$. Let
				\[ \cD_v \coloneqq \cK_{\cC_v} + \sum_{\nu \in \T_v(\Gamma)} (x^\nu_v) - \sum_{\nu \in \T_v(\Gamma)} \slope_\nu f(v) \, (x^\nu_v) \]
				that we rewrite in the form
				\[ \cD_v = \cK_{\cC_v} + \sum_{\nu \in \T_v(\Gamma)} (1 - s^\nu_v) \, (x^\nu_v) + \sum_{\nu} m_\nu \, (x^\nu_v) \]
				with $m_\nu$ denoting the (weighted) number of points in the support of $E_n$ tending to $v$ through the tangential direction $\nu$. Note that we have $\sum_{\nu \in \T_v(\Gamma)} m_\nu = E(v)$. Let
				\[ \cD_v' \coloneqq \cK_{\cC_v} + \sum_{\nu \in \T_v(\Gamma)} (1 - s^\nu_v) \, (x^\nu_v). \]
				
				Two cases can happen. Either, some of the slopes $s^\nu_v$, for $\nu \in \T_v(\Gamma)$, are not positive, that is, $v$ is not an isolated local minimum of $f_n$. In this case, the divisor $\cD_v'$ has degree at least $\g(v)$, which implies that it has non-negative rank. Or, all the slopes $s^\nu_v$, for $\nu \in \T_v(\Gamma)$, are positive, that is, $v$ is an isolated local minimum of $f_n$ for all $n$. In this case, the divisor $\cD_v'$ has degree at least $\g(v) - 1$, and is the difference of $\cK_{\cC_v}$ and an effective divisor on $\cC_v$. So again, it has non-negative rank.
				
				In either case, we conclude that the divisor $\cD'_v = \cD_ v- \sum_{\nu \in \T_v(\Gamma)} m_\nu \, (x^\nu_v)$ has non-negative rank. Since the points $x^\nu_v$ are assumed to be in general position on $\cC_v$, it follows that the divisor $\cD_v$ has rank at least $E(v) = \sum_{\nu \in \T_v(\Gamma)} m_\nu$. This shows the existence of a rational function $f_v$ on $\cC_v$ such that $\cD_v - \cE_v + \div(f_v) \geq 0$, with $\cE_v$ being the part of $\cE$ supported in $\cC_v$. We conclude with the existence of a rational function $\f = \mleft(f, f_v \st v \in V\mright)$ that verifies $\fK - \cE + \div(\f) \geq 0$. This implies the inequality $r_\fC(\fK) \geq r$, and finishes the proof of our theorem.
			\end{proof}

\section{Tropicalization of Weierstrass points} \label{sec:tropicalization_Weierstrass_points}
		
	In this section we describe the tropicalization of the Weierstrass divisor of a line bundle on a smooth curve over a non-Archimedean field.
	
	Let $\K$ be an algebraically closed complete non-Archimedean field with a non-trivial valuation denoted by $\valuation$. Let $\vR$, $\fm$, and $\k = \vR/\fm$ be the valuation ring, the maximal ideal of $\vR$, and the residue field, respectively. We also denote by $|\cdot|$ the corresponding norm on $\K$, so that $\valuation(\cdot) = -\log|\cdot|$.
	
	Let $X$ be a smooth proper curve defined over $\K$. Let $\cD$ be a divisor of degree $d$ on $X$ and let $\cL = \cO(\cD)$ be the corresponding line bundle, with $\cO = \cO_X$, the structure sheaf of $X$. Denote by $\omega!_X$ the canonical line bundle on $X$.
	
	Let $H \subseteq H^0(X, \cL)$ be a space of sections of dimension $r + 1$ and denote by $\cW = \cW(\cD, H)$ the corresponding Weierstrass divisor. We assume that the gap sequence of $H$ is the standard sequence $0, 1, \dots, r$, that is, for a general point $x \in X(\K)$, the orders of vanishing of sections of $\cL$ in $H$ are $0, 1, \dots, r$. If $\K$ is of characteristic zero, this is automatic.
	
	\subsection{Tropicalization} \label{sec:tropicalizationB}
		
		We denote by $X^{\an}$ the Berkovich analytification of $X$. We assume familiarity with the theory of Berkovich analytic curves, and refer to~\cite[\S~4]{AB15} and~\cite[\S~5]{BPR} that contain what we need.
		
		A semistable vertex set for $X^{\an}$ is a finite set of type 2 points $V$ in $X^{\an}$ such that the complement $X^{\an} \setminus V$ is a disjoint union of finitely many open annuli and infinitely many open disks. A semistable model for $X$ is an integral proper relative curve $\fX$ over $R$ with generic fiber $\fX_\eta = X$ and special fiber $\fX_0$ that is reduced and has nodal singularities. Any irreducible component of the special fiber $\fX_0$ of a semistable model $\fX$ gives a valuation on $\K(X)$ and defines a point of type 2 in $X^{\an}$. The set $V$ of points in $X^{\an}$ associated to the irreducible components of $\fX_0$ is a semistable vertex set for $X^{\an}$. This process provides, in fact, a bijection between semistable vertex sets of $X^{\an}$ and semistable models of $\fX$ (see~\cite[Thm.~5.38]{BPR}). Moreover, each point of type 2 appears in a semistable vertex set.
		
		A semistable vertex set $V$ gives rise to a skeleton $\Gamma$ for $X^{\an}$, defined as the union in $X^{\an}$ of $V$ and the skeletons of the open annuli in $X^{\an} \setminus V$.
		The canonical metric on the skeletons of the open annuli gives the skeleton a metric graph structure, naturally embedded in $X^{\an}$.
		
		The underlying graph $G = (V, E)$ has vertex set $V$ and edge set $E$ in bijection with the set of open annuli in $X^{\an} \setminus V$. There is an edge between a pair of vertices $v$ and $u$ in $V$ for each open annulus whose closure contains the points $v$ and $u$. Moreover, the edge length function $\ell \colon E \to (0, +\infty)$ associates to each edge of $G$ the modulus of the corresponding annulus. Using the correspondence between semistable models and semistable vertex sets, the graph $G$ is identified with the dual graph of $\fX_0$, the special fiber of $\fX$, with vertices in bijection with the irreducible components of $\fX_0$, and edges in bijection with the nodes of $\fX_0$. There is an edge $e = uv$ in $G$ for each node that lies on the irreducible components associated to $u$ and $v$. The length of an edge corresponds to the singularity degree in $\fX$ of the corresponding node.
		
		For each point $x$ of type 2 in $X^{\an}$, the extension $\k(x)/\k$ is of transcendence degree one. We denote by $\rmC_x$ the corresponding smooth proper curve over $\k$. In a semistable model $\fX$ in which $x$ is in the vertex set, $\rmC_x$ is the normalization of the irreducible component in $\fX_0$ associated to $x$, and $\k(x)$ is the function field of this component.
		
		We denote by $\B_+$ the standard open ball in the Berkovich affine line $\ssA^{1, \an}$.
		The complement of $\Gamma$ in $X^{\an}$ is a disjoint union of open balls $B_\nu$ in bijection with $\nu \in \T_x(X^{\an}) \setminus \T_x(\Gamma)$ for all points $x$ of type 2 in $\Gamma$, each isomorphic to $\B_+$. 
		For a given ball $B_{\nu}$ in $X^{\an} \setminus \Gamma$, the corresponding point $x$ is the unique point in $\Gamma$ that lies in the closure of $B_{\nu}$.
		Denote by $p_x^\nu$ the point of $\rmC_{x}(\k)$ corresponding to $\nu \in \T_x(X^{\an})\setminus \T_x(\Gamma)$.
		
		Let $\Gamma$ be a metric graph skeleton of $X^{\an}$ with underlying graph $G = (V, E)$ and denote by $\tau \colon X^{\an} \to \Gamma$ the canonical retraction map. We call $\tau$ the tropicalization map. In the notation of the previous paragraph, the tropicalization map sends all the points in $B_\nu$ to the point $x$. The restriction of $\tau$ to $X(\K) \subset X^{\an}$ is compatible with the specialization map from the generic fiber $\fX_\eta$ to $\fX_0$, that is, a point specialized to a node is sent by $\tau$ to a point in the corresponding edge, and a point specialized to a smooth point of $\fX_0$ is sent by $\tau$ to the vertex of $G$ corresponding to this component. We get a tropicalization map $\tau_* \colon \Div(X) \to \Div(\Gamma)$ that sends a divisor $\cD = \sum_{x \in X(\K)} a_x (x)$ on $X$ to the divisor $\tau_*(\cD) = \sum_{x \in X(\K)} a_x (\tau(x))$.
		
		We denote by $\valuation_x \colon \K(X) \to \R \cup \{+\infty\}$ the valuation of a point $x \in X^{\an} \setminus X(\K)$ with $\valuation_x(f) = +\infty$ only if $f = 0$.
		The residue field of this valuation is denoted by $\k(x)$. We also denote by $|\cdot|_x = \exp(-\valuation_x)$ the corresponding norm.
		
		For each non-zero $f \in \K(X)$, we define the tropicalization of $f$, denoted $\Trop(f) \colon \Gamma \to \R$, as the map that sends each $x \in \Gamma \subset X^{\an}\setminus X(\K)$ to $\valuation_x(f)$. This induces a tropicalization map $\Trop \colon \K(X) \setminus \{0\} \to \Rat(\Gamma)$.
		
		For a vector subspace $H \subseteq \K(X)$, we call $M = \Trop(H \setminus \{0\})$ the tropicalization of $H$, and denote it, by a slight abuse of notation, by $\Trop(H)$.
		
		We define the genus function $\g$ on $X^{\an}$ to be the genus of $\rmC_x$ for a point of type 2, extended by zero everywhere else. The restriction of $\g$ to $\Gamma$ gives an augmented metric graph of genus $g$ equal to that of $X$. We denote by $K$ the canonical divisor of the augmented metric graph $(\Gamma, \g)$, with $K(x) = 2 \g(x) - 2 + \val(x)$ for all $x \in \Gamma$.
	
	\subsection{Reduction}
		
		For a type 2 point, the valuation $\valuation_x$ has the same value group as $\valuation$. For each nonzero $ f \in \K(X)$, choosing $a \in \K$ with $\valuation(a) = \valuation_x(f)$, we get that $\ssa^{-1}f$ has valuation $\valuation_x(\ssa^{-1} f) = 0$, and therefore gives an element in the residue field $\k(x)$ that we denote by $\sstildef_x$. We call this the reduction of $f$ at $x$, which is nonzero and defined only up to multiplication by a non-zero scalar in $\k$. For a vector subspace $H \subseteq \K(X)$ of dimension $r + 1$, denote by $\sstildeH_x \subseteq \k(x)$ the $\k$-vector subspace spanned by the reductions $\sstildef_x$ of elements $f \in H$ \cite[\S~4.4]{AB15}. By \cite[Lem.~4.3]{AB15}, $\sstildeH_x$ has dimension $r + 1$ over $\k$.

	\subsection{Slopes} \label{app:slopes}
		
		For a point $x$ in $\Gamma$ of type 2 in $X^{\an}$, each unit tangent direction $\nu \in \T_x(\Gamma)$ gives a point $\ssp_x^\nu \in \rmC_x(\k)$. By the slope formula~\cite{BPR}, for any non-zero $f \in \K(X)$, we have $\slope_\nu(\Trop(f)) = \ord_{\ssp_x^\nu}\mleft(\sstildef_x\mright)$. Moreover,
		\[ \tau_*(\div(f)) = \div(\Trop(f)). \]
		
		If $H \subseteq \K(X)$ is a $\K$-vector subspace of dimension $r + 1$, for any point $x \in \Gamma$ and unit tangent vector $\nu \in \T_x(\Gamma)$, we get a collection of integers $\slope_\nu(\Trop(f)) = \ord_{\ssp_x^\nu}(\sstildef_x)$, $f \in H$. Since $\sstildeH_x$ has dimension $r + 1$, this collection has size $r + 1$. This means that the collection of slopes $\slope_{\nu}(h)$, for $h \in M = \Trop(H)$, has size $r + 1$.
		In particular, Property $(\star)$ in Section~\ref{sec:weierstrass_locus_combinatorial_linear_series} is satisfied by $M = \Trop(H)$.
		
		For each unit tangent vector $\nu$, we order the slopes $\slope_{\nu}(h)$, for $h \in M$, in the form $s^\nu_0 < s^\nu_1 < \cdots < s_r^{\nu}$.
		Since elements of $\Trop(H)$ are piecewise affine linear, adding more points of $\Gamma$ to the semistable vertex set, we can suppose that the set of slopes
		\[ s^\nu_0 < s^\nu_1 < \cdots < s_r^{\nu} \]
		is constant in the interior of any edge of $G = (V, E)$ for parallel tangent directions $\nu$ at the point of the edge that point in the same direction.

	\subsection{Weierstrass divisor and Wronskian} \label{app:wronskianB}
		
		Let $X$ be a smooth proper curve defined over $\K$. Let $\cD$ be a divisor of degree $d$ on $X$ and let $\cL = \cO(\cD)$ be the corresponding line bundle, with $\cO = \cO_X$, the structure sheaf of $X$. Denote by $\omega!_X$ the canonical line bundle on $X$.
		
		Let $H \subseteq H^0(X, \cL)$ be a space of sections of dimension $r + 1$ and denote by $\cW = \cW(\cD, H)$ the corresponding Weierstrass divisor. The Weierstrass divisor $\cW$ is the divisor of a global section of the line bundle $\omega!_{X}^{\otimes \frac{r(r + 1)}2} \otimes \cL^{\otimes (r + 1)}$ called the \emph{Wronskian}. It is described as follows.
		
		In local coordinates, for any point $p \in X(\K)$, the local ring $\cO_{p}$ is a discrete valuation ring.
		We choose a uniformizer that we denote by $\ssunif_p$. We have $\cL_p \simeq \cO_p$ as an $\cO_p$-module. Taking the generator $\ssg_p = \ssunif_p^{\cD(p)}$ of $\cL_p$, each global section $f$ of $\cL$ can be written in the form $f = \ssf_{p}\ssg_p$ with $\ssf_{p} \in \cO_p$.
		
		We define the Hasse derivative of order $j$, for $j \in \Z_{\geq 0}$, on $\K[\ssunif_p]$ by
		\[ \Der{j} \ssunif^m = \binom{m}{j} \ssunif^{m - j} \qquad \textrm{ for } m > 0, \]
		and extend it by linearity to all $\K[\ssunif]$, and then to all $\K(\ssunif)$. Since the extension $\K(X)/\K(\ssunif)$ is separable, $\Der{j}$ is extended to $\K(X)$. Note that if $\K$ has characteristic zero, we can recursively define for any $j \geq 0$, the $j$-th derivative $\ssf_{p}^{(j)}$ by
		\[ \ssf_{p}^{(j)} = \frac {\mathrm d}{\mathrm d \ssunif_p} \ssf_p^{(j - 1)} \]
		with $\ssf_{p}^{(0)} = \ssf_{p}$. In this case, we have $j!\,\Der{j} \ssf_{p} = \ssf_{p}^{(j)}$.
		
		Let $\cF \coloneqq \{\ssf_0, \dots, f_r\}$ be a basis for $H \subseteq H^0(X, \cL)$, and for each $i$, write $f_i = \ssf_{i, p} \ssg_p$. Viewing $\Wron_{\cF}$ as a meromorphic section of $\omega!_X^{\otimes \frac {r(r + 1)}2}$, the stalk of the Wronskian $\Wron_{\mathcal F}$ at $p$ is given by
		\[ \Wron_{\cF, p} = \det \ssub{\mleft(\Der{j} \ssf_{i,p}\mright)}!_{0 \leq i, j \leq r} (\mathrm d \ssunif_p)^{\frac {r(r + 1)}2} \in \omega!_p^{\otimes \frac {r(r + 1)}2}. \]
		
		We have 
		\[ \cW = (r + 1) \cD + \div(\Wron_{\cF}). \]
		
		We note, without going into details, that the Wronskian $\Wron_{\cF}$ can also be defined without local coordinates, in terms of a filtration of the jet bundle and the diagonal embedding of $X$.
		
		
		
	
	\subsection{Slope formula for meromorphic differentials}
		
		We denote by $\|{\cdot}\|$ the K\"{a}hler norm introduced by Temkin in~\cite{Tem16} on the sheaf of differentials $\omega!_{X}$ that at any point $x \in X^{\an}$ associates to any section $\alpha$ of $\omega!_{X}$ the real number $\|\alpha\|_x$. For each positive integer $n$, the K\"{a}hler norm $\|\cdot\|$ induces a metric on $\omega!^{\otimes n}_{X}$, which, by an abuse of notation, we still denote by $\|\cdot\|$. Given a meromorphic section $\alpha$ of $\omega!^{\otimes n}_{X}$, the tropicalization of $\alpha$ denoted by $\Trop(\alpha)$ is the map
		\begin{align*}
			\Trop(\alpha) \colon \Gamma \to \R, \qquad x \mapsto -\log\|\alpha\|_x.
		\end{align*}
		The tropicalization of any meromorphic $n$-form on $X$ is a rational function on $\Gamma$, that is, $\Trop(\alpha) \in \Rat(\Gamma)$. Moreover, the following slope formula holds.
		
		\begin{lemma}[Slope formula for meromorphic differentials] \label{lem:slope_diff}
			For any meromorphic section $\alpha$ of $\omega!_X^{\otimes n}$, we have 
			\[ \tau_*(\div(\alpha)) = nK + \div(\Trop(\alpha)). \]
			Moreover, for any point $x \in \Gamma$ of type 2 and for any $\nu \in \T_x(\Gamma)$, we have $\slope_{\nu}(\Trop(\alpha)) = \ord_{\ssp_x^\nu}(\sstildealpha_x) + n$.
		\end{lemma}
		
		Here, $\sstildealpha_x$ is the reduction of $\alpha$ at $x$, and is a meromorphic form on $\rmC_{x}$, see~\cite[\S~2]{TT22}.
		
		\begin{proof}
			The case $n = 1$ is proved in~\cite{TT22}, see also~\cite[Thm.~2.6]{KRZ16}. The proof in general is similar and we only sketch it. The second statement, that $\slope_{\nu}(\Trop(\alpha)) = \ord_{\ssp_x^\nu}(\sstildealpha_x) + n$ for $\nu \in \T_x(\Gamma)$, is the analogue of \cite[Prop.~2.3.3]{TT22}, and both are special cases of~\cite[Lem.~3.3.2]{BT20}.
			
			Keeping the same notation as in Section~\ref{sec:tropicalizationB}, for each ball $B_\nu \simeq \B_+$, for $\nu \in \T_x(X^{\an}) \setminus \T_x(\Gamma)$, $\alpha$ can be written in the form $f \mathrm d \varT^n$ for a meromorphic function $f$ on $B_\nu$, with $\varT$ a parameter for the ball $B_\nu\simeq \B_+$, and the slope formula applied to $f$ implies
			\[ \sum_{a \in B_\nu(\K)} \ord_a(\alpha) = \ord_{\ssp_x^\nu} \sstildealpha_x. \]
			Write $\tau_*(\div(\alpha)) = \sum_{x \in \Gamma} a_x(x)$, for $a_x \in \Z$.
			For $x$ not of type 2, $a_x = 0$, and for $x$ of type 2, we have
			\begin{align*}
				a_x &= \sum_{\nu \in \T_x(X^{\an}) \setminus \T_x(\Gamma)} \mleft(\sum_{a \in B_\nu(\K)} \ord_{a}(\alpha)\mright) = \sum_{\nu \in \T_x(X^{\an}) \setminus \T_x(\Gamma)} \ord_{\ssp_x^\nu} \sstildealpha_x \\
				&= n(2 \g(x) - 2) - \sum_{\nu \in \T_x(\Gamma)} \ord_{\ssp_x^\nu} \sstildealpha_x \qquad \textrm{(since $\sstildealpha_x$ is a meromorphic $n$-form on $\rmC_x$)} \\
				&= n(2 \g(x) - 2) + n \val(x) - \sum_{\nu \in \T_x(\Gamma)} \slope_\nu(\Trop(\alpha)) = n K(x) + \ord_x(\Trop(\alpha)),
			\end{align*}
			as required.
		\end{proof}
	
	\subsection{Tropicalization of the Wronskian}
		
		Using the notation of Section~\ref{app:wronskianB}, let $F \coloneqq \Trop(\Wron_{\cF})$ be the tropicalization of the Wronskian $\Wron_{\cF}$, which is a meromorphic $(r(r + 1)/2)$-form. Let $W = \tau_*(\cW)$ be the tropicalization of $\cW$ to $\Gamma$. Let $D = \tau_*(\mathcal D)$. The following result is a direct consequence of Lemma~\ref{lem:slope_diff}, with $\alpha = \Wron_{\cF}$.
		
		\begin{thm} \label{thm:redW_general}
			Keeping the same notation as above, we have
			\[ W(x) = (r + 1) D(x) + \frac{r(r + 1)}2 K(x) - \sum_{\nu \in \T_x(\Gamma)} \slope_{\nu} F, \]
			where $F = \Trop(\Wron_{\cF})$. Furthermore, $\slope_{\nu} F = \frac{r(r + 1)}2 + \ord_{\ssp_x^\nu} \widetilde{\Wron_{\cF}}_{_x}$.
		\end{thm}
		
		Here, $\widetilde{\Wron_{\cF}}_{_x}$ denotes the reduction of $\Wron_{\cF}$ at $x$.
	
	\subsection{Wronskian of analytic functions on annuli} \label{sec:annuli}
		
		Let $\ssA^1 = \Spec(\K[\varT])$ and $\ssA^{1,\an}$ be its Berkovich analytification. Let $A(\rho)$ be the closed annulus in $\ssA^{1, \an}$ of center $0$ with outer radius one and inner radius $\rho \in (0, 1)$,
		\[ A(\rho) = \mleft\{x \in \ssA^{1, \an} \, \bigl| \, \rho \leq |\varT|_x \leq 1 \mright\}. \]
		Let $R(\rho)$ be the ring of analytic functions on $A(\rho)$. An analytic function $f$ on $A(\rho)$ admits a formal power series expansion
		\[ f = \sum_{n \in \Z} a_n \varT^n \]
		with $\lim_{n \to \pm \infty} |a_n| s^n = 0$ for all $s \in [\rho, 1]$.
		The skeleton of $A(\rho)$ is a closed interval, which can be identified with $I \coloneqq [0, -\log\rho]$: each point $q$ in this interval corresponds to the norm $|\cdot|_{\zeta_q}(f) = \sup_{n \in \Z} |a_n| \exp(-qn) = \max_{n \in \Z} |a_n| \exp(-qn)$ on any analytic function $f$ as above.
		The tropicalization of an analytic function $f$ is the function $\Trop(f)$ on the interval $I$ given by
		\[
			\Trop(f)(q) = \min \{\valuation(a_n) + nq \, \bigl| \, n \in \Z \} \qquad \forall q \in I.
		\]
		
		Let $\zeta = \zeta_0 \in A(\rho)$ be the boundary point corresponding to the extremity 0 of $I$, that is, $|f|_\zeta = \max_{n \in \Z} |a_n|$ on any analytic function $f$ as above. The reduction at $\zeta$ of an analytic function $f$ with $|f|_\zeta = 1$ is a Laurent polynomial
		\[ \sstildef_{\zeta} = \sum_{\substack{n \in \Z \\ |a_n| = 1}} \sstildea_n \vart^{n} \]
		where $\vart$ is the reduction of $\varT$ at $\zeta$ and $\sstildea_n \in \k$ is $\ssa_n$ modulo $\fm$. The slope of $\Trop(f)$ at 0 along the unit tangent direction $\nu \in \T_0(I)$ is the minimum exponent that appears in $\sstildef_{\zeta}$.
		
		\medskip
		
		Let $f_0, \dots, f_r$ be $r + 1$ $\K$-linearly independent analytic functions on $A(\rho)$ with
		\[ f_i = \sum_{n \in \Z} a_{i, n} \varT^n, \quad a_{i, n} \in \K \]
		the analytic expansion of $f_i$.
		Suppose that $|f_i|_\zeta = 1$ for all $i = 0, 1, \ldots, r$ and that $\Trop(f_0), \dots$, $\Trop(f_r)$ have slopes $s_0 < \cdots < s_r$ at $0$ along the unit tangent direction $\nu \in \T_0(I)$.
		This means the reduction $\sstildef_{i, \zeta}$ of $f_i$ at $\zeta$ has initial term $\vart^{s_i}$.
		
		Consider the analytic function on $A(\rho)$ defined by
		\[ h \coloneqq \det \ssub{\mleft(\Der{j} f_i\mright)}!_{i, j = 0}^r, \]
		where $\Der{j} f_i$ denotes the Hasse derivative of $f_i$.
		The following proposition describes the slope of $\Trop(h)$ at the point $0$ along the direction $\nu$.
		
		\begin{prop} \label{prop:slope_wronskian}
			Using the above notation, assume
			\begin{itemize}
				\item either, the residue field $\k$ is of characteristic zero,
				
				\item or, the sequence $s_0, \dots, s_r$ forms an interval, that is, $s_j = s_0 + j$ for all $j = 0, \dots, r$.
			\end{itemize}
			Then, we have
			\[ \slope_\nu(\Trop(h)) = s_0 + \dots + s_r - \frac{r(r + 1)}2. \]
		\end{prop}
		
		\begin{proof}
			The coefficient $a_{s_i, i}$ of $f_i$ has norm 1, so replacing $f_i$ with $a_{s_i, i}^{-1} f_i$, we can assume that $a_{i, s_i} = 1$. We write $f_i = \varT^{s_i} + f'_i$, so that all the coefficients of $\varT^n$ for $n \leq s_i$ that show up in the power series expansion of $f'_i$ have norm strictly smaller than one. We have
			\[ h = \det \ssub{\mleft(\Der{j} \varT^{s^\nu_i} + \Der{j} f'_i\mright)}!_{i, j = 0}^r = \det \ssub{\mleft(\binom{s_i}{j} \varT^{s^\nu_i - j} + \Der{j} f'_i\mright)}!_{i, j = 0}^r. \]
			Developing the determinant, we observe that
			\[ h = \det \ssub{\mleft(\binom{s_i}{j} \varT^{s^\nu_i - j}\mright)}!_{i,j = 0}^r + h' = C \, \varT^{s_0 + \dots + s_r - \frac{r(r + 1)}2} + h' \]
			with $h'$ an analytic function on $A(\rho)$ whose power series expansion in $\varT$ has the property that the coefficient of $\varT^{n}$ for $n \leq s_0 + \dots + s_{r} - \frac{r(r + 1)}2$ has norm strictly smaller than one, and
			\[ C = \det \ssub{\mleft(\binom{s_i}{j}\mright)}!_{i, j = 0}^r. \]
			Note that $C$ is a non-zero integer. If the residue field $\k$ has characteristic zero, then $C$ has norm one in $\K$, and the reduction of $h$ has minimum exponent $s_0 + \dots + s_r - \frac{r(r + 1)}2$.
			
			On the other hand, if $\k$ has positive characteristic, then by assumption $s_j = s_0 + j$ for all $j = 0, \dots, r$, and we use the identity for any $m \in \Z$
			\[
			 \det \begin{pmatrix} \binom{m}{0} & \binom{m}{1} & \dots & \binom{m}{r} \\
				 &&& \\
				 \binom{m + 1}{0} & \binom{m + 1}{1} & \dots & \binom{m + 1}{r} \\
				 &&& \\
				 \vdots & \vdots & \ddots & \vdots \\
				 &&& \\
				 \binom{m + r}{0} & \binom{m + r}{1} & \dots & \binom{m + r}{r} \\
			 \end{pmatrix} = 1
			\]
			(which is proved by induction on $r$ using repeated applications of Pascal's formula) to infer that $C = 1$ so again the reduction of $h$ at $\zeta$ has minimum exponent $s_0 + \dots + s_r - \frac{r(r + 1)}2$.
			
			We have shown that in either case, the slope of $\Trop(h)$ along $\nu \in \T_0(I)$ is $s_0 + \dots + s_r - \frac{r(r + 1)}2$.
		\end{proof}
	
	\subsection{Order of vanishing of the reduction of Wronskian}
		
		A consequence of Proposition~\ref{prop:slope_wronskian} is the following description of the slopes appearing in Theorem~\ref{thm:redW_general}.
		
		\begin{prop} \label{prop:reduction_inteval}
			Let $x$ be a point of type 2 in $\Gamma$ and $\nu \in \T_x(\Gamma)$. Denote by $s^\nu_0, \dots, s^{\nu}_r$ the sequence of slopes associated by tropicalization of $H$ to $\nu$. Assume
			
			\begin{itemize}
				\item either, the residue field $\k$ is of characteristic zero,
				
				\item or, the sequence $s^\nu_0, \dots, s^{\nu}_r$ forms an interval, that is, $s^\nu_j = s^\nu_0 + j$.
			\end{itemize}
			
			Then, we have
			\[ \ord_{\ssp_x^\nu} \widetilde{\Wron_{\cF}}_{_x} = s^\nu_0 + \dots + s^\nu_r - \frac{r(r + 1)}2. \]
		\end{prop}
		
		\begin{proof}
			Consider a segment $I$ in $\Gamma$ that contains $x$ with $\nu \in \T_x(I)$. $I$ is the skeleton of a closed annulus $A$ in $X^{\an}$ isomorphic to $A(\rho)$ for $0 < \rho < 1$. The elements of $\cF$ give a collection of analytic functions $f_0, \dots, f_r$ on $A \simeq A(\rho)$ that are $\K$-linearly independent, and whose reductions have the initial term $\vart^{s_i}$. The statement now follows by applying Proposition~\ref{prop:slope_wronskian} to the analytic function $h \coloneqq \det \ssub{\mleft(\Der{i} f_i\mright)}!_{i, j = 0}^r$, and observing that the restriction of $\Wron_{\cF}$ to $A$ coincides with $h (\mathrm d \varT)^{\frac{r(r + 1)}2}$ and that $\ord_{\ssp_x^\nu} \widetilde{\Wron_{\cF}}_{_x} = \slope_\nu(\Trop(h))$.
		\end{proof}
	
	\subsection{Reduction of Weierstrass points in equal characteristic zero}
		
		Assume the residue field $\k$ has characteristic zero. As in Section~\ref{app:slopes}, denote by $s_i^\nu$, $i = 0, \dots, r$, the slopes of functions of the form $\Trop(f) \in \Trop(H)$, for $f \in H$.
		
		\begin{thm} \label{thm:redW_zero}
			Let $W = \tau_*(\cW)$. We have
			\[ W(x) = (r + 1) D(x) + \frac{r(r + 1)}{2} K(x) - \sum_{\nu \in \T_x(\Gamma)} \sum_{i = 0}^r s_i^\nu. \]
		\end{thm}
		
		\begin{proof}
			This follows by combining Proposition~\ref{prop:reduction_inteval} with Theorem~\ref{thm:redW_general}.
		\end{proof}

\section{Tropical versus algebraic Weierstrass loci} \label{sec:weierstrass_locus_combinatorial_linear_series}
	
	In the first sections of this paper, we associated a Weierstrass locus to a fixed divisor $D$ on a metric graph, and then generalized this to a closed sub-semimodule of $M$ the space $\Rat(D)$ on an augmented metric graph. However, those semimodules that come from tropicalization verify an extra set of properties, in particular, the following important one (see Section~\ref{app:slopes}):
	\begin{enumerate}
		\item[($\star$)] 
		for each point $x$ in $\Gamma$ and any unit tangent direction $\nu \in \T_x(\Gamma)$, the set of slopes $\fS^\nu(M)$ taken by functions in $M$ has size $r + 1$.
	\end{enumerate}
	
	In this section, we associate to any pair $(M, D)$ consisting of a divisor $D$ and a closed sub-semimodule $M \subseteq \Rat(D)$ that verifies ($\star$) a refined notion of Weierstrass divisor. It is inspired from the formula given in Theorem~\ref{thm:redW_zero}, with the slopes being directly retrieved from $M$ 
	using property ($\star$). We then provide a comparison of this definition with that of Section~\ref{sec:weierstrass_locus_generalized_setup}. 
	Using this link, we prove the main result of this section, Theorem~\ref{thm:weierstrass_tropicalization_extended}, which relates tropical Weierstrass loci studied in the previous sections to tropicalization of Weierstrass divisors.
	We deduce then Theorem~\ref{thm:weierstrass_tropicalization_intro} as a special case of this result.
	
	\smallskip
	
	In the following, by an abuse of terminology, we refer to any pair $(M, D)$ as above as a combinatorial limit linear series (clls). The terminology is borrowed from~\cite{AG22}, however, the precise definition of combinatorial linear series requires more properties for the semimodule $M$. In our setting, we only need property ($\star$). The results can be thus applied more generally, in particular to the setting of tropical linear series developed in~\cite{JP22}. (No knowledge of \cite{AG22, JP22} is required in this paper.)

	\subsection{Weierstrass divisor of a combinatorial limit linear series}
		
		Let $D$ be a divisor on an augmented metric graph $(\Gamma, \g)$ and $M \subseteq \Rat(D)$ a closed sub-semimodule that verifies $(\star)$. Since $M \subseteq \Rat(D)$ is closed, we can apply the machinery of Section~\ref{sec:weierstrass_locus_generalized_setup}. This point of view on Weierstrass loci however results in a loss in information provided by the slopes of $M$, unless the Weierstrass locus is finite. The following definition relies on the knowledge of the slopes along edges of $G$ prescribed by $M$. 
		
		\begin{defi} \label{defi:clls_weierstrass}
			Suppose $D$ is a divisor of degree $d$ and $M$ is a closed sub-semimodule in $\Rat(D)$ such that $M$ verifies ($\star$). The {\em clls Weierstrass divisor} of $(M, D)$ is the divisor $W^\alg(M, D, \g)$ defined as
			\[ W^\alg(M, D, \g) \coloneqq \sum_{x \in \Gamma} \mu^\alg_\W(x) \, (x) \]
			where \emph{the clls Weierstrass weight} $\mu^\alg_\W(x)$ of $x$ is defined by
			\begin{equation} \label{eq:clls_weierstrass_divisor}
				\mu^\alg_\W(x) \coloneqq (r + 1) \, D(x) + \frac{r(r + 1)}{2} \, (\val(x) + 2 \g(x) - 2) - \sum_{\nu \in \T_x(\Gamma)} \, \sum_{j = 0}^r s^\nu_j(M).
			\end{equation}
			We write $W^\alg(M, D, \g)$ simply as $W^\alg(M, \g)$, the {\em clls Weierstrass divisor} of $M$, if $D$ is understood from the context.
			If the genus function is trivial, $\g = 0$, then we abbreviate $W^\alg(M, \g)$ to $W^\alg(M)$.
		\end{defi}
		
		Note that $W^\alg (M, \g)$ has finite support. Indeed, since elements in $M$ are piecewise affine linear and $M$ satisfies $(\star)$, by a compactness argument, we can find a model $G = (V, E)$ of $\Gamma$ such that $ s^\nu_0 < s^\nu_1 < \cdots < s_r^{\nu}$ is constant in the interior of any edge of $G = (V, E)$ for parallel unit tangent vectors $\nu$ based at points of the edge and pointing in the same direction (see, for example, \cite[Lem.~5.1]{JP22}). It follows that if $x \notin V$ and $x$ is outside the supports of $D$ and $\g$, then $\mu^\alg_\W(x) = 0$. Also note that the central term in the expression of $\mu^\alg_\W(x)$ above is equal to $\frac12{r(r + 1)} K(x)$, where $K$ is the canonical divisor on $(\Gamma, \g)$.
		
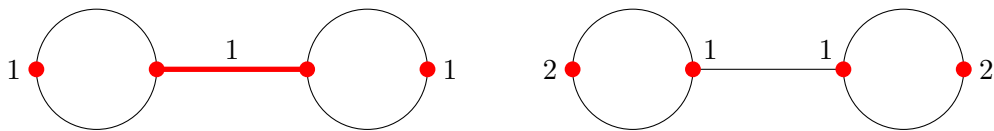
\begin{figure}[h!]
			\begin{minipage}{0.45\textwidth}
			\centering
			\begin{tikzpicture}[scale=.65]
				\coordinate (A) at (0,0);
				\coordinate (B) at (.8,0);
				\coordinate (C) at (2.8,0);
				\coordinate (D) at (3.6,0);
				\coordinate (E) at (-0.8,0);
				\coordinate (F) at (4.4,0);
				
				\draw[color=red,line width=1.2pt] (B) -- (C);
				\draw[radius=0.8] (A) circle;
				\draw[radius=0.8] (D) circle;
				
				\foreach \c in {B,C,E,F} {
					\fill[color=red] (\c) circle (3pt);
				}
				\node[above] at (1.8,0) {$1$};
				\node[left=2pt] at (E) {$1$};
				\node[right=2pt] at (F) {$1$};
			\end{tikzpicture}
			\end{minipage}
			\begin{minipage}{0.45\textwidth}
			\centering
			\begin{tikzpicture}[scale=.65]
				\coordinate (A) at (0,0);
				\coordinate (B) at (.8,0);
				\coordinate (C) at (2.8,0);
				\coordinate (D) at (3.6,0);
				\coordinate (E) at (-0.8,0);
				\coordinate (F) at (4.4,0);
				
				\draw (B) -- (C);
				\draw[radius=0.8] (A) circle;
				\draw[radius=0.8] (D) circle;
				
				\foreach \c in {B,C,E,F} {
					\fill[color=red] (\c) circle (3pt);
				}
				
				\node[left=2pt] at (E) {$2$};
				\node[above right] at (B) {$1$};
				\node[above left] at (C) {$1$};
				\node[right=2pt] at (F) {$2$};
			\end{tikzpicture}
			\end{minipage}
			\caption{The tropical Weierstrass locus $\Wloc(M)$ (left) and the clls Weierstrass divisor $W^\alg(M)$ (right) on the barbell graph.}
			\label{fig:barbell_clls_divisor}
		\end{figure}

		\begin{example}
			Consider the non-augmented barbell graph $\Gamma$ with edges of arbitrary length, see Figure~\ref{fig:barbell_clls}. This metric graph has genus two and the canonical divisor has rank one. We define a sub-semimodule $M \subseteq \Rat(K)$ of rank one on $\Gamma$ by prescribing the slopes $\{-1, 1\}$ on the middle edge and, for $i = 1, 2$, slopes $\{0, 1\}$ on both oriented edges $u_i v_i$. Then, we define $M$ as the set of all functions in $\Rat(K)$ that, along any unit tangent vector at a given point of $\Gamma$, take one of the two prescribed slopes. It is easy to see that $M$ is closed and verifies $(\star)$.
			
			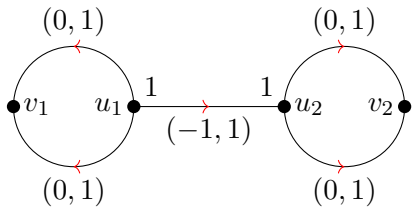
\begin{figure}[h!]
				\centering
				\begin{tikzpicture}[scale=.75]
					\coordinate (A) at (0,0);
					\coordinate (B) at (.8,0);
					\coordinate (C) at (2.8,0);
					\coordinate (D) at (3.6,0);
					\coordinate (E) at (-0.8,0);
					\coordinate (F) at (4.4,0);
					
					\begin{scope}[decoration={markings,mark=at position 0.5 with {\arrow[red]{>}}},
					] 
						\draw[postaction={decorate}] (B) -- (C);
						\draw[postaction={decorate}] (B) arc (0:180:0.8);
						\draw[postaction={decorate}] (B) arc (0:-180:0.8);
						\draw[postaction={decorate}] (C) arc (180:0:0.8);
						\draw[postaction={decorate}] (C) arc (-180:0:0.8);
					\end{scope}
					
					\foreach \c in {B,C,E,F} {
						\fill (\c) circle (2.5pt);
					}
					
					\node[right] at (E) {$v_1$};
					\node[left] at (B) {$u_1$};
					\node[right] at (C) {$u_2$};
					\node[left] at (F) {$v_2$};
					
					\node[above right] at (B) {$1$};
					\node[above left] at (C) {$1$};
					
					\node[below] at (1.8,0) {$\{-1, 1\}$};
					\node[above] at (0,0.8) {$\{0, 1\}$};
					\node[below] at (0,-0.8) {$\{0, 1\}$};
					\node[above] at (3.6,0.8) {$\{0, 1\}$};
					\node[below] at (3.6,-0.8) {$\{0, 1\}$};
				\end{tikzpicture}
				\caption{The barbell graph, the canonical divisor and the slope structure $\fS$.}
				\label{fig:barbell_clls}
			\end{figure}
			
			The clls Weierstrass divisor is
			\[ W^\alg(M) = (u_1) + (u_2) + 2 \, (v_1) + 2 \, (v_2) \]
			(see Figure~\ref{fig:barbell_clls_divisor}, right). For comparison, the tropical Weierstrass locus $\Wloc(M, \g)$ of $M$ with trivial genus function $\g = 0$ (as defined in Section~\ref{sec:weierstrass_locus_generalized_setup}), is shown on the same figure (left). It turns out to be identical to the tropical Weierstrass locus $\Wloc(K)$ of the complete linear series $\Rat(K)$ (see Example~\ref{ex:barbell_graph}).
		\end{example}
		
	\subsection{Effectivity} \label{sec:tangential_ramification}
		
		Unlike the tropical Weierstrass divisors defined earlier in this paper, the Weierstrass divisor defined in Definition~\ref{defi:clls_weierstrass} is not automatically effective. We can rewrite the Weierstrass weight as
		\begin{align*}
			\mu^\alg_\W(x) &= (r + 1) \mleft(D(x) + \frac{r}{2} \val(x) + (\g(x) - 1) \, r\mright) 
			 - \sum_{\nu \in \T_x(\Gamma)} \, \sum_{j = 0}^r \mleft(s^\nu_0 + j + \mleft(s^\nu_j - s^\nu_0 - j\mright)\mright) \\
			&= r(r + 1) \g(x) + \underbrace{(r + 1) \mleft(D^M_x(x) - r\mright)}_{\geq 0} - \underbrace{\sum_{\nu \in \T_x(\Gamma)} \, \sum_{j = 0}^r \mleft(s^\nu_j - s^\nu_0 - j\mright)}_{\geq 0}.
		\end{align*}
		
		We refer to the sequence
		\[
			\alpha_j^\nu(M) \coloneqq s_j^\nu(M) - s_0^\nu(M) - j, \qquad \text{for } j = 0, 1, \ldots, r
		\]
		as the {\em ramification sequence} of $M$ at $x$ along the tangential direction $\nu$. This sequence is non-decreasing.
		
		Let $\g$ be a genus function on $\Gamma$. The combinatorial limit linear series $M$ is called {\em $\g$-effective} if $W^\alg(M, \g)$ is effective. That is, for all $x \in \Gamma$,
		\begin{equation} \label{eq:inequality_g_effective_linear_series}
			r(r + 1) \g(x) + (r + 1) \mleft(D^M_x(x) - r\mright) \geq \sum_{\nu \in \T_x(\Gamma)} \, \sum_{j = 0}^r \alpha^\nu_j(M). \qedhere
		\end{equation}
		
		We say that $M \subseteq \Rat(\Gamma, \g)$ is {\em realizable} if there exists a smooth proper curve $X$ of genus $g$ over $\K$, a line bundle $\cL = \cO(\cD)$ of degree $d$ and a subspace $H \subseteq H^0(X, \cL)$ of rank $r$ such that $(\Gamma, \g)$ is a skeleton of $X^\mathrm{an}$, and $M = \mleft\{\tropic f \st f \in H \setminus \{0\}\mright\}$.
		If this happens over $\K$ of equicharacteristic zero, we say $M$ is realizable in equicharacteristic zero.
		
	
		\begin{prop} \label{prop:realizability_necessary_conditions}
			If $M$ is realizable in equicharacteristic zero, then the following hold:
			
			\begin{enumerate}[(i)]
				\item $W^\alg(M, \g)$ is effective, i.e., $M$ is $\g$-effective;
				
				\item the divisor of degree zero
				\[ W^\alg(M, \g) - (r + 1) \, D - \frac{r(r + 1)}{2} \, K = \sum_{x \in \Gamma} \mleft(\sum_{\nu \in \T_x(\Gamma)} \, \sum_{j = 0}^r s^\nu_j\mright) \, (x) \]
				is principal.
			\end{enumerate}
		\end{prop}
		
		\begin{proof}
			Both statements follow from Theorem~\ref{thm:redW_zero}.
		\end{proof} 
		
		\begin{example}
			Consider the non-augmented metric graph $\Gamma$ below and its canonical divisor $K$.
			We consider the following combinatorial limit linear series $M \subseteq \Rat(K)$. For each bridge edge oriented outwards (toward the adjacent circle), allow slopes $\{-1, 1, 3\}$. Divide each circle into three equal parts, in a way compatible with the position of the attachment points. On the two edges adjacent to the attachment points, allow slopes $\{0, 1, 2\}$ away from the attachment points, and on the remaining edges, allow slopes $\{-1, 0, 1\}$ (see Figure~\ref{fig:three_cycles}).
			
			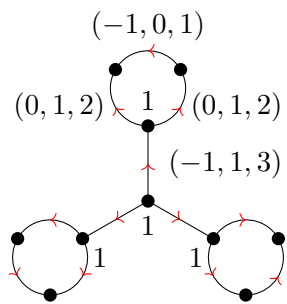
\begin{figure}[h!]
				\centering
				\begin{tikzpicture}[scale=.75]
					\coordinate (A) at (0,0);
					\coordinate (B) at (0,1);
					\coordinate (C) at (-0.87,-0.5);
					\coordinate (D) at (0.87,-0.5);
					\coordinate (Bc) at (0,1.5);
					\coordinate (Cc) at (-1.3,-0.75);
					\coordinate (Dc) at (1.3,-0.75);
					\coordinate (Bf) at ($(Bc) + (30:0.5)$);
					\coordinate (Bg) at ($(Bc) + (150:0.5)$);
					\coordinate (Cf) at ($(Cc) + (150:0.5)$);
					\coordinate (Cg) at ($(Cc) + (-90:0.5)$);
					\coordinate (Df) at ($(Dc) + (270:0.5)$);
					\coordinate (Dg) at ($(Dc) + (30:0.5)$);
					
					\draw[postaction=decorate,decoration={markings,mark=at position 0.5 with{\arrow[red]{>}},mark=at position 0.5 with {\node[right] {$\{0, 1, 2\}$};}}] (B) arc (-90:30:0.5);
					\draw[postaction=decorate,decoration={markings,mark=at position 0.5 with{\arrow[red]{>}},mark=at position 0.5 with {\node[left] {$\{0, 1, 2\}$};}}] (B) arc (270:150:0.5);
					\draw[postaction=decorate,decoration={markings,mark=at position 0.5 with{\arrow[red]{>}},mark=at position 0.5 with {\node[above] {$\{-1, 0, 1\}$};}}] (Bf) arc (30:150:0.5);
					
					\begin{scope}[decoration={markings,mark=at position 0.5 with {\arrow[red]{>}}}] 
						\draw[postaction={decorate}] (A) -- node[label=right:{$\{-1, 1, 3\}$}] {} (B);
						\draw[postaction={decorate}] (A) -- (C);
						\draw[postaction={decorate}] (A) -- (D);
						\draw[postaction={decorate}] (C) arc (30:150:0.5);
						\draw[postaction={decorate}] (C) arc (30:-90:0.5);
						\draw[postaction={decorate}] (Cf) arc (150:270:0.5);
						\draw[postaction={decorate}] (D) arc (150:30:0.5);
						\draw[postaction={decorate}] (D) arc (150:270:0.5);
						\draw[postaction={decorate}] (Df) arc (-90:30:0.5);
					\end{scope}
						
					\foreach \c in {A,B,C,D,Bf,Bg,Cf,Cg,Df,Dg} {
						\fill (\c) circle (2.5pt);
					}
					\node[below=2pt] at (A) {$1$};
					\node[above=2pt] at (B) {$1$};
					\node[below right] at (C) {$1$};
					\node[below left] at (D) {$1$};
				\end{tikzpicture}
				\caption{Three-cycle graph with a specified slope structure on $\Rat(K)$, defining a combinatorial limit linear series $M \subseteq \Rat(K)$.}
				\label{fig:three_cycles}
			\end{figure}
			
			We can define a suitable closed sub-semimodule $M \subseteq \Rat(K)$ of rank two of functions compatible with this choice of slopes. The tropical Weierstrass locus $\Wloc(M)$ of the semimodule $M$, in the sense of Section~\ref{subsec:augmented_metric_graphs_semimodules} (with $\g = 0$), contains the bridge edges and the points of coordinates $\frac 1 3$ and $\frac 2 3$ on the circles (see Figure~\ref{fig:three_cycles_weierstrass}, left). In particular, $M$ is not \wfinite{}. The clls Weierstrass divisor $W^\alg(M)$ is also shown in the figure (right). In particular, $M$ is not $\g$-effective. \qedhere
			
			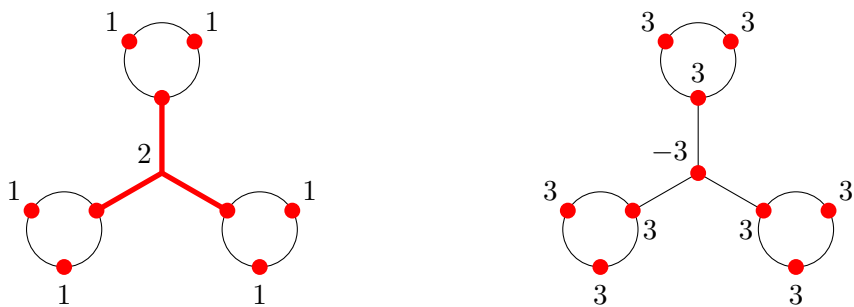
\begin{figure}[h!]
				\begin{minipage}{0.45\textwidth}
				\centering
				\begin{tikzpicture}[scale=.65]
					\coordinate (A) at (0,0);
					\coordinate (B) at (0,1);
					\coordinate (C) at (-0.87,-0.5);
					\coordinate (D) at (0.87,-0.5);
					\coordinate (Bc) at (0,1.5);
					\coordinate (Cc) at (-1.3,-0.75);
					\coordinate (Dc) at (1.3,-0.75);
					\coordinate (Bf) at ($(Bc) + (30:0.5)$);
					\coordinate (Bg) at ($(Bc) + (150:0.5)$);
					\coordinate (Cf) at ($(Cc) + (150:0.5)$);
					\coordinate (Cg) at ($(Cc) + (-90:0.5)$);
					\coordinate (Df) at ($(Dc) + (270:0.5)$);
					\coordinate (Dg) at ($(Dc) + (30:0.5)$);
					
					\draw[color=red,line width=1.2pt] (A) -- (B);
					\draw[color=red,line width=1.2pt] (A) -- (C);
					\draw[color=red,line width=1.2pt] (A) -- (D);
					\draw[radius=0.5] (Bc) circle;
					\draw[radius=0.5] (Cc) circle;
					\draw[radius=0.5] (Dc) circle;
					
					\foreach \c in {B,C,D,Bf,Bg,Cf,Cg,Df,Dg} {
						\fill[color=red] (\c) circle (3pt);
					}
					\node[above right] at (Bf) {$1$};
					\node[above left] at (Bg) {$1$};
					\node[above left] at (Cf) {$1$};
					\node[below=3pt] at (Cg) {$1$};
					\node[below=3pt] at (Df) {$1$};
					\node[above right] at (Dg) {$1$};
					\node[above left] at (A) {$2$};
				\end{tikzpicture}
				\end{minipage}
				\begin{minipage}{0.45\textwidth}
				\centering
				\begin{tikzpicture}[scale=.65]
					\coordinate (A) at (0,0);
					\coordinate (B) at (0,1);
					\coordinate (C) at (-0.87,-0.5);
					\coordinate (D) at (0.87,-0.5);
					\coordinate (Bc) at (0,1.5);
					\coordinate (Cc) at (-1.3,-0.75);
					\coordinate (Dc) at (1.3,-0.75);
					\coordinate (Bf) at ($(Bc) + (30:0.5)$);
					\coordinate (Bg) at ($(Bc) + (150:0.5)$);
					\coordinate (Cf) at ($(Cc) + (150:0.5)$);
					\coordinate (Cg) at ($(Cc) + (-90:0.5)$);
					\coordinate (Df) at ($(Dc) + (270:0.5)$);
					\coordinate (Dg) at ($(Dc) + (30:0.5)$);
					
					\draw (A) -- (B);
					\draw (A) -- (C);
					\draw (A) -- (D);
					\draw[radius=0.5] (Bc) circle;
					\draw[radius=0.5] (Cc) circle;
					\draw[radius=0.5] (Dc) circle;
					
					\foreach \c in {A,B,C,D,Bf,Bg,Cf,Cg,Df,Dg} {
						\fill[color=red] (\c) circle (3pt);
					}
					\node[above=2pt] at (B) {$3$};
					\node[above right] at (Bf) {$3$};
					\node[above left] at (Bg) {$3$};
					\node[below right] at (C) {$3$};
					\node[above left] at (Cf) {$3$};
					\node[below=3pt] at (Cg) {$3$};
					\node[below left] at (D) {$3$};
					\node[below=3pt] at (Df) {$3$};
					\node[above right] at (Dg) {$3$};
					\node[above left] at (A) {$-3$};
				\end{tikzpicture}
				\end{minipage}
				\caption{The tropical Weierstrass locus $\Wloc(M)$ (left) and the clls Weierstrass divisor $W^\alg(M)$ (right).}
				\label{fig:three_cycles_weierstrass}
			\end{figure}
		\end{example}
	
	\subsection{Comparison with the tropical Weierstrass locus}
		
		The following proposition shows that the notion of clls Weierstrass divisor can be viewed as a refinement of the tropical Weierstrass locus defined in Section~\ref{sec:weierstrass_locus_generalized_setup}.
		
		\begin{prop}[Comparison of the tropical and clls Weierstrass loci] \label{prop:comparison_clls_tropical_weierstrass_weights}
			Suppose $M \subseteq \Rat(D)$ is a combinatorial limit linear series of rank $r$ with clls Weierstrass divisor $W^{\alg}(M, \g)$. 
			Let $\Wloc(M, \g)$ denote its Weierstrass locus, defined as in Section~\ref{subsec:augmented_metric_graphs_semimodules}. 
			If $A \subseteq \Gamma$ is closed, connected, and $\Wloc(M, \g)$-measurable, then we have the equality
			\[
				\degrest{W^\alg(M, \g)}{A} = (r + 1) \, \hatweight(A; M, \g).
			\]
			In particular, if $M$ is \wfinite{} as a semimodule, then the following equality holds:
			\[
				W^\alg(M, \g) = (r + 1) \, W(M, \g).
			\]
		\end{prop}
		
		\begin{proof}
			We have
			\[
				\degrest{W^\alg(M, \g)}{A} = (r + 1) \sum_{x \in A} D(x) + \frac{r(r + 1)}{2} \sum_{x \in A} K(x) - \sum_{x \in A} \mleft(\sum_{\nu \in \T_x(\Gamma)} \, \sum_{j = 0}^r s_j^\nu \mright)
			\]	
			where $K$ denotes the canonical divisor on $(\Gamma, \g)$ (see Definition~\ref{def:canonical_divisor_augmented_graph}) and $s_j^\nu = s_j^\nu(M)$.
			The terms $(r + 1) \, D(x)$ add up to the term $(r + 1) \, \degrest{D}{A}$. Remark~\ref{rem:canonical_sum_augmented} yields that the terms $K(x)$ add up to $2 g(A) - 2 + 2 \sum_{x \in A} \g(x) + \outval(A)$, where $\outval(A) \coloneqq \abs{\partialout A}$ is the number of outgoing branches from $A$.
			
			The terms in the third part can be rearranged as a sum over directed edges of $A$, using some compatible model. Each edge has two in-going tangent directions, and the slope sums cancel out for this pair $(\nu, \overline \nu)$ of opposing in-going directions since $s_j^\nu + s_{r - j}^{\overline \nu} = 0$. The only terms that do not cancel are the tangent directions that point out of $A$, i.e.,
			\[
				\sum_{x \in A} \mleft( \sum_{\nu \in \T_x(\Gamma)} \, \sum_{j = 0}^r s_j^\nu \mright) = \sum_{\nu \in \partialout A} \mleft(\sum_{j = 0}^r s_j^\nu \mright).
			\]
			
			Combining these terms, we have
			\begin{align*}
				\degrest{W^\alg(M, \g)}{A} &= (r + 1) \degrest{D}{A} + \frac{r(r + 1)}{2} \mleft( 2 g(A, \g) - 2 + \outval(A)\mright) \\
				&\qquad\qquad - \sum_{\nu \in \partialout A} \sum_{j = 0}^r s_j^\nu \\
				&= (r + 1) \degrest{D}{A} + r(r + 1)(g(A, \g) - 1) - \sum_{\nu \in \partialout A} \sum_{j = 0}^r (s_j^\nu - j).
			\end{align*}
			Finally, we use the fact that $s_j^\nu = j + s_0^\nu$ for every $j$ and for tangent directions $\nu$ outside the Weierstrass locus $\Wloc(M, \g)$, by Theorem~\ref{thm:description_slopes_semimodules}. Thus,
			\[
				\degrest{W^\alg(M, \g)}{A} = (r + 1) \mleft(\degrest{D}{A} + (g(A, \g) - 1) \, r - \sum_{\nu \in \partialout A} s_0^\nu\mright),
			\]
			which, using the same technique as in the proof of Theorem~\ref{thm:weight_measure}, gives the first statement. The second statement follows from the first by the expression of the Weierstrass weight of a connected component of the tropical Weierstrass locus that is reduced to a point.
		\end{proof}
		
		We have the following extension of the above proposition, using the notion of tangential ramifications introduced later in Section~\ref{sec:tangential_ramification}.
		In particular, the statement holds even if $A$ is not $\Wloc(M, \g)$-measurable.
		
		\begin{prop}
			Keeping the same notation as in Proposition~\ref{prop:comparison_clls_tropical_weierstrass_weights}, for any closed, connected $A \subseteq \Gamma$, the following equality holds
			\begin{align*}
				\degrest{W^\alg(M, \g)}{A} &= (r + 1) \mleft( \degrest{D}{A} + (g(A, \g) - 1) \, r - \sum_{\nu \in \partialout A} s_0^\nu(M) \mright) \\
				&\qquad\qquad - \sum_{\nu \in \partialout A} \, \sum_{j = 0}^r \alpha_j^\nu(M),
			\end{align*}
			where $\alpha_j^\nu(M) \coloneqq s_j^\nu(M) - j - s_0^\nu(M)$ are the tangential ramifications along $\nu$.
		\end{prop}
	
	\subsection{Tropicalization of Weierstrass loci} \label{subsec:clls_tropicalization}
		
		The goal of this section is to prove Theorem~\ref{thm:weierstrass_tropicalization_extended}, using the machinery developed for semimodules on augmented metric graphs (see Section~\ref{subsec:augmented_metric_graphs_semimodules}). This provides a precise link between tropical Weierstrass loci and the tropicalization of Weierstrass divisors on algebraic curves. Using this result, we will deduce Theorem~\ref{thm:bound_weierstrass_tropicalization_intro}.
		
		Let $X$ be a smooth proper curve of genus $g$ over an algebraically closed non-Archimedean field $\K$ of arbitrary characteristic with a non-trivial valuation. Let $\cL = \mathcal{O}(\cD)$ be a line bundle of positive degree $d$ on $X$. Let $H$ be a vector subspace of global sections of $\cL$ of rank $r$ (i.e., $\dim H = r + 1$), that we naturally view in the function field of $X$. When $\K$ has positive characteristic, we additionally suppose that $\cL$ is classical~\cite{Lak81, Nee84}, that is, the gap sequence of $H$ is the standard sequence $\{0, 1, \dots, r\}$. We denote by $\cW = \cW(H)$ the corresponding Weierstrass divisor on $X$. Recall that $\cW$ is the zero divisor of a global section, called the Wronskian, of the line bundle $\omega!_X^{\otimes {r(r + 1)}/2} \otimes \cL^{\otimes (r + 1)}$, see~\cite{Lak81} and Section~\ref{app:wronskianB}. In particular, we have
		\[ \deg(\cW) = \frac{r(r + 1)}{2} \, (2g - 2) + (r + 1) \, d = (r + 1) \, (d - r + r g). \]
		
		Let $(\Gamma, \g)$ be a skeleton of $X^{\mathrm{an}}$, and let $\tau \colon X^\mathrm{an} \to \Gamma$ denote the specialization map. Let $W \coloneqq \tau_*(\cW)$ be the specialization of $\cW$ to $\Gamma$. Note that $(\Gamma, \g)$ is an augmented metric graph. We let $D \coloneqq \tau_*(\cD)$ be the specialization of $\cD$ to $\Gamma$, and let $M \subseteq \Rat(D)$ be the sub-semimodule consisting of the tropicalizations of non-zero rational functions in $H$. It follows from the slope formula that the divisorial rank of $M$ is equal to the rank of $H$, see~\cite[Thm.~8.3]{AG22} and~\cite[Prop.~4.1]{JP22}.
		
		The following theorem compares the algebraic Weierstrass divisor of $H$ on the curve $X$ with the tropical Weierstrass divisor of $M$ on the augmented metric graph $(\Gamma, \g)$.
		
		\begin{thm}[Algebraic versus tropical weights: general case] \label{thm:weierstrass_tropicalization_extended}
			Keeping the above notation, let $A$ be a closed, connected, $\Wloc \mleft(M, \g\mright)$-measurable subset of $\Gamma$. Then, the total weight of Weierstrass points of $\cW$ that tropicalize to $A$ is given by
			\[
				\degrest{\cW}{\tau^{-1}(A)} = (r + 1) \, \hatweight(A; M, \g)
			\]
			where
			\[
				\hatweight(A; M, \g) = \degrest{D}{A} + \mleft(g(A) + \sum_{x \in A} \g(x) - 1\mright) \, r - \sum_{\nu \in \partialout A} s_0^\nu(M).
			\]
			In particular, if $M$ is \wfinite, we have the following equality of divisors on $\Gamma$:
			\[
				\tau_*\mleft(\cW\mright) = (r + 1) \, W(M, \g).
			\]
		\end{thm}
		
		Before proceeding to the proof, a remark is in order.
		
		\begin{remark} \label{rem:weierstrass_tropicalization_extended}
			By Proposition~\ref{prop:coherence}, Theorem~\ref{thm:weierstrass_tropicalization_extended} holds in a slightly more general setting. Let $M'$ be any closed sub-semimodule of $\Rat(D)$ of divisorial rank $r$ containing $M$. Then, we have
			\[
				\degrest{\cW}{\tau^{-1}(A)} = (r + 1) \, \hatweight(A; M', \g)
			\] 
			for every $\Wloc \mleft(M', \g\mright)$-measurable subset $A$ of $\Gamma$.
		\end{remark}
				
		\begin{proof}[Proof of Theorem~\ref{thm:weierstrass_tropicalization_extended}]
			In the case the residue field of $\K$ has characteristic zero, we use Theorem~\ref{thm:redW_zero}, which provides a description of the divisor $W = \tau_*(\cW)$ in terms of slope structures. As explained in Section~\ref{app:slopes}, the slopes at any point $x$ and any unit tangent vector $\nu \in \T_x(\Gamma)$ of elements of the tropicalization $M$ of $H$ form a set of $r + 1$ integers $s_0^\nu, s_1^\nu, \dots, s_r^\nu$. The definition of the Weierstrass divisor associated to a tropical linear series is chosen to ensure the equality $W = W^\alg(M, \g)$, which implies
			\[
				\degrest{\cW}{\tau^{-1}(A)} = \degrest{W^\alg(M, \g)}{A}.
			\]
			Proposition~\ref{prop:comparison_clls_tropical_weierstrass_weights} states that if $A$ is $\Wloc(M, \g)$-measurable, then
			\[
				\degrest{W^\alg(M, \g)}{A} = (r + 1) \, \hatweight(A; M, \g),
			\]
			from which the result follows.
			
			In the general case, when the characteristic of $\K$ is arbitrary and the gap sequence of $H$ is standard, we use the description of the reduction of the Weierstrass divisor to the skeleton given in Theorem~\ref{thm:redW_general}. Using the notation of Section~\ref{sec:tropicalization_Weierstrass_points}, letting $W = \tau_*(\cW)$, we have
	 		\[ W(x) = (r + 1) D(x) + \frac{r(r + 1)}2 K(x) - \sum_{\nu \in \T_x(\Gamma)} \slope_{\nu} F, \]
			with $F = \Trop(\Wron_{\cF})$. Furthermore, $\slope_{\nu} F = \frac{r(r + 1)}2 + \ord_{\ssp_x^\nu} \widetilde{\Wron_{\cF}}_{_x}$.
			
			Since the slopes along the unit tangent vectors $\nu \in \T_x(\Gamma)$ that are outgoing from $A$ form a consecutive sequence of integers, by Proposition~\ref{prop:reduction_inteval} we infer that the quantity $\ord_{\ssp_x^\nu} \widetilde{\Wron_{\cF}}_{_x}$ is equal to $s^\nu_0 + \dots + s^\nu_r$. Using Theorem~\ref{thm:redW_general}, we get $\slope_{\nu} F = s^\nu_0 + \dots + s^\nu_r$.

			
			Moreover, since $F$ belongs to $\Rat(\Gamma)$, the total sum of the slopes of $F$ for the edges that appear in the interior of $A$ vanishes. We infer that
			\begin{align*}		
				\degrest{\cW}{\tau^{-1}(A)} &= (r + 1) \sum_{x \in A} D(x) + \frac{r(r + 1)}{2} \sum_{x \in A} K(x) - \sum_{x \in A} \sum_{\nu \in \T_x(\Gamma)} \slope_\nu F \\
				&= (r + 1) \sum_{x \in A} D(x) + \frac{r(r + 1)}{2} \sum_{x \in A} K(x) - \sum_{x \in A} \mleft(\sum_{\nu \in \T_x(\Gamma)}\sum_{j = 0}^r s_j^\nu \mright) \\
				&= \degrest{W^\alg(M, \g)}{A} = (r + 1) \, \hatweight(A; M, \g),
			\end{align*}
			as required. 
		\end{proof}

	\subsection{Proofs of Theorems~\ref{thm:weierstrass_tropicalization_intro} and~\ref{thm:bound_weierstrass_tropicalization_intro}} \label{subsec:proofs_tropicalization_theorems}
		
		We deduce Theorem~\ref{thm:weierstrass_tropicalization_intro} from Theorem~\ref{thm:weierstrass_tropicalization_extended}. 
				
		\begin{proof}[Proof of Theorem~\ref{thm:weierstrass_tropicalization_intro}]
			Since $D$ and $\cD$ have the same rank $r$, we can plug $H \coloneqq \Rat(\cD)$ and $M' \coloneqq \Rat(D)$ into Remark~\ref{rem:weierstrass_tropicalization_extended}, following Theorem~\ref{thm:weierstrass_tropicalization_extended}, to get
			\begin{equation} \label{eq:tropicalization_divisor}
				\degrest{\cW(\cD)}{\tau^{-1}(A)} = (r + 1) \mleft(\degrest{D}{A} + r\mleft(g(A, \g) - 1\mright) - \sum_{\nu \in \partialout A} s_0^\nu(D) \mright).
			\end{equation}
			In the context of Theorem~\ref{thm:weierstrass_tropicalization_intro}, $\g = 0$. The result follows.
		\end{proof}
		
		Using this result, we can prove Theorem~\ref{thm:bound_weierstrass_tropicalization_intro}.
		
		\begin{proof}[Proof of Theorem~\ref{thm:bound_weierstrass_tropicalization_intro}]
			This follows from the combination of Theorem~\ref{thm:number_w_points_M}, Proposition~\ref{prop:bound_red_degree_semimodule} and Theorem~\ref{thm:weierstrass_tropicalization_extended}.
		\end{proof}
	
			
\section{Further examples and discussions} \label{sec:further_examples_discussions}
	
	We here discuss several examples and questions in order to illustrate the results of the previous sections.
	
	\subsection{Dipole graph} \label{subsec:dipole_graph}
		
\begin{figure}[h!]
			\begin{minipage}{0.45\textwidth}
			\centering
			\begin{tikzpicture}[scale=.65]
				\coordinate (A) at (0,0);
				\coordinate (B) at (4,0);
				
				\draw (A) .. controls (0,1.5) and (4,1.5) .. (B);
				\draw (A) .. controls (0,-1.5) and (4,-1.5) .. (B);
				\draw (A) .. controls (1,0.5) and (3,0.5) .. (B);
				\draw (A) .. controls (1,-0.5) and (3,-0.5) .. (B);
				
				\foreach \c in {A,B} {
					\fill (\c) circle (2.5pt);
				}
				\node[left] at (A) {$2$};
				\node[right] at (B) {$2$};
			\end{tikzpicture}
			\end{minipage}
			\begin{minipage}{0.45\textwidth}
			\centering
			\begin{tikzpicture}[scale=.65]
				\coordinate (A) at (0,0);
				\coordinate (B) at (4,0);
				\coordinate (C) at (1.1,1);
				\coordinate (D) at (2.9,1);
				\coordinate (E) at (1.3,0.33);
				\coordinate (F) at (2.7,0.33);
				
				\draw (A) .. controls (0,1.5) and (4,1.5) .. (B);
				\draw (A) .. controls (0,-1.5) and (4,-1.5) .. (B);
				\draw (A) .. controls (1,0.5) and (3,0.5) .. (B);
				\draw (A) .. controls (1,-0.5) and (3,-0.5) .. (B);
				\draw[color=red,line width=1.4pt] (C) ..
						controls (1.6,1.15) and (2.4,1.15) .. (D);
				\draw[color=red,line width=1.4pt] (E) ..
						controls (2,0.4) .. (F);
				\foreach \c in {C,D,E,F} {
					\fill[color=red] (\c) circle (3pt);
				}

				\begin{scope}[rotate around={180:(2,0)}]
					\coordinate (Cr) at (1.1,1);
					\coordinate (Dr) at (2.9,1);
					\coordinate (Er) at (1.3,0.33);
					\coordinate (Fr) at (2.7,0.33);
				
					\draw[color=red,line width=1.4pt] (Cr) ..
							controls (1.6,1.15) and (2.4,1.15) .. (Dr);
					\draw[color=red,line width=1.4pt] (Er) ..
							controls (2,0.4) .. (Fr);
					\foreach \c in {Cr,Dr,Er,Fr} {
						\fill[color=red] (\c) circle (3pt);
					}
				\end{scope}
								
				\node[above] at (2,1.1) {$2$};
				\node[above] at (2,0.39) {$2$};
				\node[below] at (2,-1.1) {$2$};
				\node[below] at (2,-0.39) {$2$};
			\end{tikzpicture}
			\end{minipage}
			\caption{Dipole graph of genus $g = 3$ and its Weierstrass locus $\Wloc(K)$.}
			\label{fig:dipole}
		\end{figure}
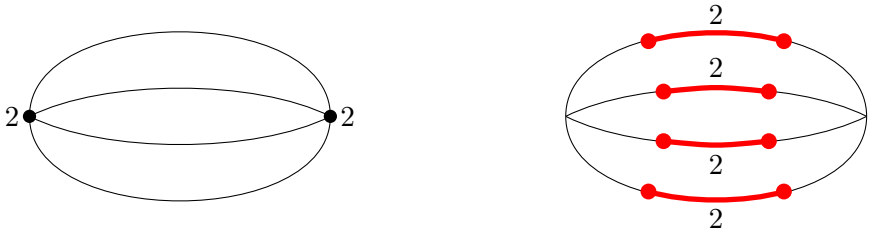

		Suppose $\Gamma$ is a dipole graph of genus $g \geq 2$ (also known as a ``banana'' graph), consisting of two vertices joined by $g + 1$ edges, possibly of different lengths. The canonical divisor $K$ has coefficient $g - 1$ on each vertex. The Weierstrass locus $\Wloc(K)$ consists of the interval $[\ell/g, \, \ell - \ell/g]$ on every edge, with $\ell$ the length of that edge (see Figure~\ref{fig:dipole} for $g = 3$). Each component $C \subseteq \Wloc(K)$ has two outgoing directions, and in each outgoing direction, the minimum slope is $s_0^\nu = - (g - 1)$. By Theorem~\ref{thm:canonical_weight}, the Weierstrass weight of each component is $g - 1$, and the total Weierstrass weight is $g^2 - 1$, as expected. 
		
	\subsection{Cube graph}
		
\begin{figure}[h!]
			\begin{minipage}{0.45\textwidth}
			\centering
			\begin{tikzpicture}[scale=.65]
				\coordinate (A) at (0,0);
				\coordinate (B) at (3.2,0);
				\coordinate (C) at (3.2,3.2);
				\coordinate (D) at (0,3.2);
				\coordinate (E) at (0.8,0.8);
				\coordinate (F) at (2.4,0.8);
				\coordinate (G) at (2.4,2.4);
				\coordinate (H) at (0.8,2.4);
				
				\draw (A) -- (B) -- (C) -- (D) -- cycle;
				\draw (E) -- (F) -- (G) -- (H) -- cycle;
				\foreach \c/\d in {A/E, B/F, C/G, D/H} {
					\draw (\c) -- (\d);
				}
				
				\foreach \c in {A,B,C,D,E,F,G,H} {
					\fill (\c) circle (3pt);
				}
			\end{tikzpicture}
			\end{minipage}
			\begin{minipage}{0.45\textwidth}
			\centering
			\begin{tikzpicture}[scale=.65]
				\coordinate (A) at (0,0);
				\coordinate (B) at (3.2,0);
				\coordinate (C) at (3.2,3.2);
				\coordinate (D) at (0,3.2);
				\coordinate (E) at (0.8,0.8);
				\coordinate (F) at (2.4,0.8);
				\coordinate (G) at (2.4,2.4);
				\coordinate (H) at (0.8,2.4);
				
				\draw (A) -- (B) -- (C) -- (D) -- cycle;
				\draw (E) -- (F) -- (G) -- (H) -- cycle;
				\foreach \c/\d in {A/E, B/F, C/G, D/H} {
					\draw (\c) -- (\d);
				}
				
				\foreach \p/\q in {{(1.28,0)/(1.92,0)}, {(1.28,3.2)/(1.92,3.2)}, {(0,1.28)/(0,1.92)}, {(3.2,1.28)/(3.2,1.92)}, {(1.44,0.8)/(1.76,0.8)},{(1.44,2.4)/(1.76,2.4)}, {(0.8,1.44)/(0.8,1.76)}, {(2.4,1.44)/(2.4,1.76)}, {(0.32,0.32)/(0.48,0.48)}, {(2.72,2.72)/(2.88,2.88)}, {(2.88,0.32)/(2.72,0.48)}, {(0.32,2.88)/(0.48,2.72)}} {
					\draw[color=red,line width=1.2pt] \p -- \q;
					\fill[color=red] \p circle (2.1pt);
					\fill[color=red] \q circle (2.1pt);
				}
			\end{tikzpicture}
			\end{minipage}
			\caption{Cube graph with its Weierstrass locus $\Wloc(K)$.}
			\label{fig:cube}
		\end{figure}
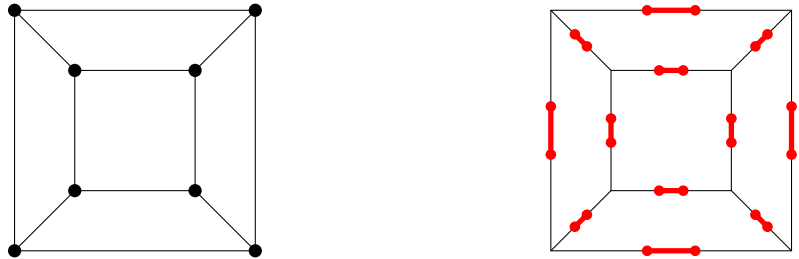

		The cube graph is shown in Figure~\ref{fig:cube}, with all edges of length one. It has genus $5$ and the canonical divisor $K$ has rank $4$. The Weierstrass locus $\Wloc(K)$ consists of the closed segment $[2/5, \, 3/5]$ on each edge, and excludes the vertices.
		Each component $C$ of $\Wloc(K)$ has out-valence $2$, with minimum slopes in $\Rat(K)$ in each outgoing direction equal to $-3$. Theorem~\ref{thm:canonical_weight} gives $\weight(C) = 2$.
		There are $12$ components, so the total weight is $24$.
		
	\subsection{Bridge edges} \label{subsec:bridge_edges}
		
		Earlier in Example~\ref{ex:barbell_graph}, we found that the Weierstrass locus $\Wloc(K)$ of the barbell graph contains the entire bridge edge.
		Here we show that this behavior holds in general.
		
		
		\begin{thm}[Weierstrass loci and bridge edges]
			Let $\Gamma$ be a metric graph that has a bridge edge $e$ such that each component of $\Gamma \setminus \mathring e$ has positive genus. Then, the edge $e$ is contained in the canonical Weierstrass locus $\Wloc(K)$.
		\end{thm}
		
		\begin{proof}
			To show this, let $u_1$ and $u_2$ denote the endpoints of $e$, and $\Gamma_1$ and $\Gamma_2$ be the components of $\Gamma \setminus \mathring e$ containing $u_1$ and $u_2$, respectively. If $g$, $g_1$ and $g_2$ are the genera of $\Gamma$, $\Gamma_1$ and $\Gamma_2$ respectively, then $g = g_1 + g_2$. Let $r = g - 1$ be the rank of the canonical divisor on $\Gamma$. We want to show that we can move $r + 1 = g_1 + g_2$ coefficients to every point $x \in e$. For $i = 1, 2$, denoting by $K_i$ the canonical divisor of $\Gamma_i$, we have $r_{\Gamma_i}\mleft(K \rest{\Gamma_i} - (u_i)\mright) = r_{\Gamma_i}(K_i) = g_i - 1 \geq 0$, which implies that, using only functions on $\Gamma$ that are constant outside $\Gamma_i$, we can move $g_i$ coefficients to $u_i$. It is easy to see that we can move coefficients along $e$ to put $g_1 + g_2$ coefficients at $x$.
		\end{proof}
		
	\subsection{Cases where the whole graph is Weierstrass} \label{subsec:graph_all_weierstrass}
		
		
		Example~\ref{ex:complete_graph_4} treated the complete graph on four vertices with unit edge lengths. Now consider the case $\Gamma$ is the complete graph on $n \geq 5$ vertices with unit edge lengths. This graph has genus $g = \frac{n^2 - 3 n + 2}{2}$, and the canonical divisor $K$ has rank $\frac{n^2 - 3 n}{2}$. The Weierstrass locus of $K$ turns out to be the whole $\Gamma$, as can be verified in a straightforward way.
		
\begin{figure}[h!]
			\begin{minipage}{0.45\textwidth}
			\centering
			\begin{tikzpicture}[scale=.65]
				\coordinate (A) at (0,1);
				\coordinate (D1) at (-2,0);
				\coordinate (D2) at (-0.5,0);
				\coordinate (D3) at (2,0);
				\coordinate (E4) at (0.5,0);
				\coordinate (E5) at (0.7,0);
				\coordinate (E6) at (0.9,0);
				\coordinate (M) at (-30:1.3);
				\coordinate (N) at (-45:1.3);
				\coordinate (P) at (-60:1.3);
				
				\foreach \c in {D1,D2,D3} {
					\draw (A) -- (\c);
					\draw[radius=0.5] (\c) + (0,-0.5) circle;
				}
				
				\foreach \c in {E4,E5,E6} {
					\draw (A) -- (\c);
					\draw[radius=0.01] (\c) + (0,-0.2) circle;
				}

				\fill (A) circle (3pt);
				\node[above right] at (A) {$d$};
			\end{tikzpicture}
			\end{minipage}
			\begin{minipage}{0.45\textwidth}
			\centering
			\begin{tikzpicture}[scale=.65]
				\coordinate (A) at (0,1);
				\coordinate (D1) at (-2,0);
				\coordinate (D2) at (-0.5,0);
				\coordinate (D3) at (2,0);
				\coordinate (E4) at (0.5,0);
				\coordinate (E5) at (0.7,0);
				\coordinate (E6) at (0.9,0);
				\coordinate (M) at (-30:1.3);
				\coordinate (N) at (-45:1.3);
				\coordinate (P) at (-60:1.3);
				
				\foreach \c in {D1,D2,D3} {
					\draw[color=red,line width=1.2pt] (A) -- (\c);
					\draw[radius=0.5,color=red,line width=1.2pt] (\c) + (0,-0.5) circle;
				}
				
				\foreach \c in {E4,E5,E6} {
					\draw[color=red,line width=1.2pt] (A) -- (\c);
					\fill[radius=1pt,color=red] (\c) + (0,-0.2) circle;
				}
				\fill[color=red] (A) circle (1pt);
			\end{tikzpicture}
			\end{minipage}
			\caption{The generalized barbell graph, the divisor $D$ and its Weierstrass locus $\Wloc(D)$.}
			\label{fig:multipod}
		\end{figure}
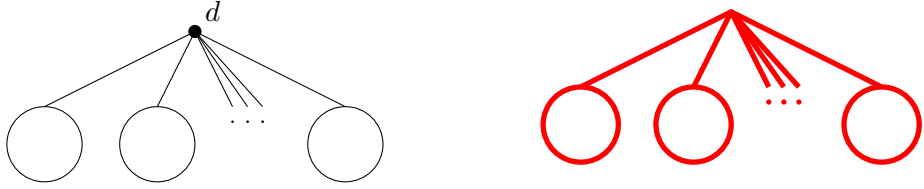

		We give a second such family for which the choice of the length function is free.
		See also~\cite[Ex.~4.6]{Richman18}. Let $\Gamma$ be the metric graph generalizing the barbell graph (Example~\ref{ex:barbell_graph}) to any number of cycles. More precisely, take $g \geq 2$ cycles of arbitrary length and join them all to a central vertex $v$ with a bridge edge of positive length, as in Figure~\ref{fig:multipod}. Consider the divisor $D = d \, (v)$, with $d \geq 3$. The rank of $D$ can be computed by straightforward verification or using the more general result~\cite[Prop.~4.8]{AB15}. For $d \leq 2 g - 2$, one finds $r = \lfloor \frac d 2\rfloor$, and for higher values of $d$, we have $r = d - g$. Since a divisor of positive degree on a cycle has rank one less than the degree, and since divisor coefficients can move freely on bridge edges, it is easy to show $D_x(x) \geq d - 1 \geq r + 1$ for every $x \in \Gamma$. Therefore, the Weierstrass locus is the whole graph.
		
		Note that in this family the quantity $\min_{x \in \Gamma} \, (D_x(x) - r)$ can be arbitrarily large. The existence of these two families, with different combinatorial properties, suggests the following.
		
		\begin{question}
			Provide a classification of all graphs $G$ that admit a length function and a divisor with Weierstrass locus the whole metric graph. Among them, what are the ones for which this property holds for every choice of edge lengths?
		\end{question}
		
		The treatment of Weierstrass points for curves over positive characteristic fields suggests the following possible modification of the theory of tropical Weierstrass points in the isolated cases where the whole graph is Weierstrass. We replace the rank $r$ with the integer
		\[ b = b(\Gamma, D) \coloneqq \min_{x \in \Gamma} D_x(x), \]
		and define the Weierstrass locus as the subset of points $x \in \Gamma$ verifying $D_x(x) \geq b + 1$. The weight of a connected component $C$ of this modified Weierstrass locus is modified by setting
		\[ \weight(C; D) \coloneqq \degrest{D}{C} + \mleft(g(C) - 1\mright) \, b - \sum_{\nu \in \partialout C} s_0^\nu(D). \]
		
		This leads to a consistent theory on the tropical side, with the weights of components of the Weierstrass locus adding up to $d - b + b g$ (instead of $d - r + r g$). This is reminiscent of the setting of curves in the situation where the standard sequence of vanishing orders differs from the sequence $0, 1, \dots, r$, cf.~\cite{Lak81}. However, at this point, we are not aware of any geometric meaning to this tropical count.
	
	\subsection{Augmented cycle with one point of positive genus} \label{subsec:augmented_cycle}
		
		We compute Weierstrass loci for the canonical divisor with respect to the canonical and generic linear systems on an augmented cycle on which one point has positive genus, generalizing Examples~\ref{ex:augmented_cycle_one_point_can_2} and~\ref{ex:augmented_cycle_one_point_gen_2}. The canonical case recovers a result of Diaz~\cite[Thm.~A2.1]{Dia85}: the generic non-separating node on a uninodal stable curve is a limit of exactly $g(g - 1)$ Weierstrass points on nearby smooth curves.
		
	 	Let $a$ be a positive integer, and consider the augmented metric graph $(\Gamma, \g)$ where $\Gamma$ is the cycle of length one, parametrized by the interval $[0, \, 1]$, the single vertex $v$ coincides with the endpoints $v = 0 = 1$, and $\g(v) = a$. The genus of this augmented metric graph is $g = a + 1$.
		
		\subsubsection{The case of the canonical linear system} \label{subsubsec:augmented_cycle_one_can}
			
			We expand on Example~\ref{ex:augmented_cycle_one_point_can_2} for which $a = 2$ was fixed. Consider the canonical divisor $K$ and the associated canonical semimodule $\KRat(\g)$, as defined in Section~\ref{subsec:weierstrass_points_augmented_graphs_can}. The rank is $r = g - 1 = a$ according to Theorem~\ref{thm:rank_semimodule_can}, and the total weight of the Weierstrass locus is $g^2 - 1 = a^2 + 2 a$. The Weierstrass locus consists of the vertex $v$ and all the points of the form $\frac k{a + 1}$ for $k = 1, \dots, a$. The Weierstrass weights are $\weight(v; K, \g) = a^2 + a$ and $\weight \mleft(\frac k{a + 1}; K, \g\mright) = 1$. Figure~\ref{fig:aug_cycle_can} shows the canonical divisor and its (canonical) Weierstrass locus depicted in the middle.
		
		\subsubsection{The case of the generic linear system} \label{subsubsec:augmented_cycle_one_gen}
			
			In the second case, we generalize Example~\ref{ex:augmented_cycle_one_point_gen_2} and consider the same divisor $K$ as above, but take the generic semimodule $\Ratgen(K, \g)$ as defined in Section~\ref{subsec:weierstrass_points_augmented_graphs_gen}. In this case, the rank is $r = g - 2 = a - 1$ (see Proposition~\ref{prop:rank_canonical_divisor_generic_setting}) and the total weight of the Weierstrass locus is $a^2 + a$. The Weierstrass points are $v$, and all the points $\frac k a$ with $k = 1, \dots, a - 1$. The weights are $\weight(v; K, \g) = a^2 + 1$ and $\weight \mleft(\frac k a; K, \g\mright) = 1$. Figure~\ref{fig:aug_cycle_can} shows the canonical divisor on the left, and its generic Weierstrass locus depicted on the right.
			
			We note that the Weierstrass loci are different even though they are both finite. The total weights are also different, as the underlying semimodules have different ranks.
			
			\begin{figure}[h!]
				\begin{minipage}{0.3\textwidth}
				\centering
				\begin{tikzpicture}[scale=.65]
					\coordinate (A) at (-1,0);
					\coordinate (B) at (0,0);
					
					\draw[radius=1.0] (B) circle;
					\fill (A) circle (2.5pt);
					\fill[opacity=0.2] (A) circle (5pt);
					
					\node[above left] at (A) {$\g(v) = a$};
					\node[below left] at (A) {$K(v) = 2 a$};
				\end{tikzpicture}
				\end{minipage}
				\begin{minipage}{0.3\textwidth}
				\centering
				\begin{tikzpicture}[scale=.65]
					\coordinate (A) at (-1,0);
					\coordinate (B) at (0,0);
					\coordinate (C) at (252:1);
					\coordinate (D) at (324:1);
					\coordinate (E) at (36:1);
					\coordinate (F) at (108:1);
					
					\draw[radius=1.0] (B) circle;
					\fill[opacity=0.2] (A) circle (5pt);
					
					\foreach \c in {A,C,D,E,F} {
						\fill[color=red] (\c) circle (3pt);
					}
	
					\node[left=3pt] at (A) {$a^2 + a$};
					\node[below left] at (C) {$1$};
					\node[right=2pt] at (D) {$1$};
					\node[right=2pt] at (E) {$1$};
					\node[above left] at (F) {$1$};
				\end{tikzpicture}
				\end{minipage}
				\begin{minipage}{0.3\textwidth}
					\centering
					\begin{tikzpicture}[scale=.65]
						\coordinate (A) at (-1,0);
						\coordinate (B) at (0,0);
						\coordinate (C) at (270:1);
						\coordinate (D) at (0:1);
						\coordinate (E) at (90:1);
						
						\draw[radius=1.0] (B) circle;
						\fill[opacity=0.2] (A) circle (5pt);
						
						\foreach \c in {A,C,D,E} {
							\fill[color=red] (\c) circle (3pt);
						}
		
						\node[left=3pt] at (A) {$a^2 + 1$};
						\node[below=3pt] at (C) {$1$};
						\node[right=2pt] at (D) {$1$};
						\node[above=2pt] at (E) {$1$};
					\end{tikzpicture}
					\end{minipage}
				\caption{An augmented cycle graph, with its canonical Weierstrass locus $\Wloc(K, \g)$ in the middle, and its Weierstrass locus $\Wlocgen(K, \g)$ on the right. The drawing is made for $a = 4$.}
				\label{fig:aug_cycle_can}
			\end{figure}
	
	\subsection{Augmented dipole graph} \label{sec:dipole_augmented}
		
		We now consider an augmented dipole graph made up of two vertices $u$ and $v$ joined by $n = h + 1$ edges of arbitrary lengths, where $h$ is the genus of the corresponding metric graph.
		We assume $\g$ has support in $\{u, v\}$, and denote by $a$ and $b$ the genus of $u$ and $v$, respectively, with $a \leq b$.
		This metric graph is the one that appears in the work by Esteves and Medeiros~\cite{EM02}.
		As we explained previously, the canonical linear series reflects the genericity of the points of intersection on each of the two components.
		
		The canonical divisor has coefficients $K(u) = h - 1 + 2 a$ and $K(v) = h - 1 + 2 b$.
		%
		%
		%
		The determination of the Weierstrass locus $\Wloc(K, \g)$ turns out to be quite involved in general, and its shape depends on the values of $a$, $b$, $h$ and the edge lengths.
		We illustrate the computation in one concrete example.
		
		Suppose $a = 3$, $b = 5$, $h = 2$, and all the edges have unit length. 
		The augmented graph has total genus $g = 10$, and its canonical divisor has rank $r = 9$.
		The Weierstrass locus is made up of the vertex $v$ (weight $50$), the union of the three segments $[0, \, 1/10]$ lying on each edge (weight $34$), the point of coordinate $6/10$ on each edge (weight $1$), and the segments $[3/10, \, 4/10]$ and $[8/10, \, 9/10]$ on each edge (each of weight $2$). See Figure~\ref{fig:dipole_augmented_second_case}. The total weight is $50 + 34 + 3 \cdot (2 + 1 + 2) = 99 = g^2 - 1$.
		
		\begin{figure}[h!]
			\begin{minipage}{0.45\textwidth}
			\centering
			\begin{tikzpicture}[scale=0.75]
				\coordinate (A) at (0,0);
				\coordinate (B) at (4,0);
				
				\draw (A) .. controls (0,1.5) and (4,1.5) .. (B);
				\draw (A) .. controls (0,-1.5) and (4,-1.5) .. (B);
				\draw (A) -- (B);
				
				\foreach \c in {A,B} {
					\fill (\c) circle (3pt);
					\fill[opacity=0.2] (\c) circle (5pt);
				}
				\node[above left] at (A) {$\g(u) = 3$};
				\node[below left] at (A) {$K(u) = 7$};
				\node[above right] at (B) {$\g(v) = 5$};
				\node[below right] at (B) {$K(v) = 11$};
			\end{tikzpicture}
			\end{minipage}
			\begin{minipage}{0.45\textwidth}
			\centering
			\begin{tikzpicture}[scale=0.75]
				\coordinate (A) at (0,0);
				\coordinate (B) at (4,0);
				
				\coordinate (Cc) at (4/10,0);
				\coordinate (Dc) at (12/10,0);
				\coordinate (Ec) at (16/10,0);
				\coordinate (Fc) at (24/10,0);
				\coordinate (Gc) at (32/10,0);
				\coordinate (Hc) at (36/10,0);
				
				\coordinate (Ct) at (0.16,0.45);
				\coordinate (Dt) at (1,0.98);
				\coordinate (Et) at (1.45,1.09);
				\coordinate (Ft) at (2.55,1.09);
				\coordinate (Gt) at (3.35,0.85);
				\coordinate (Ht) at (3.78,0.55);
				
				\coordinate (Cb) at (0.16,-0.45);
				\coordinate (Db) at (1.0,-0.98);
				\coordinate (Eb) at (1.45,-1.09);
				\coordinate (Fb) at (2.55,-1.09);
				\coordinate (Gb) at (3.35,-0.85);
				\coordinate (Hb) at (3.78,-0.55);
				
				\draw (A) .. controls (0,1.5) and (4,1.5) .. (B);
				\draw (A) .. controls (0,-1.5) and (4,-1.5) .. (B);
				\draw (A) -- (B);
				
				\foreach \c in {A, B} {
					\fill[opacity=0.2] (\c) circle (5pt);
				}
				
				\foreach \c in {B,Cc,Dc,Ec,Fc,Gc,Hc,Ct,Dt,Et,Ft,Gt,Ht,Cb,Db,Eb,Fb,Gb,Hb} {
					\fill[color=red] (\c) circle (3pt);
				}
				
				\draw[color=red,line width=1.2pt] (A) -- (Cc);
				\draw[color=red,line width=1.2pt] (Dc) -- (Ec);
				\draw[color=red,line width=1.2pt] (Gc) -- (Hc);
				
				\draw[color=red,line width=1.2pt] (A) .. controls (0,0.3) and (0.1,0.35) .. (Ct);
				\draw[color=red,line width=1.2pt] (Dt) .. controls (1.15,1.06) and (1.3,1.06) .. (Et);
				\draw[color=red,line width=1.2pt] (Gt) .. controls (3.5,0.78) and (3.6,0.74) .. (Ht);
				
				\draw[color=red,line width=1.2pt] (A) .. controls (0,-0.3) and (0.1,-0.35) .. (Cb);
				\draw[color=red,line width=1.2pt] (Db) .. controls (1.15,-1.06) and (1.3,-1.06) .. (Eb);
				\draw[color=red,line width=1.2pt] (Gb) .. controls (3.5,-0.78) and (3.6,-0.74) .. (Hb);
			
				\node[left=3pt] at (A) {$34$};
				\node[below=3pt] at (1.1,-0.95) {$2$};
				\node[below=3pt] at (Fb) {$1$};
				\node[below right=1pt] at (3.5,-0.65) {$2$};
				\node[right=3pt] at (B) {$50$};
			\end{tikzpicture}
			\end{minipage}
			\caption{Augmented dipole graph with $h = 2$, $a = 3$ and $b = 5$, edges all of unit length, and its Weierstrass locus $\Wloc(K, \g)$.}
			\label{fig:dipole_augmented_second_case}
		\end{figure}
	
	\subsection{Weierstrass divisor of a combinatorial limit linear series}
		
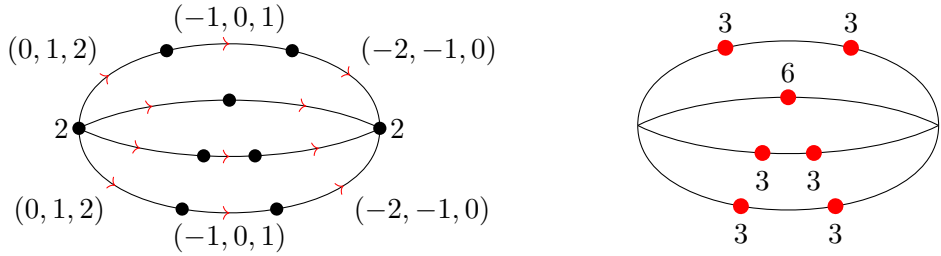
\begin{figure}[h!]
			\begin{minipage}{0.45\textwidth}
			\centering
			\begin{tikzpicture}[scale=.75]
				\coordinate (A) at (0,0);
				\coordinate (B) at (4,0);
				
				\draw[postaction=decorate,decoration={markings,mark=at position 1/6 with {\arrow[red]{>}; \node[above left]{$\{0, 1, 2\}$};},mark=at position 1/3 with {\fill circle[radius=2.5pt];},mark=at position 1/2 with {\arrow[red]{>}; \node[above]{$\{-1, 0, 1\}$};},mark=at position 2/3 with {\fill circle[radius=2.5pt];},mark=at position 5/6 with {\arrow[red]{>}; \node[above right]{$\{-2, -1, 0\}$};}}] (A) .. controls (0,1.5) and (4,1.5) .. (B);
				\draw[postaction=decorate,decoration={markings,mark=at position 3/16 with {\arrow[red]{>}; \node[below left]{$\{0, 1, 2\}$};},mark=at position 3/8 with {\fill circle[radius=2.5pt];},mark=at position 1/2 with {\arrow[red]{>}; \node[below]{$\{-1, 0, 1\}$};},mark=at position 5/8 with {\fill circle[radius=2.5pt];},mark=at position 13/16 with {\arrow[red]{>}; \node[below right]{$\{-2, -1, 0\}$};}}] (A) .. controls (0,-1.5) and (4,-1.5) .. (B);
				\draw[postaction=decorate,decoration={markings,mark=at position 1/4 with {\arrow[red]{>};},mark=at position 1/2 with {\fill circle[radius=2.5pt];},mark=at position 3/4 with {\arrow[red]{>};}}] (A) .. controls (1,0.5) and (3,0.5) .. (B);
				\draw[postaction=decorate,decoration={markings,mark=at position 5/24 with {\arrow[red]{>};},mark=at position 5/12 with {\fill circle[radius=2.5pt];},mark=at position 1/2 with {\arrow[red]{>};},mark=at position 7/12 with{\fill circle[radius=2.5pt];},mark=at position 19/24 with {\arrow[red]{>};}}] (A) .. controls (1,-0.5) and (3,-0.5) .. (B);
				
				\foreach \c in {A,B} {
					\fill (\c) circle (3pt);
				}
				\node[left] at (A) {$2$};
				\node[right] at (B) {$2$};
			\end{tikzpicture}
			\end{minipage}
			\begin{minipage}{0.45\textwidth}
			\centering
			\begin{tikzpicture}[scale=.75]
				\coordinate (A) at (0,0);
				\coordinate (B) at (4,0);
				
				\draw[postaction=decorate,decoration={markings,mark=at position 1/3 with {\fill[color=red] circle[radius=2.5pt]; \node[above=2pt]{$3$};},mark=at position 2/3 with {\fill[color=red] circle[radius=2.5pt]; \node[above=2pt]{$3$};}}] (A) .. controls (0,1.5) and (4,1.5) .. (B);
				\draw[postaction=decorate,decoration={markings,mark=at position 3/8 with {\fill[color=red] circle[radius=2.5pt]; \node[below=3pt]{$3$};},mark=at position 5/8 with {\fill[color=red] circle[radius=2.5pt]; \node[below=3pt]{$3$};}}] (A) .. controls (0,-1.5) and (4,-1.5) .. (B);
				\draw[postaction=decorate,decoration={markings,mark=at position 1/2 with {\fill[color=red] circle[radius=2.5pt]; \node[above=2pt]{$6$};}}] (A) .. controls (1,0.5) and (3,0.5) .. (B);
				\draw[postaction=decorate,decoration={markings,mark=at position 5/12 with {\fill[color=red] circle[radius=2.5pt]; \node[below=3pt]{$3$};},mark=at position 7/12 with {\fill[color=red] circle[radius=2.5pt]; \node[below=3pt]{$3$};}}] (A) .. controls (1,-0.5) and (3,-0.5) .. (B);
			\end{tikzpicture}
			\end{minipage}
			\caption{Dipole graph and its clls Weierstrass divisor $W^\alg(M)$.}
			\label{fig:dipole_clls_divisor}
		\end{figure}

		We go back to the non-augmented dipole graph with four edges (of unit length to simplify the notation), a particular case of the class of examples presented in Section~\ref{subsec:dipole_graph}. The genus is $g = 3$ and the rank of the canonical divisor $K$ is $r = 2$. Denote by $u$ and $v$ the two vertices and by $e_1, e_2, e_3$ and $e_4$ the four edges of $\Gamma$ (see Figure~\ref{fig:dipole_clls_divisor}, left). For $i = 1, 2, 3, 4$, let $t_i \in \mleft[0, \, \frac 1 6\mright]$. For each choice of the $t_i$'s, we will construct a semimodule $M$ of rank two that verifies condition $(\star)$ from Section~\ref{sec:weierstrass_locus_combinatorial_linear_series}, and compute its clls Weierstrass divisor $W^\alg(M)$ (here, $\g = 0$).
		
		Assume $t_i$'s are fixed. We prescribe the set of slopes taken by functions in $M$ as in Figure~\ref{fig:dipole_clls_divisor}. For each $i$, we endow the edge $e_i$ with the slope sets $\{0, 1, 2\}$ on the interval $\mleft[0, \, \frac 1 2 - t_i\mright]$, slopes $\{-1, 0, 1\}$ on the interval $\mleft[\frac 1 2 - t_i, \, \frac 1 2 + t_i\mright]$, and slopes $\{-2, -1, 0\}$ on the interval $\mleft[\frac 1 2 + t_i, \, 1\mright]$.
		
		We define $M$ as the sub-semimodule of $\Rat(K)$ consisting of all the functions that take one of the prescribed slopes above at any point of $\Gamma$ along any unite tangent vector. We thus get a 4-parameter family of pairs $(M,D)$ of rank two, $M \subseteq \Rat(K)$, that verify $(\star)$.
		
		The clls Weierstrass divisor is given by~\eqref{eq:clls_weierstrass_divisor}, and yields $W^\alg(M) = 3 \sum_{i = 1}^4 \mleft(\mleft(x_i\mright) + \mleft(y_i\mright)\mright)$, where $x_i$ and $y_i$ denote the points of coordinates $1/2 - t_i$ and $1/2 + t_i$ on the edge $e_i$, respectively. Figure~\ref{fig:dipole_clls_divisor} gives a visual rendering of $W^\alg(M)$ for the choice $(t_1, t_2, t_3, t_4) = \mleft(\frac 1 6, 0, \frac 1 {12}, \frac 1 8\mright)$. 
		
	\subsection{Weierstrass points of random combinatorial graphs} \label{subsec:Weierstrass_points_random_combinatorial_graphs}
		
		There exist combinatorial graphs without any Weierstrass points among their vertices. The dipole graph is an example of such a graph, see Figure~\ref{fig:dipole}. So is the cube graph, see Figure~\ref{fig:cube}.
		Such graphs are interesting from the point of view of arithmetic geometry, see~\cite[\S~4]{baker2008specialization} and~\cite{Ogg78, LN64, Atk67, AP03}.
		
		This seems, however, to be a rare phenomenon, as a computer verification of examples indicates. Examples of random graphs were created and visualized using Python, Matplotlib~\cite{matplotlib}, and NetworkX~\cite{networkx}.
		
		\begin{question} 
			What is the proportion of combinatorial graphs that do not have any Weierstrass point among their vertices? That is, what is the probability that a combinatorial graph on $n$ vertices has no Weierstrass point?
		\end{question}
		
		\begin{figure}[h!]
			\includegraphics[scale=0.42]{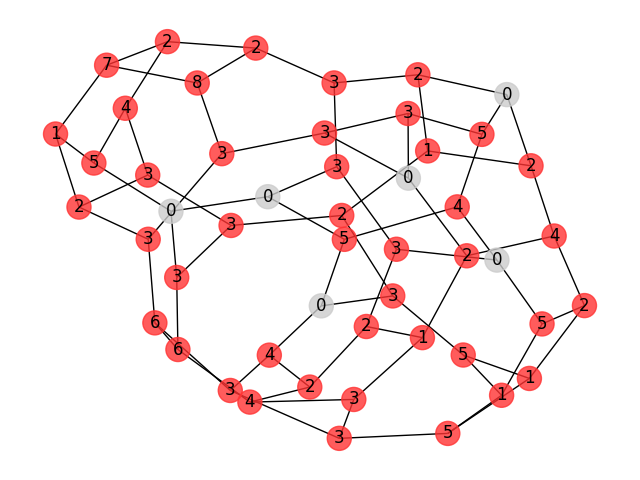}
			\caption{Random trivalent graph and its Weierstrass locus $\Wloc(K)$. The graph has genus $26$, and the vertex labels indicate the coefficients $K_v(v) - 25$.}
			\label{fig:random_trivalent}
		\end{figure}
		
		Randomness is understood within a class of graphs, for example regular graphs of given degree, or Erd\"os--R\'enyi random graphs. This is related to the following question of Baker.
		
		\begin{question}[Baker~\cite{baker2008specialization}]
			Provide a classification of combinatorial graphs without Weierstrass points among their vertices.
		\end{question}

\bibliography{Bibliography}

\end{document}

%% file: proof-figure2.pdf_t
\begin{picture}(0,0)%
\includegraphics{proof-figure2.pdf}%
\end{picture}%
\setlength{\unitlength}{3947sp}%
\begingroup\makeatletter\ifx\SetFigFont\undefined%
\gdef\SetFigFont#1#2#3#4#5{%
  \reset@font\fontsize{#1}{#2pt}%
  \fontfamily{#3}\fontseries{#4}\fontshape{#5}%
  \selectfont}%
\fi\endgroup%
\begin{picture}(12211,6739)(498,-8297)
\put(4726,-8161){\makebox(0,0)[lb]{\smash{{\SetFigFont{20}{24.0}{\familydefault}{\mddefault}{\updefault}{\color[rgb]{0,0,0}$C_3$}%
}}}}
\put(526,-4261){\makebox(0,0)[lb]{\smash{{\SetFigFont{20}{24.0}{\familydefault}{\mddefault}{\updefault}{\color[rgb]{0,0,0}$C_1$}%
}}}}
\put(4426,-1861){\makebox(0,0)[lb]{\smash{{\SetFigFont{20}{24.0}{\familydefault}{\mddefault}{\updefault}{\color[rgb]{0,0,0}$V'$}%
}}}}
\put(5251,-3511){\makebox(0,0)[lb]{\smash{{\SetFigFont{20}{24.0}{\familydefault}{\mddefault}{\updefault}{\color[rgb]{0,0,0}$C_2$}%
}}}}
\end{picture}%

%% file: tropical-weierstrass-points.bbl
\begin{bibdiv}
\begin{biblist}

\bib{AB15}{article}{
      author={Amini, Omid},
      author={Baker, Matthew},
       title={Linear series on metrized complexes of algebraic curves},
        date={2015},
        ISSN={0025-5831},
     journal={Math. Ann.},
      volume={362},
      number={1-2},
       pages={55\ndash 106},
         url={https://doi.org/10.1007/s00208-014-1093-8},
      review={\MR{3343870}},
}

\bib{amini2015lifting1}{article}{
      author={Amini, Omid},
      author={Baker, Matthew},
      author={Brugall\'{e}, Erwan},
      author={Rabinoff, Joseph},
       title={Lifting harmonic morphisms {I}: metrized complexes and
  {B}erkovich skeleta},
        date={2015},
        ISSN={2522-0144},
     journal={Res. Math. Sci.},
      volume={2},
       pages={Art. 7, 67},
         url={https://doi.org/10.1186/s40687-014-0019-0},
      review={\MR{3375652}},
}

\bib{AG22}{article}{
      author={Amini, Omid},
      author={Gierczak, Lucas},
       title={Limit linear series: combinatorial theory},
        date={2026},
     journal={J. Comb. Algebra},
         url={https://doi.org/10.4171/jca/124},
}

\bib{amini2013reduced}{article}{
      author={Amini, Omid},
       title={Reduced divisors and embeddings of tropical curves},
        date={2013},
        ISSN={0002-9947},
     journal={Trans. Amer. Math. Soc.},
      volume={365},
      number={9},
       pages={4851\ndash 4880},
         url={https://doi.org/10.1090/S0002-9947-2013-05789-3},
      review={\MR{3066772}},
}

\bib{Ami14}{article}{
      author={Amini, Omid},
       title={Equidistribution of {W}eierstrass points on curves over
  non-{A}rchimedean fields},
        date={2014},
     journal={Preprint arXiv:1412.0926},
}

\bib{AP03}{article}{
      author={Ahlgren, Scott},
      author={Papanikolas, Matthew},
       title={Higher {W}eierstrass points on {$X_0(p)$}},
        date={2003},
        ISSN={0002-9947},
     journal={Trans. Amer. Math. Soc.},
      volume={355},
      number={4},
       pages={1521\ndash 1535},
         url={https://doi.org/10.1090/S0002-9947-02-03204-X},
      review={\MR{1946403}},
}

\bib{Arb74}{article}{
      author={Arbarello, Enrico},
       title={Weierstrass points and moduli of curves},
        date={1974},
        ISSN={0010-437X},
     journal={Compositio Math.},
      volume={29},
       pages={325\ndash 342},
      review={\MR{360601}},
}

\bib{Atk67}{article}{
      author={Atkin, Arthur O.~L.},
       title={Weierstrass points at cusps of ${\Gamma}_0(n)$},
        date={1967},
        ISSN={0003-486X},
     journal={Ann. of Math. (2)},
      volume={85},
       pages={42\ndash 45},
         url={https://doi.org/10.2307/1970524},
      review={\MR{218561}},
}

\bib{baker2008specialization}{article}{
      author={Baker, Matthew},
       title={Specialization of linear systems from curves to graphs},
        date={2008},
        ISSN={1937-0652},
     journal={Algebra \& Number Theory},
      volume={2},
      number={6},
       pages={613\ndash 653},
         url={https://doi.org/10.2140/ant.2008.2.613},
        note={With an appendix by Brian Conrad},
      review={\MR{2448666}},
}

\bib{BG95}{article}{
      author={Ballico, Edoardo},
      author={Gatto, Letterio},
       title={On the monodromy of {W}eierstrass points on {G}orenstein curves},
        date={1995},
        ISSN={0021-8693},
     journal={J. Algebra},
      volume={175},
      number={2},
       pages={633\ndash 643},
         url={https://doi.org/10.1006/jabr.1995.1205},
      review={\MR{1339660}},
}

\bib{BJ16}{incollection}{
      author={Baker, Matthew},
      author={Jensen, David},
       title={Degeneration of linear series from the tropical point of view and
  applications},
        date={2016},
   booktitle={Nonarchimedean and tropical geometry},
      series={Simons Symp.},
   publisher={Springer, [Cham]},
       pages={365\ndash 433},
      review={\MR{3702316}},
}

\bib{BL12}{article}{
      author={Brugall\'{e}, Erwan~A.},
      author={{L\'{o}pez de Medrano}, Lucia~M.},
       title={Inflection points of real and tropical plane curves},
        date={2012},
     journal={J. Singul.},
      volume={4},
       pages={74\ndash 103},
         url={https://doi.org/10.5427/jsing.2012.4e},
      review={\MR{3044488}},
}

\bib{baker2007riemann}{article}{
      author={Baker, Matthew},
      author={Norine, Serguei},
       title={Riemann--{R}och and {A}bel--{J}acobi theory on a finite graph},
        date={2007},
        ISSN={0001-8708},
     journal={Adv. Math.},
      volume={215},
      number={2},
       pages={766\ndash 788},
         url={https://doi.org/10.1016/j.aim.2007.04.012},
      review={\MR{2355607}},
}

\bib{BPR}{article}{
      author={Baker, Matthew},
      author={Payne, Sam},
      author={Rabinoff, Joseph},
       title={Nonarchimedean geometry, tropicalization, and metrics on curves},
        date={2016},
        ISSN={2313-1691},
     journal={Algebr. Geom.},
      volume={3},
      number={1},
       pages={63\ndash 105},
         url={https://doi.org/10.14231/AG-2016-004},
      review={\MR{3455421}},
}

\bib{BS13}{article}{
      author={Baker, Matthew},
      author={Shokrieh, Farbod},
       title={Chip-firing games, potential theory on graphs, and spanning
  trees},
        date={2013},
        ISSN={0097-3165},
     journal={J. Combin. Theory Ser. A},
      volume={120},
      number={1},
       pages={164\ndash 182},
         url={https://doi.org/10.1016/j.jcta.2012.07.011},
      review={\MR{2971705}},
}

\bib{BT20}{article}{
      author={Brezner, Uri},
      author={Temkin, Michael},
       title={Lifting problem for minimally wild covers of {B}erkovich curves},
        date={2020},
        ISSN={1056-3911},
     journal={J. Algebraic Geom.},
      volume={29},
      number={1},
       pages={123\ndash 166},
         url={https://doi.org/10.1090/jag/728},
      review={\MR{4028068}},
}

\bib{caporaso}{incollection}{
      author={Caporaso, Lucia},
       title={Algebraic and tropical curves: comparing their moduli spaces},
        date={2013},
   booktitle={Handbook of moduli. {V}ol. {I}},
      series={Adv. Lect. Math. (ALM)},
      volume={24},
   publisher={Int. Press, Somerville, MA},
       pages={119\ndash 160},
      review={\MR{3184163}},
}

\bib{CEG08}{article}{
      author={Cumino, Caterina},
      author={Esteves, Eduardo},
      author={Gatto, Letterio},
       title={Limits of special {W}eierstrass points},
        date={2008},
        ISSN={1687-3017},
     journal={Int. Math. Res. Pap. IMRP},
      number={2},
       pages={Art. ID rpn001, 65},
      review={\MR{2431731}},
}

\bib{CF91}{article}{
      author={Cukierman, Fernando},
      author={Fong, Lung-Ying},
       title={On higher {W}eierstrass points},
        date={1991},
        ISSN={0012-7094},
     journal={Duke Math. J.},
      volume={62},
      number={1},
       pages={179\ndash 203},
         url={https://doi.org/10.1215/S0012-7094-91-06208-3},
      review={\MR{1104328}},
}

\bib{Cuk89}{article}{
      author={Cukierman, Fernando},
       title={Families of {W}eierstrass points},
        date={1989},
        ISSN={0012-7094},
     journal={Duke Math. J.},
      volume={58},
      number={2},
       pages={317\ndash 346},
         url={https://doi.org/10.1215/S0012-7094-89-05815-8},
      review={\MR{1016424}},
}

\bib{Del08}{article}{
      author={Del~Centina, Andrea},
       title={Weierstrass points and their impact in the study of algebraic
  curves: a historical account from the ``{L}\"{u}ckensatz'' to the 1970s},
        date={2008},
        ISSN={0430-3202},
     journal={Ann. Univ. Ferrara Sez. VII Sci. Mat.},
      volume={54},
      number={1},
       pages={37\ndash 59},
         url={https://doi.org/10.1007/s11565-008-0037-1},
      review={\MR{2403373}},
}

\bib{CS94}{article}{
      author={de~Carvalho, C\'{\i}cero~Fernandes},
      author={St\"{o}hr, Karl-Otto},
       title={Higher order differentials and {W}eierstrass points on
  {G}orenstein curves},
        date={1994},
        ISSN={0025-2611},
     journal={Manuscripta Math.},
      volume={85},
      number={3-4},
       pages={361\ndash 380},
         url={https://doi.org/10.1007/BF02568204},
      review={\MR{1305748}},
}

\bib{Dia85}{article}{
      author={Diaz, Steven},
       title={Exceptional {W}eierstrass points and the divisor on moduli space
  that they define},
        date={1985},
        ISSN={0065-9266},
     journal={Mem. Amer. Math. Soc.},
      volume={56},
      number={327},
       pages={iv+69},
         url={https://doi.org/10.1090/memo/0327},
      review={\MR{791679}},
}

\bib{EH86}{article}{
      author={Eisenbud, David},
      author={Harris, Joe},
       title={Limit linear series: basic theory},
        date={1986},
        ISSN={0020-9910},
     journal={Invent. Math.},
      volume={85},
      number={2},
       pages={337\ndash 371},
         url={https://doi.org/10.1007/BF01389094},
      review={\MR{846932}},
}

\bib{EH87-WP}{article}{
      author={Eisenbud, David},
      author={Harris, Joe},
       title={Existence, decomposition, and limits of certain {W}eierstrass
  points},
        date={1987},
        ISSN={0020-9910},
     journal={Invent. Math.},
      volume={87},
      number={3},
       pages={495\ndash 515},
         url={https://doi.org/10.1007/BF01389240},
      review={\MR{874034}},
}

\bib{EH87-KD}{article}{
      author={Eisenbud, David},
      author={Harris, Joe},
       title={The {K}odaira dimension of the moduli space of curves of genus
  {$\geq 23$}},
        date={1987},
        ISSN={0020-9910},
     journal={Invent. Math.},
      volume={90},
      number={2},
       pages={359\ndash 387},
         url={https://doi.org/10.1007/BF01388710},
      review={\MR{910206}},
}

\bib{EM02}{article}{
      author={Esteves, Eduardo},
      author={Medeiros, Nivaldo},
       title={Limit canonical systems on curves with two components},
        date={2002},
        ISSN={0020-9910},
     journal={Invent. Math.},
      volume={149},
      number={2},
       pages={267\ndash 338},
         url={https://doi.org/10.1007/s002220200211},
      review={\MR{1918674}},
}

\bib{ES07}{article}{
      author={Esteves, Eduardo},
      author={Salehyan, Parham},
       title={Limit {W}eierstrass points on nodal reducible curves},
        date={2007},
        ISSN={0002-9947},
     journal={Trans. Amer. Math. Soc.},
      volume={359},
      number={10},
       pages={5035\ndash 5056},
         url={https://doi.org/10.1090/S0002-9947-07-04193-1},
      review={\MR{2320659}},
}

\bib{Esteves98}{incollection}{
      author={Esteves, Eduardo},
       title={Linear systems and ramification points on reducible nodal
  curves},
        date={1998},
      volume={14},
       pages={21\ndash 35},
        note={Algebra Meeting (Rio de Janeiro, 1996)},
      review={\MR{1662675}},
}

\bib{Gen21}{article}{
      author={Gendron, Quentin},
       title={Sur les n{\oe}uds de {W}eierstra{\ss}},
        date={2021},
     journal={Ann. H. Lebesgue},
      volume={4},
       pages={571\ndash 589},
         url={https://doi.org/10.5802/ahl.81},
      review={\MR{4275247}},
}

\bib{gathmann2008riemann}{article}{
      author={Gathmann, Andreas},
      author={Kerber, Michael},
       title={A {R}iemann--{R}och theorem in tropical geometry},
        date={2008},
        ISSN={0025-5874},
     journal={Math. Z.},
      volume={259},
      number={1},
       pages={217\ndash 230},
         url={https://doi.org/10.1007/s00209-007-0222-4},
      review={\MR{2377750}},
}

\bib{GL95}{article}{
      author={Garcia, Arnaldo},
      author={Lax, Robert~F.},
       title={On canonical ideals, intersection numbers, and {W}eierstrass
  points on {G}orenstein curves},
        date={1995},
        ISSN={0021-8693},
     journal={J. Algebra},
      volume={178},
      number={3},
       pages={807\ndash 832},
         url={https://doi.org/10.1006/jabr.1995.1379},
      review={\MR{1364344}},
}

\bib{haase2012linear}{article}{
      author={Haase, Christian},
      author={Musiker, Gregg},
      author={Yu, Josephine},
       title={Linear systems on tropical curves},
        date={2012},
        ISSN={0025-5874},
     journal={Math. Z.},
      volume={270},
      number={3-4},
       pages={1111\ndash 1140},
         url={https://doi.org/10.1007/s00209-011-0844-4},
      review={\MR{2892941}},
}

\bib{networkx}{inproceedings}{
      author={Hagberg, Aric~A.},
      author={Schult, Daniel~A.},
      author={Swart, Pieter~J.},
       title={Exploring network structure, dynamics, and function using
  {N}etwork{X}},
        date={2008},
   booktitle={Proceedings of the 7th python in science conference},
      editor={Varoquaux, Ga\"el},
      editor={Vaught, Travis},
      editor={Millman, Jarrod},
     address={Pasadena, CA USA},
       pages={11\ndash 15},
}

\bib{matplotlib}{article}{
      author={Hunter, John~D.},
       title={Matplotlib: A {2D} graphics environment},
        date={2007},
     journal={Computing in Science \& Engineering},
      volume={9},
      number={3},
       pages={90\ndash 95},
}

\bib{Hur92}{article}{
      author={Hurwitz, Adolf},
       title={{\"U}ber algebraische {G}ebilde mit eindeutigen
  {T}ransformationen in sich},
        date={1892},
     journal={Math. Ann.},
      volume={41},
      number={3},
       pages={403\ndash 442},
}

\bib{JP22}{article}{
      author={Jensen, David},
      author={Payne, Sam},
       title={Tropical linear series and tropical independence},
        date={2022},
     journal={Preprint arXiv:2209.15478},
}

\bib{KRZ16}{article}{
      author={Katz, Eric},
      author={Rabinoff, Joseph},
      author={Zureick-Brown, David},
       title={Uniform bounds for the number of rational points on curves of
  small {M}ordell--{W}eil rank},
        date={2016},
        ISSN={0012-7094},
     journal={Duke Math. J.},
      volume={165},
      number={16},
       pages={3189\ndash 3240},
         url={https://doi.org/10.1215/00127094-3673558},
      review={\MR{3566201}},
}

\bib{Lak81}{article}{
      author={Laksov, Dan},
       title={Weierstrass points on curves},
        date={1981},
      volume={87},
       pages={221\ndash 247},
      review={\MR{646822}},
}

\bib{Lax75}{article}{
      author={Lax, Robert~F.},
       title={Weierstrass points of the universal curve},
        date={1975},
        ISSN={0025-5831},
     journal={Math. Ann.},
      volume={216},
       pages={35\ndash 42},
         url={https://doi.org/10.1007/BF02547970},
      review={\MR{384809}},
}

\bib{Lax87-rational}{article}{
      author={Lax, Robert~F.},
       title={Weierstrass points on rational nodal curves},
        date={1987},
        ISSN={0017-0895},
     journal={Glasgow Math. J.},
      volume={29},
      number={1},
       pages={131\ndash 140},
         url={https://doi.org/10.1017/S0017089500006741},
      review={\MR{876157}},
}

\bib{LN64}{article}{
      author={Lehner, Joseph},
      author={Newman, Morris},
       title={Weierstrass points of ${\Gamma}_0(n)$},
        date={1964},
        ISSN={0003-486X},
     journal={Ann. of Math. (2)},
      volume={79},
       pages={360\ndash 368},
         url={https://doi.org/10.2307/1970550},
      review={\MR{161841}},
}

\bib{luo2011rank}{article}{
      author={Luo, Ye},
       title={Rank-determining sets of metric graphs},
        date={2011},
        ISSN={0097-3165},
     journal={J. Combin. Theory Ser. A},
      volume={118},
      number={6},
       pages={1775\ndash 1793},
         url={https://doi.org/10.1016/j.jcta.2011.03.002},
      review={\MR{2793609}},
}

\bib{LW90}{article}{
      author={Lax, Robert~F.},
      author={Widland, Carl},
       title={Weierstrass points on {G}orenstein curves},
        date={1990},
        ISSN={0030-8730},
     journal={Pacific J. Math.},
      volume={142},
      number={1},
       pages={197\ndash 208},
         url={http://projecteuclid.org/euclid.pjm/1102646469},
      review={\MR{1038736}},
}

\bib{Nee84}{article}{
      author={Neeman, Amnon},
       title={Weierstrass points in characteristic {$p$}},
        date={1984},
        ISSN={0020-9910},
     journal={Invent. Math.},
      volume={75},
      number={2},
       pages={359\ndash 376},
         url={https://doi.org/10.1007/BF01388569},
      review={\MR{732551}},
}

\bib{Ogg78}{article}{
      author={Ogg, Andrew~P.},
       title={On the {W}eierstrass points of ${X}_{0}(n)$},
        date={1978},
        ISSN={0019-2082},
     journal={Illinois J. Math.},
      volume={22},
      number={1},
       pages={31\ndash 35},
         url={http://projecteuclid.org/euclid.ijm/1256048830},
      review={\MR{463178}},
}

\bib{Richman18}{article}{
      author={Richman, David~Harry},
       title={The distribution of {W}eierstrass points on a tropical curve},
        date={2024},
        ISSN={1022-1824,1420-9020},
     journal={Selecta Math. (N.S.)},
      volume={30},
      number={2},
       pages={Paper No. 31, 48},
         url={https://doi.org/10.1007/s00029-024-00919-5},
      review={\MR{4712053}},
}

\bib{Tem16}{incollection}{
      author={Temkin, Michael},
       title={Metrization of differential pluriforms on {B}erkovich analytic
  spaces},
        date={2016},
   booktitle={Nonarchimedean and tropical geometry},
      series={Simons Symp.},
   publisher={Springer, [Cham]},
       pages={195\ndash 285},
      review={\MR{3702313}},
}

\bib{TT22}{article}{
      author={Temkin, Michael},
      author={Tyomkin, Ilya},
       title={Reduction and lifting problem for differential forms on
  {B}erkovich curves},
        date={2022},
        ISSN={0001-8708},
     journal={Adv. Math.},
      volume={397},
       pages={Paper No. 108208, 23},
         url={https://doi.org/10.1016/j.aim.2022.108208},
      review={\MR{4366859}},
}

\bib{Weierstrass67}{book}{
      author={Weierstrass, Karl},
       title={Mathematische werke, 7 vols (1894--1927)},
   publisher={Olms, Hildesheim},
        date={1967},
}

\end{biblist}
\end{bibdiv}
